\documentclass[3pt]{elsarticle}
\pdfoutput=1
\usepackage{amsthm}
\usepackage[intlimits]{amsmath}
\usepackage{amssymb}
\usepackage{amsthm}
\usepackage{amsfonts}
\usepackage{bm}
\usepackage{mathdots}
\usepackage[normalem]{ulem}
\usepackage[utf8]{inputenc} % for UTF-8 encoding
\usepackage[T1]{fontenc}
\usepackage{times}
\usepackage{mathptmx}  % times roman, including math
\DeclareSymbolFont{cmletters}{OML}{cmr}{m}{n}
\DeclareMathAlphabet{\mathcal}{OMS}{cmsy}{m}{n} % recover 'mathcal'

\journal{Computer and Geotechnics}
\ifdefined\copyReadyOn
\usepackage{ecrc}
\volume{00}
\firstpage{1}
\runauth{G. Beer, Ch. Duenser, V. Mallardo}
\journalname{\journal}
\jid{CG}
\jnltitlelogo{CG}
\CopyrightLine{2015}{Published by Elsevier Ltd.}
\fi

\usepackage{tikz}
\usepackage{tikz-3dplot}
\usepackage{pgfplots}
\usepackage{pgfplotstable}
\usepackage{tikzscale}
\usepackage{standalone} % for standalone compilation of tikz graphics
\usepackage{hf-tikz}

\usepackage{forloop}
\usepackage{arrayjobx}

\usepackage{algorithmic}
\usepackage{algorithm}

\usepackage{booktabs} % for tables: toprule, midrule, bottomrule
\usepackage{calc}
\usepackage{xcolor,pdfcomment,soul} % add pdf comments
\usepackage{xspace} % horizontal space after makro words
\usepackage{multirow}

\usepackage{scrlayer-scrpage}
\usepackage{setspace}

\usepackage[labelformat=simple]{subcaption}   	%% 1(a)
\usepackage{pdfpages} % to include external pdf (i.e. eidesstattliche erlärung)

\usepackage{tabto} % for tabbing within itemize

\usepackage{makeidx}
\makeindex

\usepackage{colortbl}

\usepackage{enumerate}

\usepackage{overpic}

\usepackage{listings} %% enviroment for code formating 

\usepackage[intlimits]{amsmath}
\usepackage{amssymb} 
\usepackage{amsthm}
\usepackage{amsfonts,mathrsfs}
\usepackage{bm}
\usepackage{mathtools} % for \coloneqq and \shortintertext

\usepackage{siunitx} %% units

\usepackage{xfrac}
\usepackage{listings}
\usepackage{color}
\definecolor{dkgreen}{rgb}{0,0.6,0}
\definecolor{gray}{rgb}{0.5,0.5,0.5}
\definecolor{mauve}{rgb}{0.58,0,0.82}

\usepackage[utf8]{inputenc} % for UTF-8 encoding
\usepackage[german,english]{babel}  
\usepackage[T1]{fontenc} % größerer Zeichenvorat z.B. für Umlaute
\usepackage{lmodern}

\usepackage{geometry}
\usepackage{overpic}
\graphicspath{{pics/}{tmp/}}

\newboolean{inGray}
\setboolean{inGray}{false}
\newcommand{\addtoindex}[2][]{% key=#1,text=#2
    \ifthenelse { \equal{#1}{} }
    {#2\index{#2}\xspace}%
    {#2\index{#1}\xspace}%  
}

\newcommand{\myVec}[1]{\mathbf{#1}}
\newcommand{\myVecGreek}[1]{\boldsymbol{#1}}

\newcommand\tens[2]{\mathsf{#1}} % 
\newcommand{\NURBS}{R} 			% symbol for NURBS
\newcommand\uu{\xi} 			% Local NURBS coordinates
\newcommand\vv{\eta} 			% Local NURBS coordinates
\newcommand\abscissa{\tilde{\uu}}			% abscissa
\newcommand{\domain}{\Omega}

\newcommand{\boundary}{\Gamma}

\newcommand\primary{u} 			% primary variable - displacements in elasticity
\newcommand\dual{t} 			% dual variable - tractions in elasticity
\providecommand\url[1]{\emph{#1}}

\newcommand\fund[1]{\tens{#1}{2}}

\newcommand\pt[1]{\boldsymbol{#1}}
\newcommand\sourcept{\tilde{\pt{x}}}
\newcommand\fieldpt{\hat{\pt{x}}}
\newcommand\fieldptv{\check{\pt{x}}}

\newenvironment{mytable}[4]%
{
  \begin{table}[#1]%[htb!]
    \def\mycap{#2}
    \def\mylabel{#3}
    \centering
    \begin{tabular}{#4}
      \toprule
}
{
  \bottomrule
  \end{tabular}
  \caption{\mycap}
  \label{\mylabel}
  \end{table}
}
\newcommand{\mytableheader}[1]{
       #1 \\ \midrule
}

\newenvironment{myalgorithm}[2]%
{
  \begin{algorithm}
    \caption{#1}
    \label{#2}
    \begin{algorithmic}[1]
}   
{
\end{algorithmic}
\end{algorithm}
}

\newtheoremstyle{myremark}% ⟨name⟩ 
{3pt}% Space above
{3pt}% Space below
{}% Body font
{}% Indent amount
{\itshape}% Theorem head font
{:}% Punctuation after theorem head
{.5em}% Space after theorem head
{}% Theorem head spec (can be left empty, meaning ‘normal’)
\theoremstyle{myremark}
\newtheorem*{remark}{Remark}

\newcommand{\myeqref}[1]{equation~(\ref{#1})}	
\newcommand{\myfigref}[1]{Figure~\ref{#1}}

\newcounter{footnoteNumber} % initialization of counter with 0

\DeclareCaptionFormat{listing}{#1#2#3}
\newcommand{\myalignatsinglelabel}[3]%
{%
    \begin{equation}
        \label{#2}
        \begin{alignedat}{#1} 
            #3
        \end{alignedat}
    \end{equation}
}

\newcommand{\myalignat}[2]%
{%
  \begin{alignat}{#1} 
    #2
   \end{alignat}
}

\usepackage{pgf,pgffor,pgfmath}
\pgfplotsset{compat=1.10}

\usetikzlibrary{external}
\usetikzlibrary{positioning,calc,intersections,arrows,3d}
\usetikzlibrary{shapes,snakes,fit,spy}
\usepgfplotslibrary{patchplots}
\usetikzlibrary{matrix}

\usepgfplotslibrary{groupplots}
\usepgfplotslibrary{fillbetween}

\def\pgfplotfontsizetitle{\small}
\def\pgfplotfontsizelegend{\small}
\def\pgfplotfontsize{\small}
\def\pgfplotfontsizetiny{\scriptsize}

\def\tikzfontsizetiny{\scriptsize}

\pgfplotsset{
  mystyle/.style ={%
    grid = major,
    every tick label/.append style={font=\pgfplotfontsizetiny},
    every axis label/.append style={font=\pgfplotfontsize},
    legend style={font=\pgfplotfontsizelegend},
    label style={font=\pgfplotfontsize},
    title style={font=\pgfplotfontsizetitle},
    /pgf/number format/set thousands separator = {}, % 1,000 -> 1000
  }
}% 

\makeatletter
\pgfplotsset{
    myIgnoreRowModulo2/.style args={#1}{%
        /pgfplots/x filter/.code={%
        \let\xValue\pgfmathresult % store x value of current table index
        \pgfmathparse{int(mod(int(\coordindex),int(2))} \pgfmathresult % compute modulo
        \ifnum#1=\pgfmathresult
            \def\pgfmathresult{} % ignore table entry
        \else
            \pgfmathparse{\xValue} \pgfmathresult % restore x value to \pgfmathresult macro
        \fi
        }
    } 
}
\makeatother

\colorlet{drawblue}      {blue!80!white}
\colorlet{drawred}       {red!80!white}
\colorlet{drawgray}      {gray}
\definecolor{drawgreen}  {RGB}{44,162,95}
\colorlet{drawpurple}    {purple}
\colorlet{draworange}    {orange}
\colorlet{drawlime}      {lime!80!black}
\colorlet{drawartichoke} {yellow!60!black}
\colorlet{TUGgray}{black!15}
\definecolor{TUGred}{RGB}{247,1,70}
\definecolor{IFBblue}{RGB}{51,112,169}

\definecolor{basisColor1}{RGB}{59,76,192}
\definecolor{basisColor2}{RGB}{87,117,225}
\definecolor{basisColor3}{RGB}{119,154,247}
\definecolor{basisColor4}{RGB}{152,185,255}
\definecolor{basisColor5}{RGB}{184,208,249}
\definecolor{basisColor6}{RGB}{195,209,230}	% changed for printing
\definecolor{basisColor7}{RGB}{247,200,190}	% changed for printing
\definecolor{basisColor8}{RGB}{247,187,160}
\definecolor{basisColor9}{RGB}{244,154,123}
\definecolor{basisColor10}{RGB}{229,112,88}
\definecolor{basisColor11}{RGB}{203,62,56}
\definecolor{basisColor12}{RGB}{180,4,38}

\definecolor{basisColor8sw}{RGB} {189,189,189}
\definecolor{basisColor9sw}{RGB} {150,150,150}
\definecolor{basisColor10sw}{RGB}{115,115,115}
\definecolor{basisColor11sw}{RGB}{91,91,91}
\definecolor{basisColor12sw}{RGB}{37,37,37}

\colorlet{myblue}    {blue}
\colorlet{myred}     {red}
\colorlet{mygreen}   {drawgreen}
\colorlet{mypurple}  {purple}
\colorlet{myorange}  {orange}

\ifthenelse{\boolean{inGray}}{ 
\pgfplotscreateplotcyclelist{mycyclelist}{
  basisColor12sw,semithick, every mark/.append style={solid, fill=.},mark=*           \\%
  basisColor11sw,semithick, every mark/.append style={solid, fill=.},mark=square*     \\%
  basisColor10sw,semithick, every mark/.append style={solid, fill=.},mark=triangle*   \\%
  basisColor9sw , semithick,every mark/.append style={solid, fill=.},mark=diamond*    \\%
  basisColor8sw , semithick,every mark/.append style={solid, fill=.},mark=otimes    \\%
}
}{

\pgfplotscreateplotcyclelist{mycyclelist}{
  basisColor12,semithick,   every mark/.append style={solid, fill=.},mark=*           \\%
  basisColor11,semithick,   every mark/.append style={solid, fill=.},mark=square*     \\%
  basisColor10,semithick,   every mark/.append style={solid, fill=.},mark=triangle*   \\%
  basisColor9, semithick,   every mark/.append style={solid, fill=.},mark=diamond*    \\%
  basisColor8, semithick,   every mark/.append style={solid, fill=.},mark=otimes    \\%
}
}

\ifthenelse{\boolean{inGray}}{ 
\pgfplotscreateplotcyclelist{mycyclelistcompare}{
  basisColor12sw, semithick, loosely dashdotdotted     \\%
  basisColor11sw, semithick, loosely dashed      		\\%
  basisColor10sw, semithick, dashed \\%loosely dashdotted        \\%
  basisColor9sw,  semithick, solid               	\\%
  basisColor12sw, only marks,   every mark/.append style={solid, fill=.},mark=otimes    \\%
  basisColor11sw, only marks,   every mark/.append style={solid, fill=.},mark=triangle*   \\%
  basisColor10sw, only marks,   every mark/.append style={solid, fill=.},mark=square*     \\%
  basisColor9sw,  only marks,   every mark/.append style={solid, fill=.},mark=*          \\%
}

}{
\pgfplotscreateplotcyclelist{mycyclelistcompare}{
  basisColor12, semithick, loosely dashdotdotted     \\%
  basisColor11, semithick, loosely dashed      		\\%
  basisColor10, semithick, dashed \\%loosely dashdotted        \\%
  basisColor9,  semithick, solid               	\\%
  basisColor12, only marks,   every mark/.append style={solid, fill=.},mark=otimes    \\%
  basisColor11, only marks,   every mark/.append style={solid, fill=.},mark=triangle*   \\%
  basisColor10, only marks,   every mark/.append style={solid, fill=.},mark=square*     \\%
  basisColor9,  only marks,   every mark/.append style={solid, fill=.},mark=*          \\%
}
}

\ifthenelse{\boolean{inGray}}{ 
\tikzset{mycyclelistcompareReferenceA/.style={basisColor12sw,solid}}
\tikzset{mycyclelistcompareTestA/.style={basisColor12sw,only marks,mark=otimes}}
}{
\tikzset{mycyclelistcompareReferenceA/.style={basisColor12,solid}}
\tikzset{mycyclelistcompareTestA/.style={basisColor12,only marks,mark=otimes}}
}
\tikzset{helpline/.style={thin,dashed}}
\tikzset{labelline/.style={thin}}
\tikzset{referencePath/.style={dotted,very thick}}
\tikzset{helparrow/.style={thin,arrows={-latex}}}
\tikzset{axis/.style={thin,arrows={->}}}
\tikzset{force/.style={thick,arrows={->}}}
\tikzset{forceInverse/.style={thick,arrows={<-}}}
\tikzset{Gamma/.style={ultra thick}}
\tikzset{controlPoly/.style={draw=black}}
\tikzset{GammaFill/.style={fill=lightgray,fill opacity=0.5}}
\tikzset{colorDiri/.style={drawgreen}}
\tikzset{GammaFillDiri/.style={fill=drawgreen,fill opacity=0.5}}
\tikzset{initialgrid/.style={thin,gray}}
\tikzset{addgridline/.style={dashed,gray}}
\tikzset{trimmingcurve/.style={thick}}
\tikzset{boundingbox/.style={thick, dotted}}
\tikzset{parameterSpace/.style={ }}
\tikzset{basisfunction/.style={very thick,smooth}}
\tikzset{bspline/.style={very thick,smooth,red}}
\tikzset{intersectioncurve/.style={dashed,thick}}
\tikzset{integrationRegionEdge/.style={dashed}}
\tikzset{pointer/.style={arrows={-latex}}}

\tikzstyle{anode}= [circle, inner sep=1.3pt, draw, fill=black]
\tikzstyle{gausspoint}=[shape=circle,draw=black,fill=black,inner sep=1.1pt]
\tikzstyle{controlPoint}=[shape=circle,draw=black,fill=black,thin,inner sep=0pt,minimum size=1.5mm]
\tikzstyle{abscissaPoint}=[shape=circle,draw=black,fill=white,thin,inner sep=0pt,minimum size=1.5mm]
\tikzstyle{anchorPoint}=[shape=circle,draw=black,fill=black,thin,inner sep=0pt,minimum size=1.5mm]
\tikzstyle{anchorPointDeg}=[shape=circle,draw=black,fill=TUGred,thin,inner sep=0pt,minimum size=1.5mm]
\tikzstyle{anchorPointDegD}=[shape=cross out,thick,draw=black,inner sep=0pt,minimum size=1.5mm]
\tikzstyle{trimmingIntersectionPoint}=[shape=cross out,thick,draw=black,inner sep=0pt,minimum size=1.5mm]

\tikzset{%
  highlight/.style={rectangle,rounded corners,fill=red!60,draw,fill opacity=0.125,thick,inner sep=0pt}
}

\usepackage{colortbl}

\def\trianglecolor{black}
\newcommand{\upperSlopeTriangle}[4] 	% input: #1 slope #2 vertical shift  #3 min x-value #4 max x-value
{
	\def\trianglecolor{black}
	\addplot[forget plot, domain=#3:#4,color=\trianglecolor,samples=2]{  #2 / (x^#1) } node (A1) [pos=1] {}; 
	\addplot[forget plot, domain=#3:#4,color=\trianglecolor,samples=2]{   #2  / (#3^#1)} node (A2) [pos=1] {} node [anchor=south,pos=0.5,black] {\tikzfontsizetiny $1$};
	\draw[color=\trianglecolor] (A1.center) -- (A2.center) node [anchor=west,pos=0.5,black] {\tikzfontsizetiny #1};
}

\newcommand{\lowerSlopeTriangle}[4] 	% input: #1 slope #2 vertical shift  #3 min x-value #4 max x-value
{
	\def\trianglecolor{black}
	\addplot[forget plot, domain=#3:#4,color=\trianglecolor,samples=2]{  #2 / (x^#1) } node (A1) [pos=0] {}; 
	\addplot[forget plot, domain=#3:#4,color=\trianglecolor,samples=2]{   #2  / (#4^#1)} node (A2) [pos=0] {} node [anchor=north,pos=0.5,black] {\tikzfontsizetiny $1$};
	\draw[color=\trianglecolor] (A1.center) -- (A2.center) node [anchor=east,pos=0.5,black] {\tikzfontsizetiny #1};
}

\newcommand{\myaddgraphic}[5]
{
 \node[anchor=south west,inner sep=0] (image) {\phantom{\includegraphics[#2]{#1}}};
  \begin{scope}[x={(image.south east)},y={(image.north west)}]
      
      \begin{scope}
          
          #5
          
          \node[anchor=south west,inner sep=0] {\includegraphics[#2]{#1}};
      \end{scope} 
      
      #4
      
      \pgfmathparse{int(#3)} \let\gridIndicator\pgfmathresult
      \ifthenelse{ \gridIndicator = 1 }
      {
          \draw[help lines,xstep=.1,ystep=.1] (0,0) grid (1.001,1.001);
          \foreach \x in {1,...,9} { \node [anchor=north] at (\x/10,0) {\x};}
          \foreach \y in {1,...,9} { \node [anchor=east] at (0,\y/10) {\y};}
      }{}
      
  \end{scope}    
}

\tikzstyle{reverseclip}=[insert path={(current page.north east) --
    (current page.south east) --
    (current page.south west) --
    (current page.north west) --
    (current page.north east)}
]

\newcounter{itR}
\newcommand{\bsplinevalue}[5] % input: #1 knot vector #2 order,  #3 intrinsic coordinate #4 span index #5 output name
{                    				
	\newarray\vKnots
	\newarray\vN
	\newarray\vNumeratorL
	\newarray\vNumeratorR
	\newarray\vSave

	\readarray{vKnots}{#1}
	\readarray{vN}{1}
	\readarray{vSave}{0}
	\readarray{vNumeratorL}{0}
	\readarray{vNumeratorR}{0}

	\foreach \j in {1,...,#2}
	{        
		\pgfmathparse{ int(#4+\j+1) } \checkvKnots(\pgfmathresult)
		\pgfmathsetmacro{\numR}{\cachedata-#3}
		
		\pgfmathparse{ int(#4-\j+2) } \checkvKnots(\pgfmathresult)
		\pgfmathsetmacro{\numL}{#3-\cachedata} 					

		\expandarrayelementtrue
		\pgfmathparse{ int(\j+1) }
		\vNumeratorL(\pgfmathresult)={\numL}
		\vNumeratorR(\pgfmathresult)={\numR}
		
		\forloop[1]{itR}{0}{\value{itR} < \j }
		{
			\pgfmathparse{ int(\theitR+1) } \let\tS\pgfmathresult  	
			\checkvSave(\tS)							
			\pgfmathsetmacro{\save}{\cachedata}  
		
			\pgfmathparse{int(\j-\theitR+1)} \checkvNumeratorL(\pgfmathresult)
			\pgfmathsetmacro{\tmpL}{\cachedata} 	
		    
			\pgfmathparse{ int(\theitR+2) } \checkvNumeratorR(\pgfmathresult)
			\pgfmathsetmacro{\tmpR}{\cachedata} 	       
	
			\pgfmathparse{ int(\theitR+1) } \checkvN(\pgfmathresult)
			\pgfmathparse{ \cachedata / (\tmpL+\tmpR) } \let\tmp\pgfmathresult
			
			\pgfmathparse{ \save + \tmpR * \tmp } 
			\vN(\tS)={\pgfmathresult}

			\pgfmathparse{ \tmpL * \tmp } \let\tmpsave\pgfmathresult
			\pgfmathparse{ int(\theitR+2) } \let\tS\pgfmathresult      
			\vSave(\tS)={\tmpsave}
		}
		
		\pgfmathparse{ int(\j+1) } \checkvSave(\pgfmathresult )
		\vN(\tS)={\cachedata}
	}

	\pgfmathparse{int( #2+1) } \let\lastIndex\pgfmathresult
	\checkvN(1) \pgfmathsetmacro{\first}{\cachedata}  
	\foreach \i [remember=\a as \lasta (initially \first)] in {2,...,\lastIndex}
	{
		\checkvN(\i) \def\a{\lasta,\cachedata} 
		\ifthenelse{\i=\lastIndex}{ \xdef#5{\a} }{}
	}
    
	\foreach \i in {1,...,\lastIndex}
	{
	    \clrarray{vNumeratorR}(\i)
	    \clrarray{vNumeratorL}(\i)
	    \clrarray{vN}(\i)
	    \clrarray{vSave}(\i)
	}
	\delarray\vN
	\delarray\vNumeratorL
	\delarray\vNumeratorR
	\delarray\vSave
}

\newcounter{countvalues}
\newcounter{getBasis}
\newcommand{\bsplinebasis}[4] % input: #1 knot vector, #2 order,  #3 id of basis function, #4 output
{						% example: \bsplinebasis{0&0&0&1&2&2&2}{2}{0}{\output}
    \newarray\vKnots 	% remark index start from 1
    \readarray{vKnots}{#1}

    \pgfmathparse{#3}  \let\i\pgfmathresult
    \pgfmathparse{#2}  \let\p\pgfmathresult

    \setcounter{countvalues}{0}
    \setcounter{getBasis}{\p} %\addtocounter{getBasis}{1}
    \pgfmathparse{int(\i+\p)}
    \foreach \knotspan in {\i,...,\pgfmathresult}
    {
	\pgfmathparse{ int(\knotspan+1+1) } \checkvKnots(\pgfmathresult)
	\pgfmathsetmacro{\tmpR}{\cachedata} 						%R: \tmpR
	
	\pgfmathparse{ int(\knotspan+1) } \checkvKnots(\pgfmathresult) 
	\pgfmathsetmacro{\tmpL}{\cachedata} 						%L: \tmpL

	\pgfmathparse{ \tmpR - \tmpL } \let\spansize\pgfmathresult  			%S: \spansize 
   
        \pgfmathparse{ \spansize > 0.0 } \let\bNonZero\pgfmathresult
        \ifthenelse{ \bNonZero = 1 }
        {
            \foreach \percentU in {0,10,...,100}
            {
		\pgfmathparse{\tmpL+\spansize*\percentU/100} \let\u\pgfmathresult

		\bsplinevalue{#1}{\p}{\u}{\knotspan}{\Basis} 
		\def\basisfuncarray{{\Basis}} 				% remark index start from 0
		
		\pgfmathparse{\basisfuncarray[\thegetBasis]} %\pgfmathresult
		
		\ifthenelse{\thecountvalues=0}
		{ 
			\xdef\nodeB{"\u,\pgfmathresult"}
		}{
			\xdef\nodeB{\nodeB,"\u,\pgfmathresult"}  
		}
		\addtocounter{countvalues}{1}
            }
        }{}
        \addtocounter{getBasis}{-1}
    }
    
	\xdef#4{\nodeB}

	\delarray\vKnots
}

\newcommand{\bsplinebasisspan}[6] % input: 	#1 knot vector, #2 order,  #3 segment id of basis function within (#4), 
{							% 		#4 id of knot span,  #5 id of plotted span, #6 output
    \newarray\vKnots 	% remark index start from 1
    \readarray{vKnots}{#1}

    \pgfmathparse{#5}  \let\plotknotspan\pgfmathresult
    \pgfmathparse{#4}  \let\splineknotspan\pgfmathresult
    \pgfmathparse{#3}  \let\i\pgfmathresult
    \pgfmathparse{#2}  \let\p\pgfmathresult

    \setcounter{countvalues}{0}
    \setcounter{getBasis}{\i} 
    \foreach \knotspan in {\plotknotspan}
    {
	\pgfmathparse{ int(\knotspan+1+1) } \checkvKnots(\pgfmathresult)
	\pgfmathsetmacro{\tmpR}{\cachedata} 						%R: \tmpR
	
	\pgfmathparse{ int(\knotspan+1) } \checkvKnots(\pgfmathresult) 
	\pgfmathsetmacro{\tmpL}{\cachedata} 						%L: \tmpL

	\pgfmathparse{ \tmpR - \tmpL } \let\spansize\pgfmathresult  			%S: \spansize 
   
        \pgfmathparse{ \spansize > 0.0 } \let\bNonZero\pgfmathresult
        \ifthenelse{ \bNonZero = 1 }
        {
            \foreach \percentU in {0,10,...,100}
            {
		\pgfmathparse{\tmpL+\spansize*\percentU/100} \let\u\pgfmathresult

		\bsplinevalue{#1}{\p}{\u}{\splineknotspan}{\Basis} 
		\def\basisfuncarray{{\Basis}} 				% remark index start from 0
		
		\pgfmathparse{\basisfuncarray[\thegetBasis]} %\pgfmathresult
		
		\ifthenelse{\thecountvalues=0}
		{ 
			\xdef\nodeB{"\u,\pgfmathresult"}
		}{
			\xdef\nodeB{\nodeB,"\u,\pgfmathresult"}  
		}
		\addtocounter{countvalues}{1}
            }
        }{}
        \addtocounter{getBasis}{-1}
    }
    
	\xdef#6{\nodeB}

	\delarray\vKnots
}

\newcommand{\plotbsplinebasis}[4] 	% input: #1 knot vector, #2 order,  #3 id of basis function, #4 plot settings
{							% example: \plotbsplinebasis{0&0&0&1&2&2&2}{2}{0}{green,smooth}
	\bsplinebasis{#1}{#2}{#3}{\nodeOut}
	\def\nodearray{{\nodeOut}}

	\xdef\name{ }
	\addtocounter{countvalues}{-1}
	\foreach \i in {0,...,\thecountvalues}
	{
		\pgfmathparse{\nodearray[\i]}
		\coordinate (point\i) at (\pgfmathresult);	  
		\xdef\name{ \name (point\i)  }
	}
	
	\draw[#4] plot coordinates{ \name };
	
	\xdef\name{ }
}

\newcommand{\plotbsplinesegment}[6] 	% input: 	#1 knot vector, #2 order,  #3 segment id of basis function within (#4),  
{								%		#4 id of knot span, #5 id of plotted span, #6 plot settings
	\bsplinebasisspan{#1}{#2}{#3}{#4}{#5}{\nodeOut}
	\def\nodearray{{\nodeOut}}

	\xdef\name{ }
	\addtocounter{countvalues}{-1}
	\foreach \i in {0,...,\thecountvalues}
	{
		\pgfmathparse{\nodearray[\i]}
		\coordinate (point\i) at (\pgfmathresult);	  
		\xdef\name{ \name (point\i)  }
	}
	
	\draw[#6] plot coordinates{ \name };
	
	\xdef\name{ }
}

\newcommand{\plotbsplineaccumulated}[5] 	% input: 	#1 knot vector, #2 order,  #3 subdivision coefficient,  
{						%		#4 id of knot span,  #5 plot settings
    \newarray\vKnots 	% remark index start from 1
    \readarray{vKnots}{#1}
    \newarray\vSubCoef 	% remark index start from 1
    \readarray{vSubCoef}{#3}
    
    \pgfmathparse{#4}  \let\plotknotspan\pgfmathresult
    \pgfmathparse{#4}  \let\splineknotspan\pgfmathresult
    \pgfmathparse{#2}  \let\p\pgfmathresult
    
    \setcounter{countvalues}{0}
    \pgfmathparse{int( \p+1) } \let\lastIndex\pgfmathresult
    \foreach \knotspan in {\plotknotspan}
    {
        \pgfmathparse{ int(\knotspan+1+1) } \checkvKnots(\pgfmathresult)
        \pgfmathsetmacro{\tmpR}{\cachedata} 						%R: \tmpR
        
        \pgfmathparse{ int(\knotspan+1) } \checkvKnots(\pgfmathresult) 
        \pgfmathsetmacro{\tmpL}{\cachedata} 						%L: \tmpL
        
        \pgfmathparse{ \tmpR - \tmpL } \let\spansize\pgfmathresult  			%S: \spansize 
        
        \pgfmathparse{ \spansize > 0.0 } \let\bNonZero\pgfmathresult
        \ifthenelse{ \bNonZero = 1 }
        {
            \foreach \percentU in {0,10,...,100}
            {
                \pgfmathparse{\tmpL+\spansize*\percentU/100} \let\u\pgfmathresult
                
                \bsplinevalue{#1}{\p}{\u}{\splineknotspan}{\Basis} 
                \def\basisfuncarray{{\Basis}} 				% remark index start from 0
                
                \setcounter{getBasis}{0} 
                \pgfmathparse{\basisfuncarray[\thegetBasis]} %\pgfmathresult
                \let\basisValue\pgfmathresult
                
                \checkvSubCoef(1) \pgfmathsetmacro{\coef}{\cachedata}  
                \pgfmathparse{ \basisValue * \coef } \let\first\pgfmathresult
                
                \xdef\lastx{\first}
                \foreach \i in {2,...,\lastIndex}
                { 
                    \addtocounter{getBasis}{1}       
                    \pgfmathparse{\basisfuncarray[\thegetBasis]}
                    \let\basisValue\pgfmathresult
                    
                    \checkvSubCoef(\i) \pgfmathsetmacro{\coef}{\cachedata}  
                    \pgfmathparse{ \lastx + \basisValue * \coef } \let\sum\pgfmathresult
                    
                    \xdef\lastx{\sum}
                    
                    \ifthenelse{\i=\lastIndex}
                    {
                        \ifthenelse{\thecountvalues=0}
                        { 
                            \xdef\nodeBB{"\u,\sum"}
                        }{
                            \xdef\nodeBB{\nodeBB,"\u,\sum"}  
                        }
                        \addtocounter{countvalues}{1}
                    }{}
                    
                }
            }
        }{}
    }
    
    \delarray\vKnots
    \delarray\vSubCoef
    
    \def\nodearray{{\nodeBB}}
    
    \xdef\name{ }
    \addtocounter{countvalues}{-1}
    \foreach \i in {0,...,\thecountvalues}
    {
        \pgfmathparse{\nodearray[\i]}
        \coordinate (point\i) at (\pgfmathresult);                
        \xdef\name{ \name (point\i)  }
    }
    
    \draw[#5] plot coordinates{ \name };
    
    \xdef\name{ }
}

\begin{document}
    
\title{Immersed isogeometric Boundary Elements: A user friendly method for the 3-D elasto-plastic simulation of underground excavations.}
\begin{frontmatter}

\author[ifbaddr]{Gernot Beer\corref{cor1}}
\author[ifbaddr]{Christian Duenser}

\address[ifbaddr]{Institute of Structural Analysis, Graz University
  of Technology, Lessingstraße 25/II, 8010 Graz, Austria}

% \address[unife]{Department of Architecture, University of Ferrara, Via Quartieri 8, 44121 Ferrara, Italy}

\cortext[cor1]{Corresponding author.
  mail: \url{gernot.beer@tugraz.at}, web: \url{www.ifb.tugraz.at}}

\begin{abstract}
The immersed isogeometric Boundary Element Method is presented and applied to the simulation of underground excavations. Nonuniform rational B-splines (NURBS) are used for the accurate definition of complex geometries with few parameters.
Immersed technology is applied to automatically generate cell meshes. This allows heterogeneous and anisotropic ground conditions as well as nonlinear material behaviour to be considered. On a practical example the user friendliness and accuracy is demonstrated.

\end{abstract}

\end{frontmatter}

\section{Introduction}
The Boundary Element Method (BEM) is an attractive simulation method for underground excavations, because it allows the consideration of semi-infinite or infinite domains without mesh truncation. For homogeneous, elastic domains a discretisation of the volume can be avoided and only the boundary surface needs to be defined. While the BEM has been around as long as it's big sister the Finite Element Method (FEM) its application to geomechanics problems has been very limited. Only few companies offer BEM software and then only for linear elastic domains.
This is despite the fact that early in the development of the BEM, publications appeared \cite{Telles1981,Banerjee1981}, that demonstrated that the method could be extended to deal with inelastic problems. However, a volume discretisation was required. This seems to destroy the very nature of the BEM, namely the boundary only discretisation. However, this is not entirely  true, since volume discretisation is only required for the near field and not for the whole domain. The main advantage, that the BEM can deal with infinite domains, remains.
Examples of BEM simulations for 2-D elasto-plastic problems in geomechanics were published in  \cite{Venturini,Brebbia1983}.  Examples in 3-D are rare in the literature and are mainly in other fields of applications such a mechanical engineering (see for example \cite{GaoDavies2000}). The majority of the publications use the modified Newton-Raphson method but it was demonstrated in \cite{GaoDavies2000} that it is also possible to use a Newton-Raphson approach, by computing a sort of tangent stiffness.

 In any case, the fact that a volume discretisation is required to deal with elasto-plastic and heterogeneous problems, was a major stumbling block in the application of the method to geomechanics.  In addition, the simulation of ground support such as cable dowels, rock anchors or shotcrete was not addressed.
In the past decades the main emphasis of the research group at the TU Graz was aimed at making the BEM more user friendly and suitable for geomechanics. Isogeometric analysis, first introduced  by Hughes \cite{Hughes2005a}, who proposed that NURBS should be used instead of Lagrange polynomials, showed an elegant way in which geometries could be defined very accurately with few parameters, basically avoiding mesh generation. Early publications of the isogeometric BEM (IGABEM) \cite{scott2013isogeometric,simpson2013isogeometric} showed the advantages of this approach.
In addition the use of NURBS for the approximation of the unknowns resulted in elegant refinement procedures. In \cite{marussig2016a} it was pointed out that the definition of the geometry could be uncoupled form the approximation of the unknown and this was verified in \cite{doi:10.1002/nme.5778}.
The early emphasis of our work was focussed on the elimination of the requirement for a volume mesh generation for nonlinear problems. A proposed solution was to define the volume by 3-D NURBS  \cite{Beer2016,Beer17}. The simulation of reinforcement bars and rock bolts was also addressed in \cite{BEER2020113409}, using linear inclusions.
While the definition of the volume with 3-D NURBS works well for simple geometries it is not applicable to more complicated geometries. For these problems we have to revert to cell discretisation, which brings us back to the problem of user friendliness.

The main aim of this paper is therefore to present an approach that allows cell meshes to be automatically generated for complex geometries, for example problems that contain intersections of boundary surfaces. We use immersed technology, familiar to the FEM community (see for example \cite{ZHANG20042051}). However, the unique nature of cells requires a re-engineering of published methods: No connectivity is required between the cells and - most importantly - we have to consider that volume integration is very compute intensive. The latter means that an intelligent design of the cell mesh is paramount.

In the following we first introduce the concept of the IGABEM and show how geometries, that are used for underground excavations, can be accurately (and sometimes even exactly) defined with few parameters. We also show how the use of NURBS for the approximation of the unknowns leads to elegant refinement procedures. 
Finally, we discuss the main topic of the paper namely the automatic generation of cell meshes for complex excavation geometries. 
On a practical example in oil extraction we demonstrate that the method works well and gives results comparable to ones obtained from a FEM simulation.

\section{B-splines and NURBS}

 B-splines are an attractive alternative to Lagrange polynomials and Serendipity functions predominantly used in simulation. The basis for creating the functions is the \textit{knot vector}. This is a vector containing a series of non-decreasing values of the local coordinate $\xi$:
\begin{equation}
\label{ }
\Xi=\left(\begin{array}{cccc}\xi_{0} & \xi_{1} & \cdots & \xi_{N}\end{array}\right)
\end{equation}

A B-spline basis function of order $p=0$ (constant) is given by:
\begin{equation}
N_{i,0}(\mathrm{\xi}) =	\left\{ \begin{array}{l}
			1 \quad \text{if}  \quad  \mathrm{\xi}_{i}\leqslant \mathrm{\xi} < \xi_{i+1} \\
			0 \quad \text{otherwise} \\
			\end{array} 
		\right.
\label{Bf1}
\end{equation}

\begin{figure}[H]
\begin{center}
\begin{overpic}[scale=0.5]{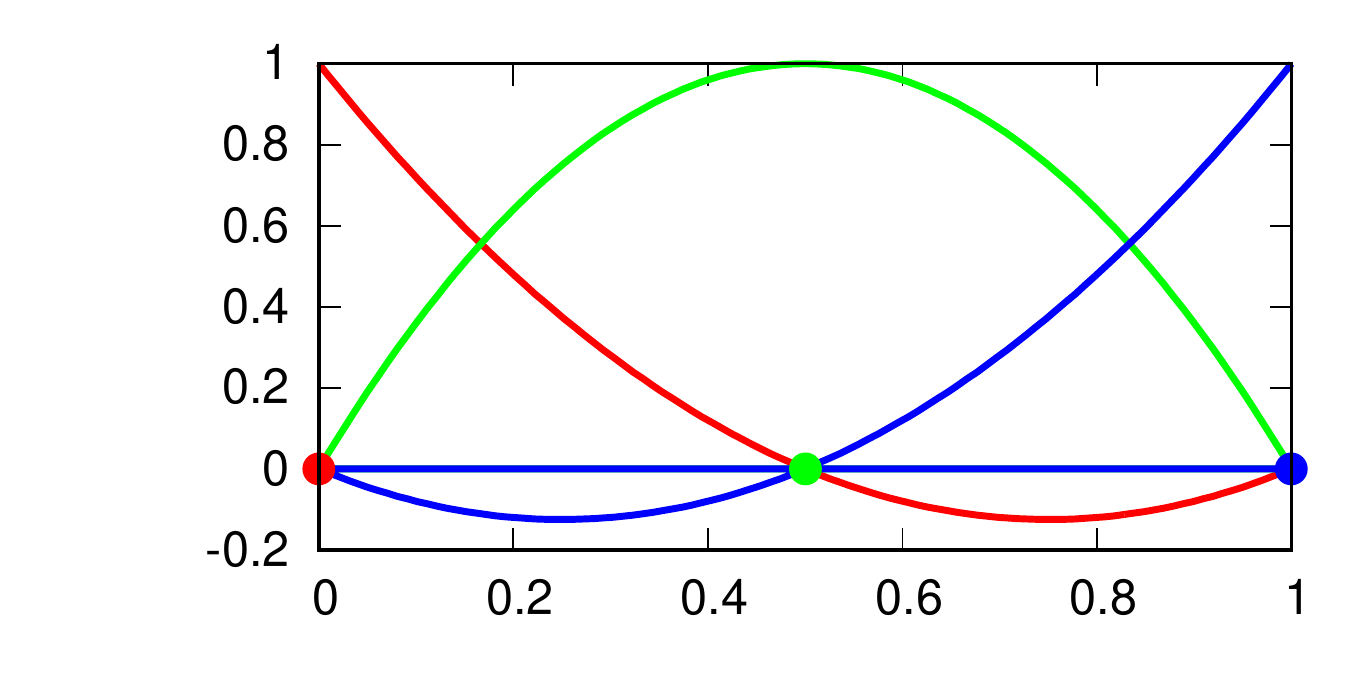}
 \put(50,0){$\xi$}
 \put(0,25){$L_{i}(\xi)$}
\end{overpic} 
\begin{overpic}[scale=0.5]{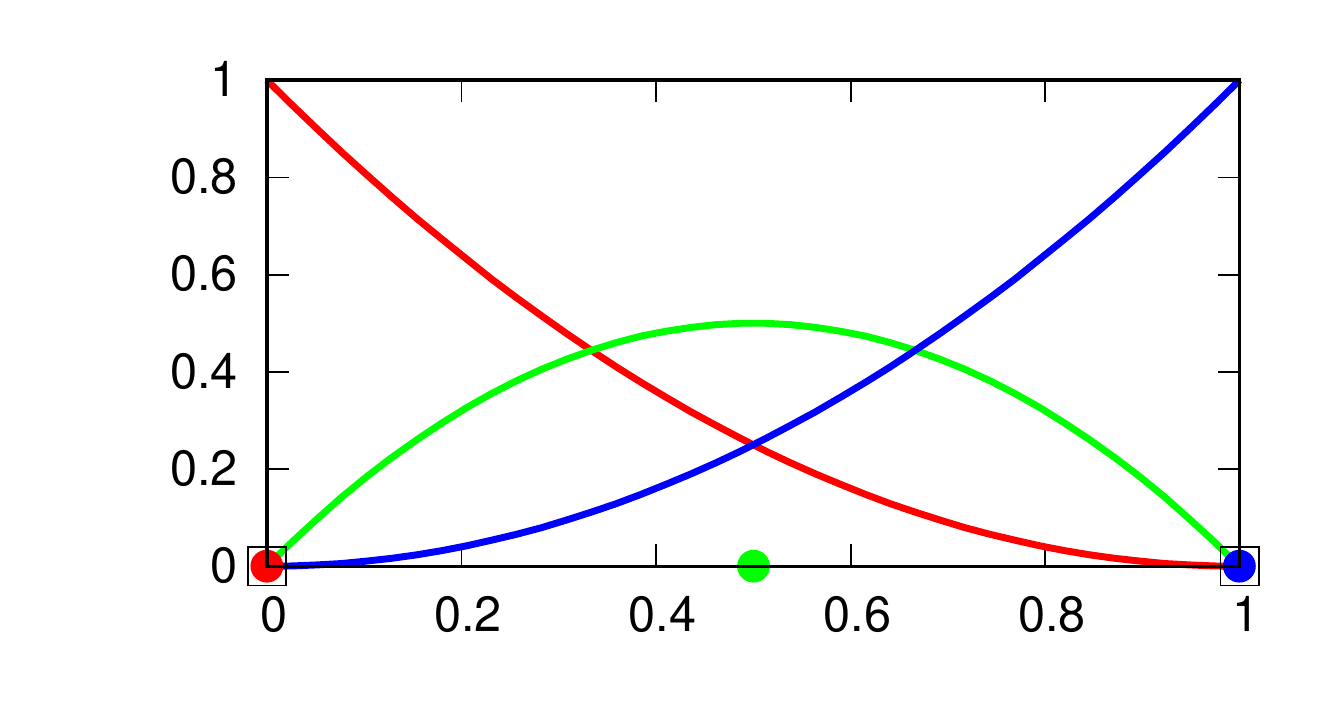}
 \put(50,0){$\xi$}
 \put(0,25){$N_{i}(\xi)$}
\end{overpic}
\caption{Comparison between Lagrange polynomials and B-splines of order 2. For the Lagrange polynomials coloured filled circles depict nodes; for the B-splines functions they depict anchors. Knots are shown as hollow squares.}
\label{LagBs}
\end{center}
\end{figure}

Higher order basis functions are defined by referencing lower order functions:
\begin{equation}
  N_{i,p}(\mathrm{\xi})=\frac{\xi-\xi_{i}}{\xi_{i+p}-\xi_{i}} \: N_{i,p-1}(\xi) 
			   + \frac{\xi_{i+p+1}-\xi}{\xi_{i+p+1}-\xi_{i+1}} \: N_{i+1,p-1}(\xi)
\label{Bf2}
\end{equation}
Figure \ref{LagBs} shows a comparison between Lagrange polynomials and B-splines. 

It can be seen that B-splines have only positive values but the main difference is that the basis functions are not associated with nodes as is the case with Lagrange polynomials. Instead we associate the functions with \textit{anchors} which are depicted as filled circles in the figure.
The location of the $i$-th anchor in the parameter space can be computed by:
\begin{equation}
\label{Greville}
% \uu(a_{i})
\abscissa_i = \frac{\uu_{i+1}+\uu_{i+2} + \dots +\uu_{i+p}}{p} \qquad i=0,1, \dots ,I-1.
\end{equation}

\begin{figure}[H]
\begin{center}
\begin{overpic}[scale=0.6]{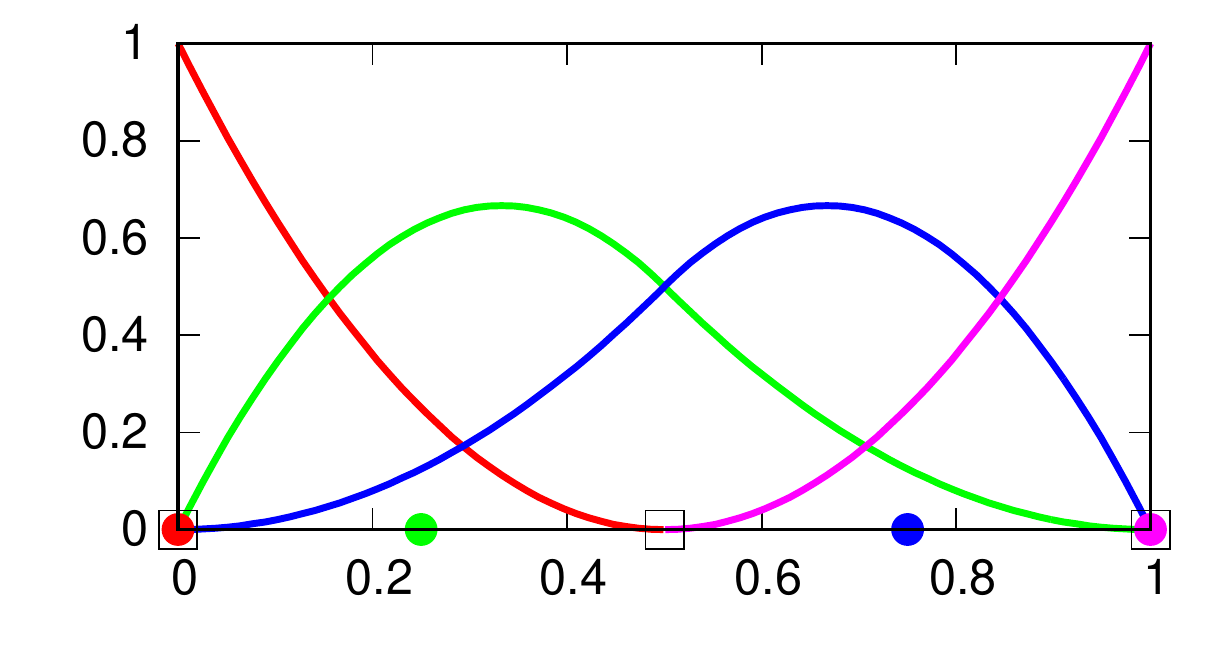}
 \put(50,0){$\xi$}
\end{overpic} 
\caption{Effect one knot insertion at $\xi=0.5$ for a B-spline of order 2. This changes the shape of the basis functions and increases their number.}
\label{Bsk}
\end{center}
\end{figure}

However the big advantage of B-splines is that the order of the function and the continuity can be changed by simply changing the knot-vector.
For the B-spline of order 2 the knot vector is given by:
\begin{equation}
\label{ }
\Xi=\left\{\begin{array}{cccccc}0 & 0 & 0 & 1 & 1 & 1\end{array}\right\}
\end{equation}
where the repetition ($p+1$)-times of the first and last entry defines the order.
Furthermore we can also change the continuity of the function by inserting knots (see \myfigref{Bsk}).
For example a knot vector:
\begin{equation}
\label{ }
\Xi=\left\{\begin{array}{ccccccc}0 & 0 & 0 & 0.5 & 1 & 1 & 1\end{array}\right\}
\end{equation}
reduces the continuity at $\xi=0.5$ by one order to $C^{1}$.

\remark{The fact that we are able to change the order and continuity of the basis functions and that they do not depend on nodal points is the most significant difference to conventional basis functions. This will be exploited later when we define the approximation of the unknown.}

\bigskip

Nonuniform rational B-splines or NURBS are based on B-splines but have improved properties for the definition of geometry.
NURBS of order $p$ are defined as:
\begin{equation}
\label{N1D}
  R_{i,p}(\xi)=\frac{N_{i,p}(\xi) \: \mathrm{w}_{i}}{\sum_{j=0}^{I}N_{j,p}(\xi) \: \mathrm{w}_{j}}  
\end{equation}
where $I+1$ is the number of basis functions and $\mathrm{w}_{i}$ are weights. 

\section{Geometry definitions by NURBS - Curves}
NURBS are ideally suited for the description of geometry. 
To define a curve we can write:
\begin{equation}
\mathbf{x} = \sum_{i=1} ^{I}R_{i}(\xi) \: \textbf{c}_{i}
\label{Map1}
 \end{equation}
 where $R_{i}(\xi)$ are the basis functions defined in \myeqref{N1D} except that the numbering starts from 1 instead of 0, $\textbf{c}_{i}$ are control point coordinates. NURBS allow certain geometrical shapes for example arcs to be described exactly.
 
 \newpage
 
 \subsection{Defining an arc}
A NURBS curve can exactly define an arc.
  \begin{figure}[H]
\begin{center}
\includegraphics[scale=0.7]{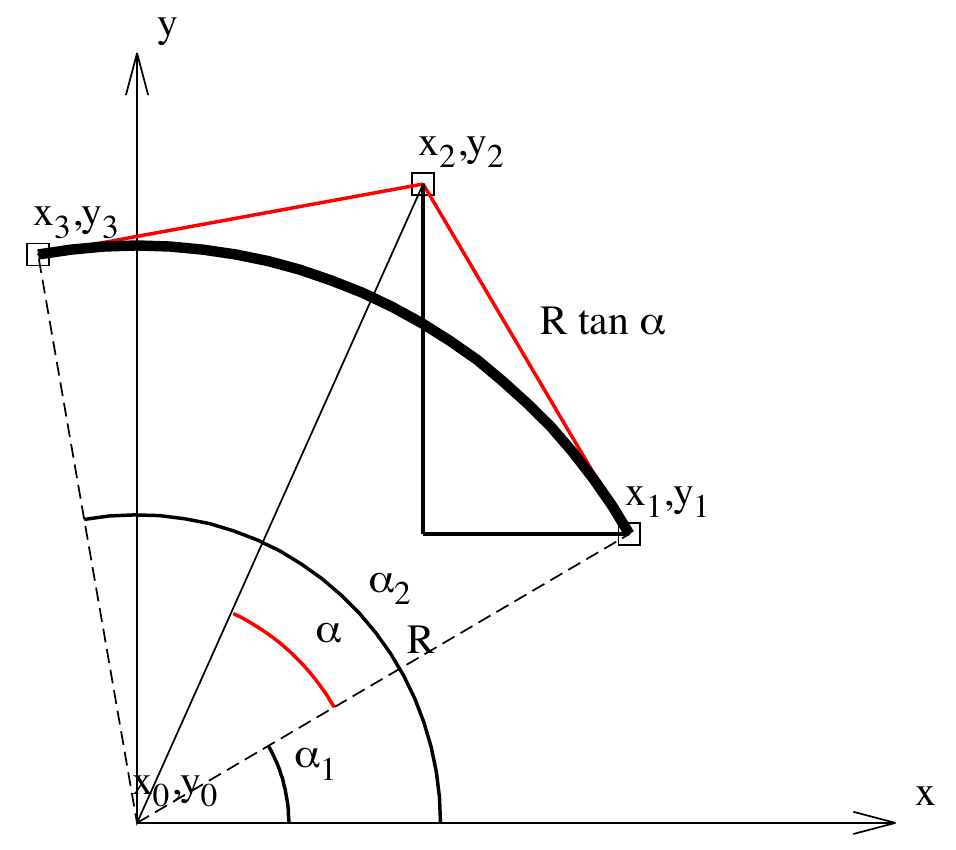}
\caption{Explanation for determining the control points for an arc.}
\label{Arc}
\end{center}
\end{figure}

The $x,y$ coordinates of the control points are given by (see \myfigref{Arc}):
\begin{eqnarray}
\label{eq1}
\mathbf{ c}_{1} &=&\left[\begin{array}{c}x_{0} + R \cdot \cos \alpha_{1} \\y_{0} + R \cdot \sin \alpha_{1} \end{array}\right] \\
\nonumber
\mathbf{ c}_{2} &=&\left[\begin{array}{c}x_{1} - R\cdot \tan \alpha \cdot \sin \alpha_{1}  \\ y_{1} + R\cdot \tan \alpha \cdot \cos \alpha_{1} \end{array}\right] \\
\nonumber
\mathbf{ c}_{3} &=&\left[\begin{array}{c}x_{0} + R \cdot \cos \alpha_{2}  \\ y_{0} + R \cdot \sin \alpha_{2} \end{array}\right] 
\end{eqnarray}
where $\alpha= (\alpha_{2} -\alpha_{1})/2$. The weights are $w(1)=1 \ ; w(2)= \cos \alpha \ ; w(3)=1$.

 \subsection{Description of an NATM tunnel}
\begin{figure}[H]
\begin{center}
\includegraphics[scale=0.8]{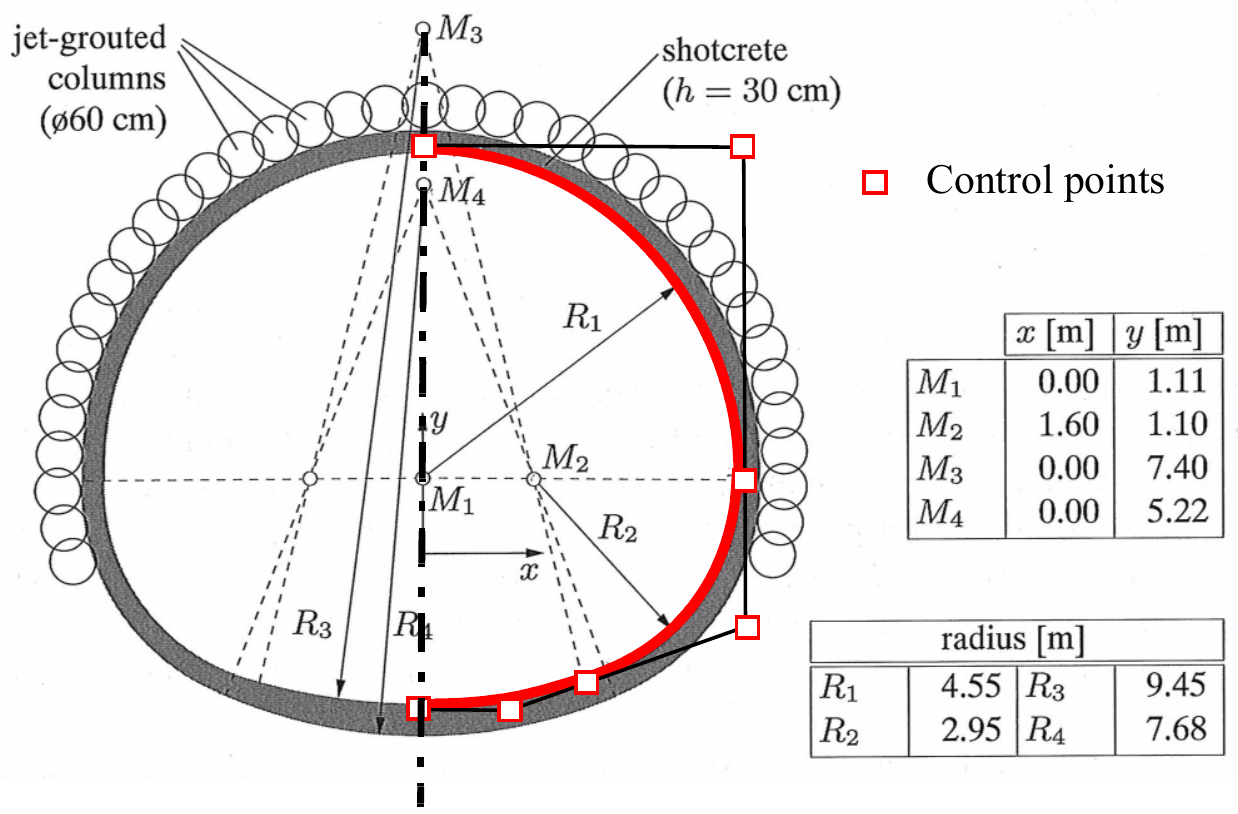}
\caption{A NATM tunnel with a table of input data describing the shape. Also shown are the computed control points as hollow squares and in red the exact definition by a NURBS curve.}
\label{NATM}
\end{center}
\end{figure}

As a practical example we show how the design shape of a NATM tunnel can be exactly described by NURBS curves.
 \myfigref{NATM} shows an example of a tunnel together with a table of values that define the shape.
Algorithm \ref{alg:NATM} shows how to compute the control point coordinates and weights of the three arcs for half the tunnel using the values in the table only.
The results are shown in  \myfigref{NATM}.
 \begin{myalgorithm}{Algorithm for computing the control point coordinates and weights for NATM Tunnel}{alg:NATM}
	\REQUIRE  M(2,3), R(3)
	\STATE knot vector= [0 0 0 1 1 1]
	\STATE dx= M(1,2) - M(1,3) ; dy= M(2,2) - M(2,3)
        \STATE $\phi= tan^{-1}(dx/dy)$
	\FOR {$n=1$ \TO 3}
	 \IF{n = 1}
          \STATE $\alpha_1= 0$ ;  $\alpha_2= \pi/2$
         \ENDIF
         \IF{n = 2}
          \STATE $\alpha_1= 3/2\pi + \phi$ ; $\alpha_2= 2\pi $
         \ENDIF
          \IF{n =3}
          \STATE $\alpha_1= 3/2\pi$ ; $\alpha_2= 3/2\pi + \phi$
         \ENDIF
         \STATE  Compute control points coordinates and weights for arc $n$ 
         \STATE  with $\left\{\begin{array}{c}x_0 \\y_0\end{array}\right\}$=M(:,n), R(n), $\alpha_1,\alpha_2$, using \myeqref{eq1}
        \ENDFOR
	\RETURN Knot vector, control point coordinates and weights for 3 arcs
\end{myalgorithm}

\newpage

\section{Geometry definition with NURBS- surfaces}
We can expand the NURBS to three dimensions using a tensor product:
\begin{equation}
\mathbf{x} = \sum_{i=0} ^{I}  \sum_{j=0} ^{J} R_{i,p}(\xi) R_{j,q}(\eta) \: \textbf{c}_{ij}
\label{Map0}
 \end{equation}
where $q$ is the order and $J+1$ is the number of basis functions in the local $\eta$ direction. $\textbf{c}_{ij}$ are control point coordinates.  
We now need an additional knot vector $H$ that defines the basis function in $\eta$ direction.

\subsection{Definition of a surface}
\label{Surf}
As an example we show in \myfigref{Surf} a surface created with the knot vectors:
\begin{eqnarray}
\label{ }
\Xi &=& \left[\begin{array}{cccccc}0 & 0 & 0 & 1 & 1 & 1\end{array}\right] \\
\nonumber
H &=& \left[\begin{array}{cccc}0 & 0 & 1 & 1\end{array}\right]
\end{eqnarray}
i.e. quadratic in $\xi$-direction and linear in $\eta$ direction. The control point coordinates and weights are shown in table \ref{tab:Surf}.
\begin{mytable}
  {H}               %% table position  
  {Control point coordinates and weights for definition of surface.}  %% caption
  {tab:Surf}  %% label
  {cccccc}         %% column layout
  \mytableheader{ $i$ & $j$ & x & y & z & w }  %% header
  %% table content
0 & 0 & 0 & 1 & 0 & 1 \\
1 & 0 & 0 & 1 & 1 & 0.707 \\
2 & 0 & 0 & 0 & 1 & 1 \\
0 & 1 & 1 & 1 & 0 & 1 \\
1 & 1 & 1 & 1 & 1 & 0.707 \\
2 & 1 & 1 & 0 & 1 & 1 \\
\end{mytable}%
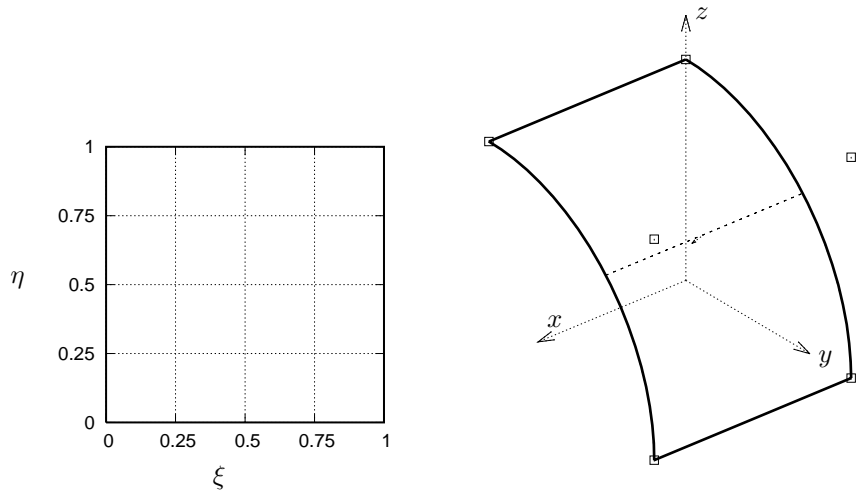
\begin{figure}
\begin{center}
\begin{overpic}[scale=0.5]{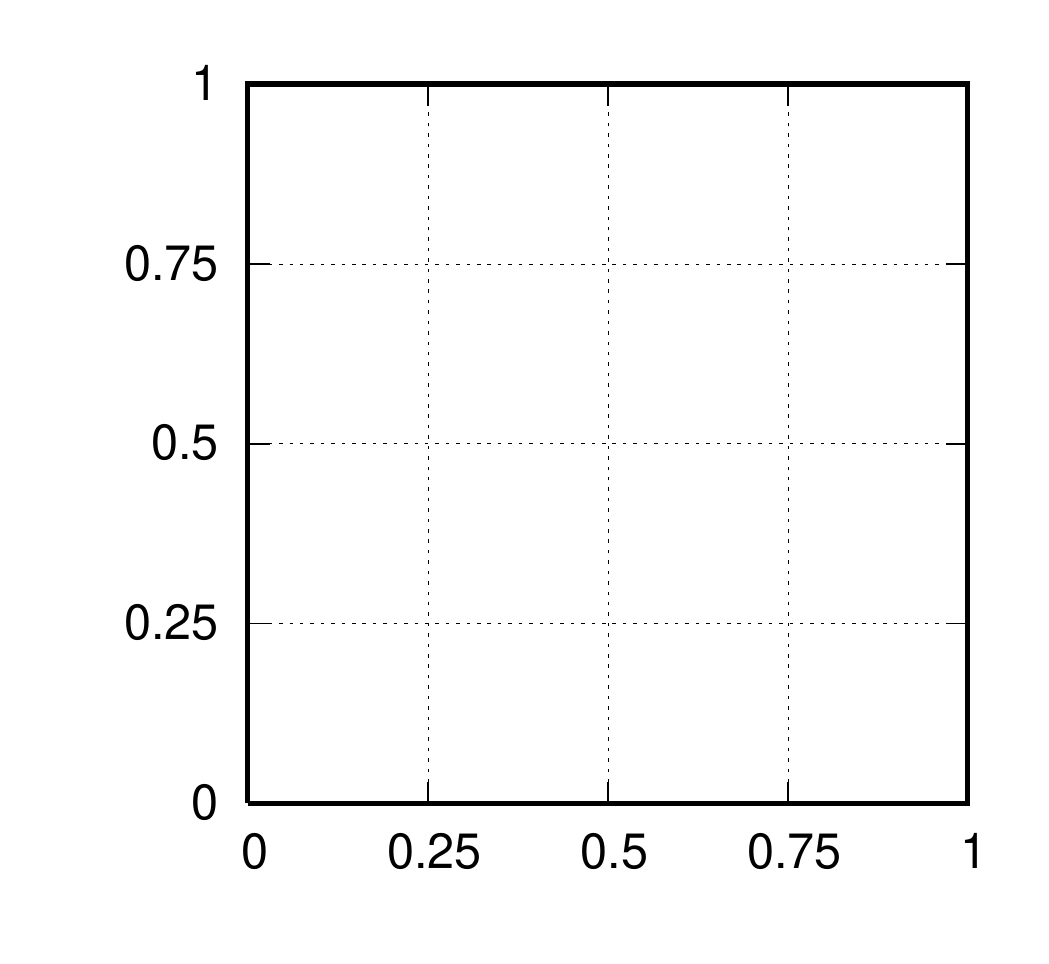}
 \put(50,0){$\xi$}
 \put(0,50){$\eta$}
\end{overpic}
\begin{overpic}[scale=0.5]{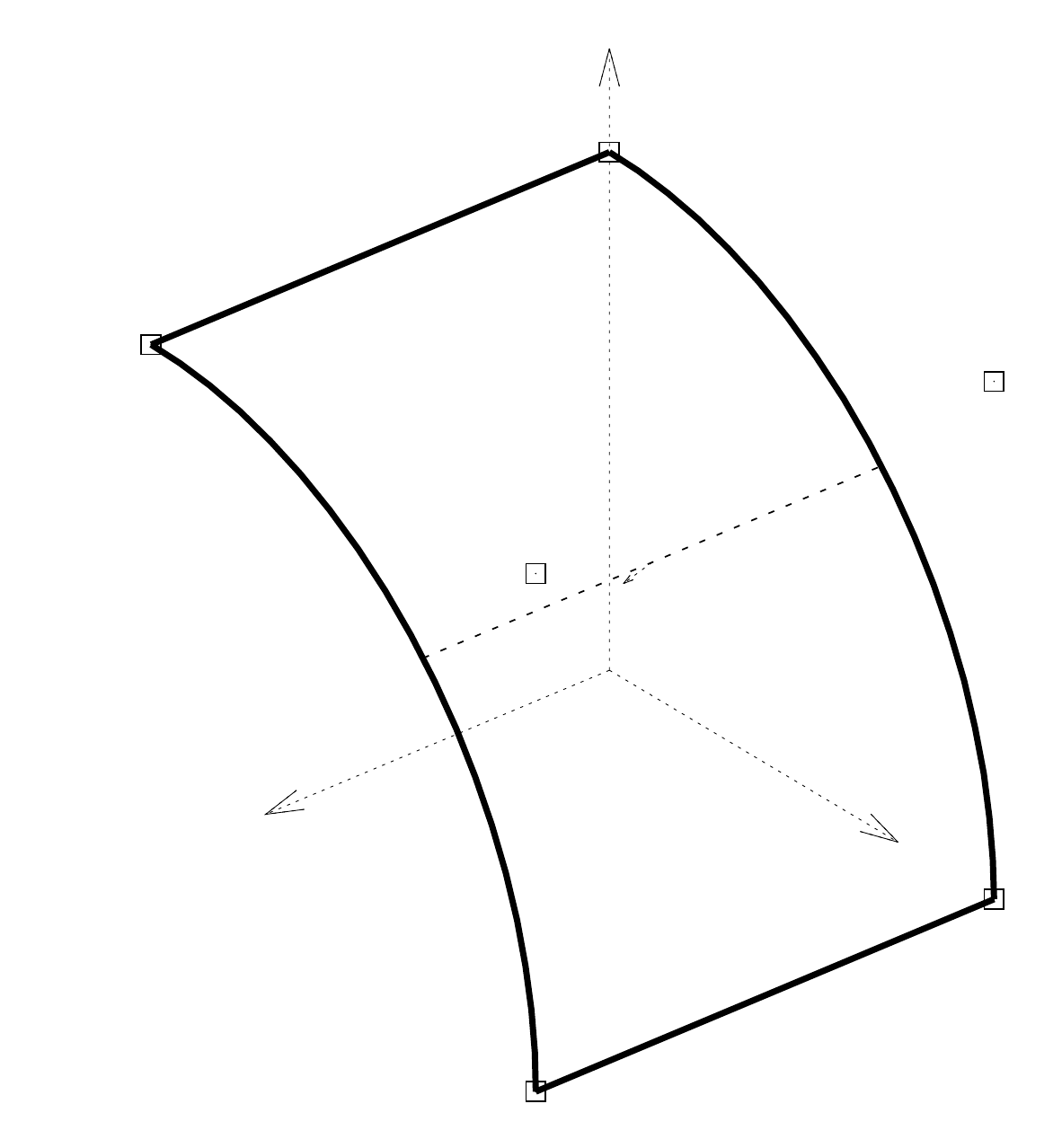}
 \put(80,25){$y$}
 \put(25,32){$x$}
 \put(55,95){$z$}
\end{overpic} 
\caption{A NURBS surface left in local, right in global coordinate system showing control points as hollow squares.}
\label{Surf}
\end{center}
\end{figure}
In the following we will show geometric manipulations with NURBS.

\newpage

\begin{figure}
\begin{center}
\begin{overpic}[scale=0.5]{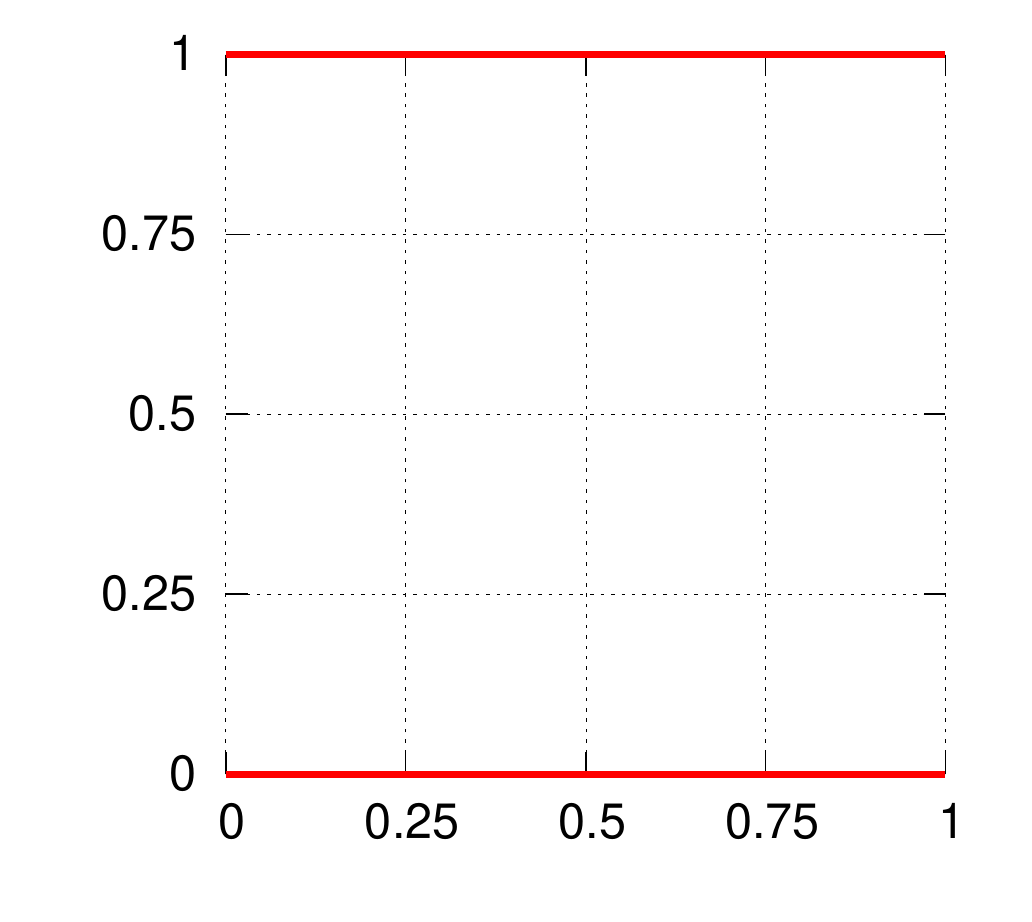}
 \put(50,0){$\xi$}
 \put(0,50){$\eta$}
\end{overpic}
\begin{overpic}[scale=0.5]{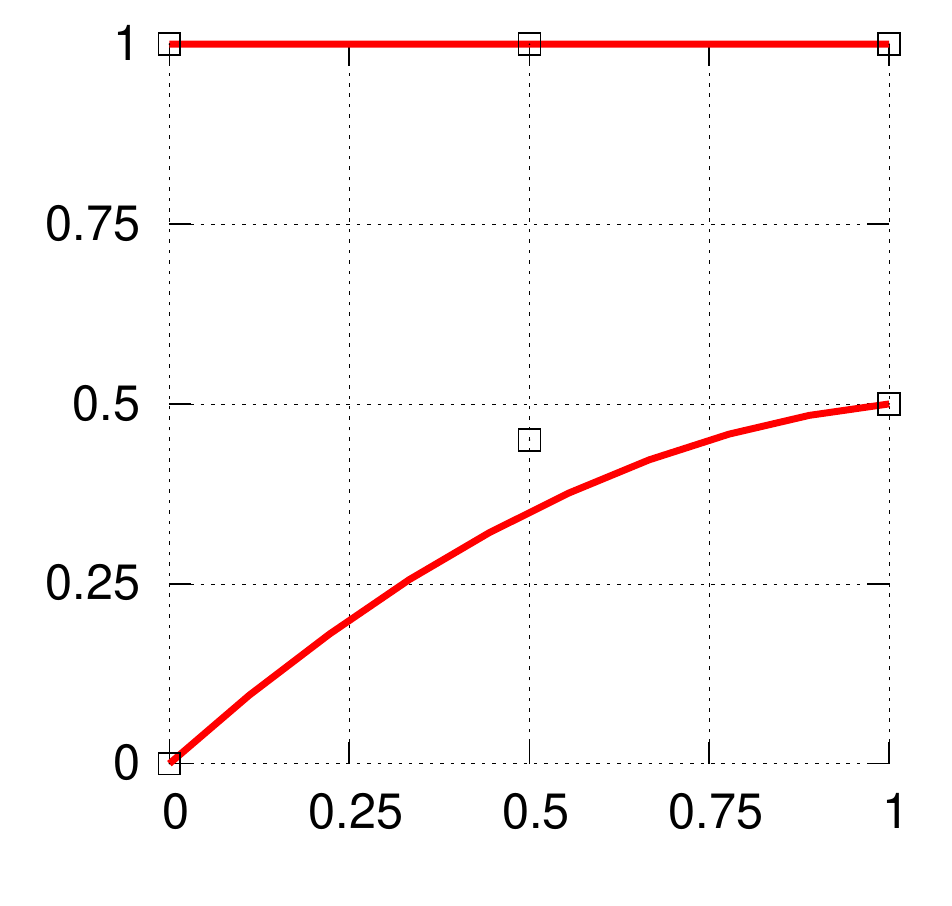}
 \put(50,0){$\hat{\xi}$}
 \put(2,50){$\hat{\eta}$}
\end{overpic}
\begin{overpic}[scale=0.5]{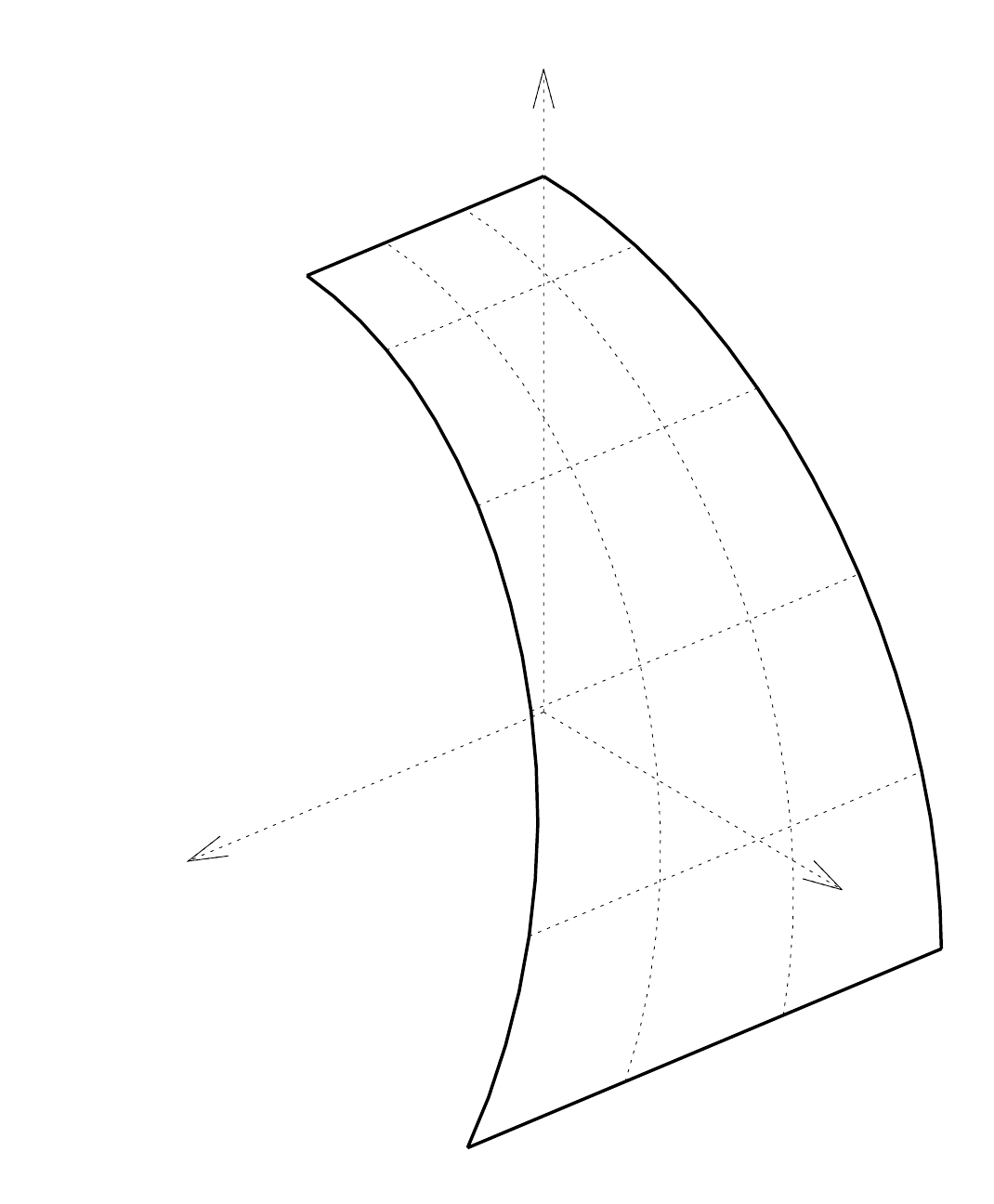}
 \put(75,25){$y$}
 \put(10,25){$x$}
 \put(50,95){$z$}
\end{overpic} 
\caption{Trimming of a surface:. From left to right: Local map, trimmed map (showing trimming curves and their control points) and global map of trimmed surface}
\label{SurfT}
\end{center}
\end{figure}
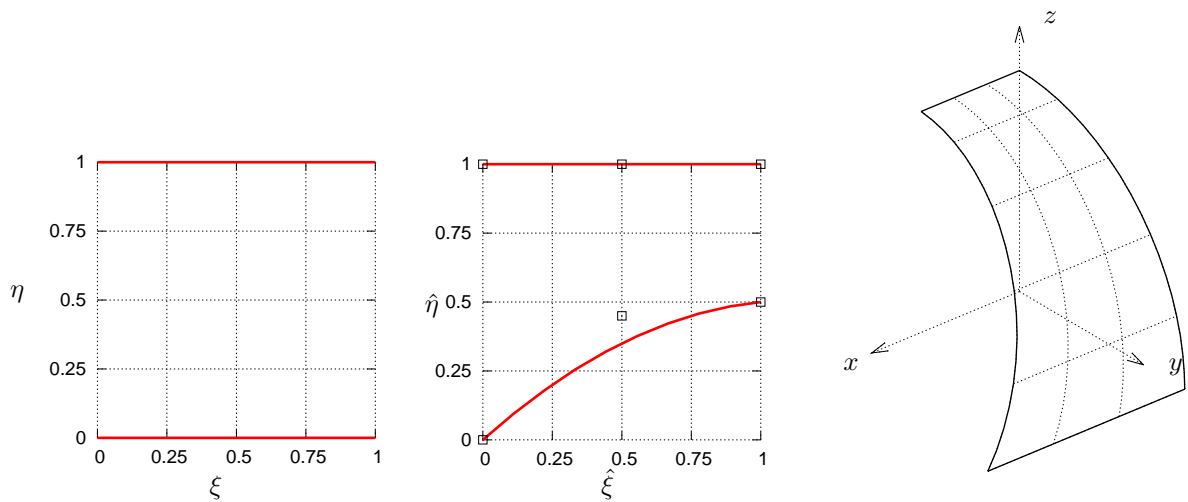

\subsection{Trimming}

We can cut off a portion of the surface by a process called trimming.

We first define the limits of the trimmed space using basis functions $\hat{R}_{i,p}(\xi) $ and  $ \hat{R}_{j,q}(\eta)$ and control points $ \hat{\textbf{c}}_{ij}$.
We then map from a $\xi,\eta$ space to a local $\hat{\xi},\hat{\eta}$ space:
\begin{equation}
\left\{\begin{array}{c}\hat{\xi} \\ \hat{\eta} \end{array}\right\} = \sum_{i=0} ^{I}  \sum_{j=0} ^{J} \hat{R}_{i,p}(\xi) \hat{R}_{j,q}(\eta) \: \hat{\textbf{c}}_{ij}
\label{}
 \end{equation}
Finally we map form the $\hat{\xi},\hat{\eta}$ coordinate system to the $x,y,z$ coordinate system using \myeqref{Map0}:
 \begin{equation}
\mathbf{x} = \sum_{i=0} ^{I}  \sum_{j=0} ^{J} R_{i,p}(\hat{\xi}) R_{j,q}(\hat{\eta}) \: \textbf{c}_{ij}
\label{}
 \end{equation}
 In \myfigref{SurfT} we show an example on how the surface shown in \myfigref{Surf} can be trimmed.

\newpage

\section{The boundary element method without volume effects}
Here we introduce the boundary element method for the case where only boundary effects are considered. This only allows simulations of elastic, homogeneous and isotropic domains. Later, we will introduce volume effects, that allow to remove these restrictions. 

Consider an infinite domain  $\domain$ bounded by $\boundary$.
To establish the integral equations we apply Betti's theorem and the collocation method.
The regularised integral equations are written in matrix notation as (for a derivation see \cite{BeerMarussig}).:
\begin{equation}
\label{Inteq}
    \begin{aligned}
\int_{\boundary} \fund{T}(\sourcept_{n},\fieldpt) ( \myVec{u}(\fieldpt) - \myVec{u}(\sourcept_{n})) \ d\boundary(\fieldpt) - \mathbf{ A}_{n} \myVec{u}(\sourcept_{n}) &=& \int_{\boundary} \fund{U}(\sourcept_{n},\fieldpt) \ \myVec{t}(\fieldpt) \ d\boundary(\fieldpt) .
\end{aligned}
\end{equation}
where $\fund{T}(\sourcept_{n},\fieldpt)$ is a matrix that contains the fundamental solutions for the tractions and $\fund{U}(\sourcept_{n},\fieldpt)$ the one for the displacements at point $\fieldpt$ due to a source at $\sourcept_n$ on the boundary of an infinite, elastic domain.
$\myVec{u}$ and $\myVec{t}$ are displacement and traction vectors on $\boundary$.
 $\mathbf{ A}_{n}$ represents the azimuthal integral that is equal to the unit matrix for infinite domain problems.
The fundamental solutions for the displacements and tractions are well published and summarised in \cite{BeerMarussig}. The points with coordinates $\sourcept_n$ are also known as \textbf{collocation points}. 

\subsection{Discretisation of integral equations}
To be able to solve the integral equations they have to be discretised.
The discretisation involves  the subdivision of the boundary into patches and the approximation of the unknown boundary values.
Patches are similar to boundary elements but have significantly more flexibility regarding the geometry they represent.
The geometrical approximation of the integral equations by discretisation can be written as:
\begin{equation}
    \begin{aligned}
        \label{Regular2}
 \int_{\boundary} \fund{U}(\sourcept_{n},\fieldpt) \ \myVec{t}(\fieldpt) \ d\boundary(\fieldpt) - \int_{\boundary} \fund{T}(\sourcept_{n},\fieldpt) ( \myVec{u}(\fieldpt) - \myVec{u}(\sourcept_{n})) \ d\boundary(\fieldpt) + \mathbf{ A}_{n} \myVec{u}(\sourcept_{n})  \approx  \\
        \sum_{e=1}^{E} \int_{\Gamma_{e}}  \fund{U}(\sourcept_{n},\fieldpt) \ \myVec{\dual}^{e}(\fieldpt)\ d\Gamma_{e}(\fieldpt)  
        - \sum_{e=1}^{E} \int_{\Gamma_{e}} \fund{T}(\sourcept_{n},\fieldpt)  \myVec{\primary}^{e}(\fieldpt) d \Gamma_{e}(\fieldpt)  \\
        +  \left[\sum_{e=1}^{E} \ \left(\int_{\Gamma_{e}} \fund{T}(\sourcept_{n},\fieldpt) d \Gamma_{e} (\fieldpt)\right) + \mathbf{A}_{n}\right] \myVec{\primary}(\sourcept_{n})
    \end{aligned}
\end{equation}
where $e$ specifies the patch number and $E$ is the total number of patches.

\subsection{Patches}
Patches are similar to boundary elements that use Lagrange polynomials or Serendipity functions, except that they can define a much larger portion of the surface accurately and sometimes exactly.
For the simulation in geomechanics it is convenient to have two type of patches: finite and semi-infinite patches. The latter are able to describe surfaces that go to infinity.

\subsubsection{Geometry definition of finite patches}
The mapping from the local patch coordinates $\myVecGreek{\xi}=(\xi,\eta)^{\mathrm{T}}=[0,1]^2$ to the global $\pt{x}$ coordinate system is given by
\begin{equation}
\label{mapatch}
\pt{x} (\xi,\eta)= \sum_{i=1}^{I} \NURBS_{i} (\xi, \eta) \pt{c}_{i}.
\end{equation}
where $\NURBS_{i} (\xi, \eta)$ are now combined NURBS basis functions (\ref{Map0}) with the control points (coordinates $\pt{c}_{i}$) numbered consecutively from one, first in the $\uu$- and then in the $\vv$-direction. 

\begin{figure}
\begin{center}
\begin{overpic}[scale=0.5]{pics/SurfaceL.pdf}
 \put(50,0){$\uu$}
 \put(0,50){$\vv$}
\end{overpic}
\begin{overpic}[scale=0.5]{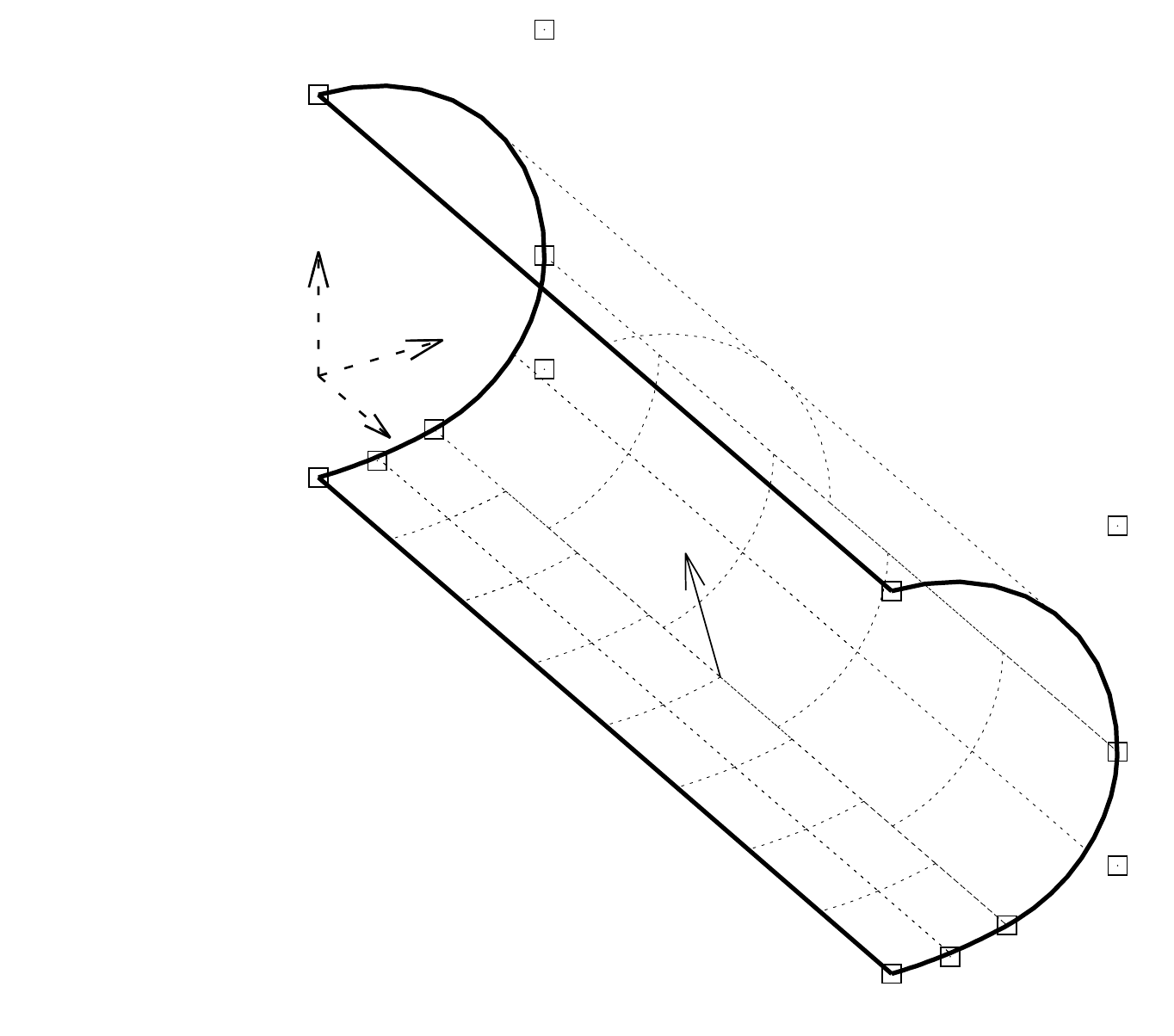}
 \put(40,60){$y$}
 \put(35,47){$x$}
 \put(27,70){$z$}
\end{overpic}
\caption{A finite patch representing half of an NATM tunnel. Left: local map, Right: global map, showing control points as hollow squares and the outward normal.}
\label{Patch3D}
\end{center}
\end{figure}
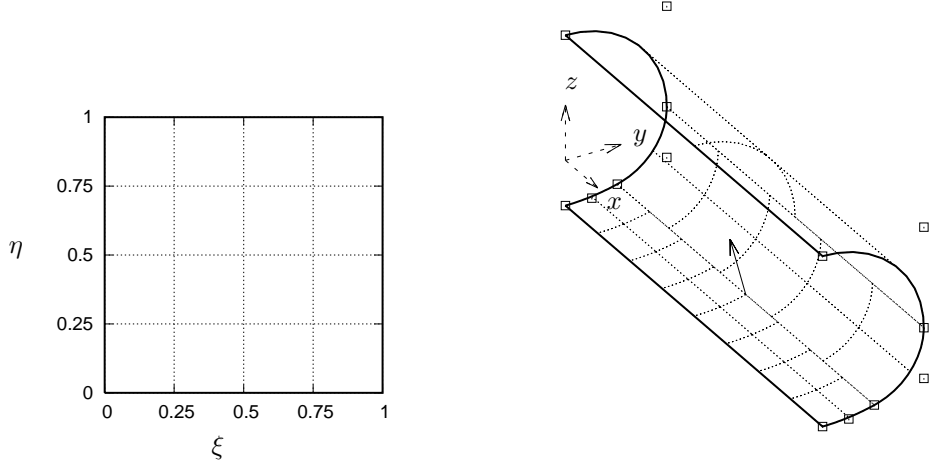

The vectors tangential to the surface are given by

\begin{align}
%\label{ }
\mathbf{ v}_{\xi}= \frac{\partial \pt{x}}{\partial \xi}= \left(\begin{array}{c}\frac{\partial x}{\partial \xi} \\ \\ \frac{\partial y}{\partial \xi} \\ \\ \frac{\partial z}{\partial \xi}\end{array}\right)
&& \text{and} && \mathbf{ v}_{\eta}= \frac{\partial \pt{x}}{\partial \eta}= \left(\begin{array}{c}\frac{\partial x}{\partial \eta} \\ \\ \frac{\partial y}{\partial \eta} \\ \\ \frac{\partial z}{\partial \eta} \end{array}\right)
\end{align}
and the unit vector normal to the surface is given by:
\begin{equation}
%\label{ }
\mathbf{ n}= \frac{\mathbf{ v}_{\xi} \times  \mathbf{ v}_{\eta}}{J}.
\end{equation}
The Jacobian is
\begin{equation}
%\label{ }
J= | \mathbf{ v}_{\xi} \times  \mathbf{ v}_{\eta} |.
\end{equation}
 The  direction of the outward normal depends on how the control points are numbered.
 To show that a patch can be much more geometrically complex than a boundary element we show in \myfigref{Patch3D} an example of a finite patch that represents exactly the geometry of  half a NATM tunnel.

\subsubsection{Geometry definition of semi-infinite patches}
\label{sec:InfinitePatches}
Here we introduce a patch definition that is able to model a surface that extends to infinity on one side that has been published in \cite{Beer2015b}. 
The mapping for a patch that extends to infinity in the $\eta$-direction is given by
\begin{equation}
%\label{ }
\pt{x}= \sum_{j=1}^{2}\sum_{i=0}^{I} \NURBS^{\infty}_{ij}(\xi,\eta) \pt{c}_{ij}
\end{equation}
where
\begin{equation}
%\label{ }
\NURBS^{\infty}_{ij}(\xi,\eta)= \NURBS_{i}(\xi) M_{j}^{\infty}(\eta)
\end{equation}
The special infinite basis functions are given by:
\begin{align}
 M_{1}^{\infty} =  \frac{1 - 2\eta}{1-\eta} && \text{and} &&  M_{2}^{\infty}  =  \frac{\eta}{1-\eta}. 
\end{align}
\begin{figure}[H]
\begin{center}
\begin{overpic}[scale=0.5]{pics/SurfaceL.pdf}
 \put(50,0){$\uu$}
 \put(0,50){$\vv$}
\end{overpic}
\begin{overpic}[scale=0.5]{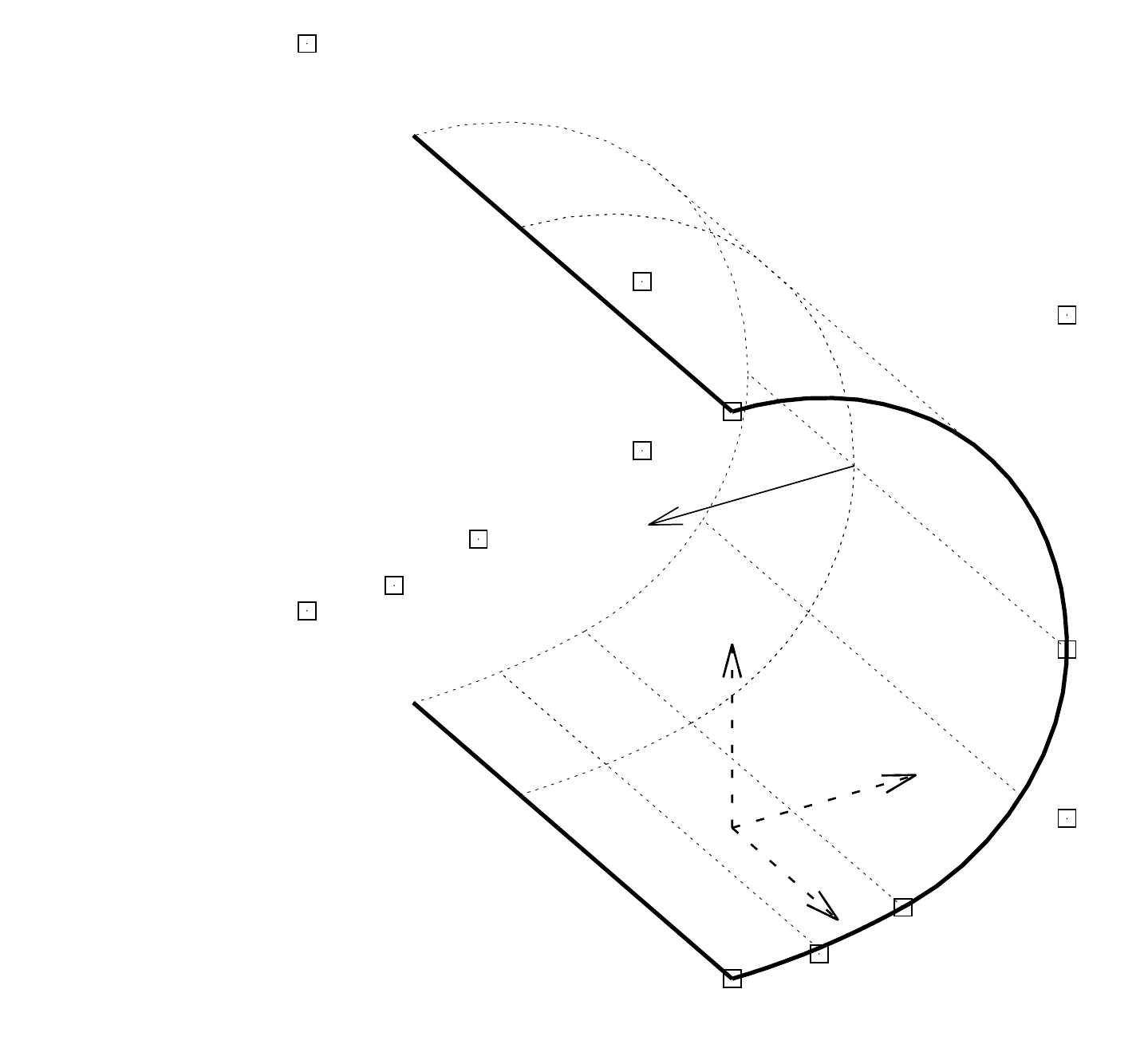}
 \put(82,25){$y$}
 \put(74,12){$x$}
 \put(60,35){$z$}
\end{overpic}
\caption{An infinite patch. Left: local map, Right: global map, showing control points and outward normal.}
\label{Patch3DI}
\end{center}
\end{figure}
\begin{figure}[H]
\begin{center}
\includegraphics[scale=0.5]{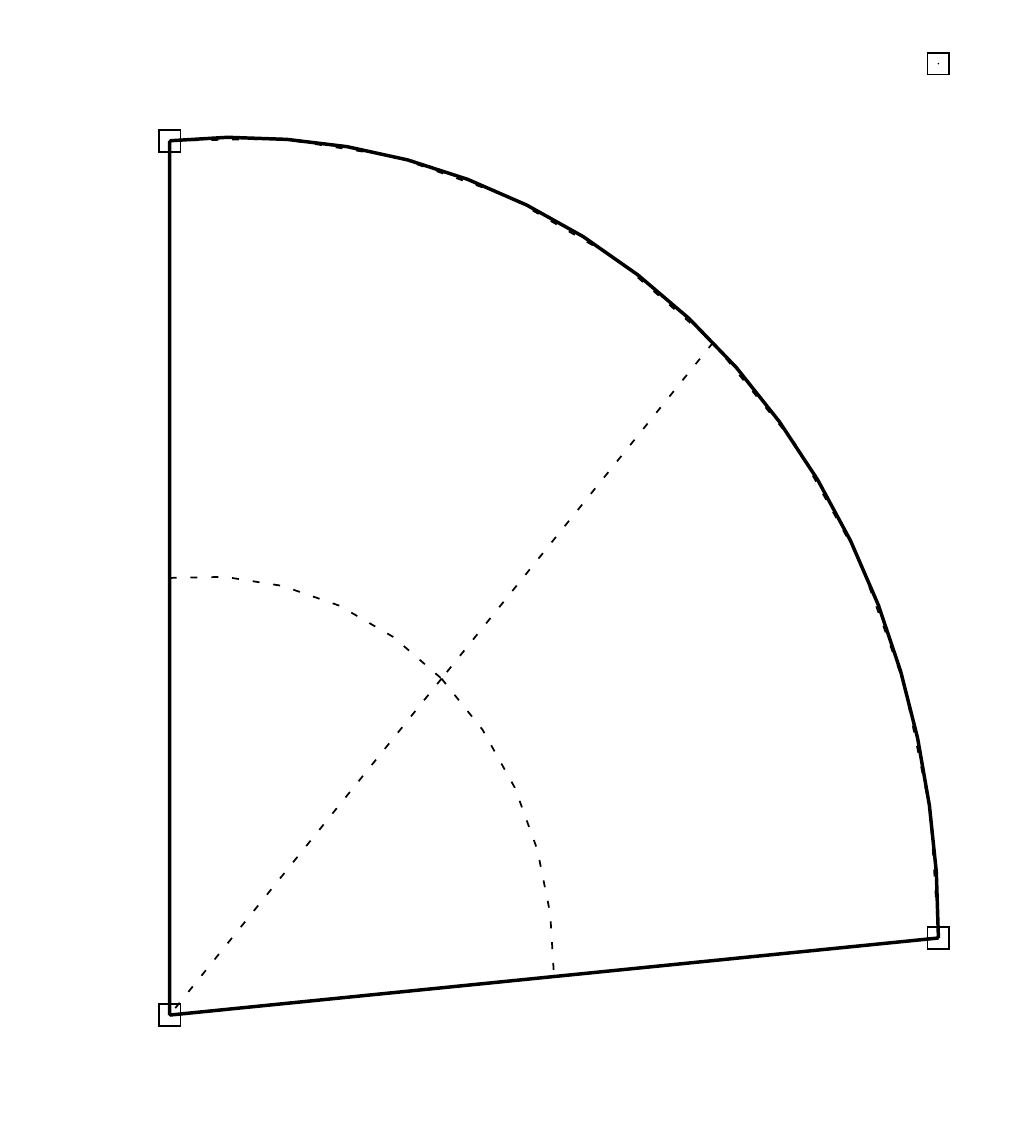}
\caption{An degenerate patch: Control points on lower left are coincidental }
\label{Patchdeg}
\end{center}
\end{figure}
The vectors in the tangential directions are given by
\begin{eqnarray}
%\label{ }
\mathbf{ v}_{\xi}= \frac{\partial \pt{x}}{\partial \xi}=\sum_{j=1}^{2}\sum_{i=0}^{I} \frac{\partial \NURBS_{i}(\xi)}{\partial \xi}  M_{j}^{\infty}(\eta) \pt{c}_{ij} \\
\mathbf{ v}_{\eta}= \frac{\partial \pt{x}}{\partial \eta}=\sum_{j=1}^{2}\sum_{i=0}^{I} \NURBS_{i}(\xi)  \frac{\partial M_{j}^{\infty}(\eta)}{\partial \eta} \pt{c}_{ij} 
\end{eqnarray}
where
\begin{align}
%\label{ }
\frac{\partial M_{1}^{\infty}}{\partial \eta}= \frac{-1}{(1-\eta)^{2}} && \text{and} && \frac{\partial M_{2}^{\infty}}{\partial \eta}= \frac{1}{(1-\eta)^{2}}.
\end{align}
The unit vector normal to the patch is given by:
\begin{equation}
%\label{ }
\mathbf{ n}= \frac{\mathbf{ v}_{\xi} \times  \mathbf{ v}_{\eta}}{J}.
\end{equation}
The Jacobian is
\begin{equation}
%\label{ }
J= | \mathbf{ v}_{\xi} \times  \mathbf{ v}_{\eta} |.
\end{equation}
 It is noted that the Jacobian $J$ tends to infinity as $\vv$ tends to 1.
\myfigref{Patch3DI} shows an example of an infinite patch representing half of an NATM tunnel.

\subsubsection{Geometry definition of degenerate patches}
Patches of triangular shape can be defined by making 3 control points of the patch coincidental (\myfigref{Patchdeg}). Although the Jacobian tends to zero at the degenerate point, this is of no consequence because integration points are not located there.

\subsection{Approximation of the boundary values}
To be able to solve the integral equations we need to approximate the boundary values.
We approximate the unknown boundary displacements in patch $e$ by:
\begin{equation}
\label{uapprox}
\mathbf{ u}^{e}(\xi,\eta)= \sum_{k=1}^{K}\hat{ \NURBS_{k} }(\xi, \eta) \mathbf{ u}^{e}_{k}.
\end{equation}
where $\hat{\NURBS}_{k} (\xi, \eta)$ are NURBS basis functions, $K$ is the number of parameters and $\mathbf{ u}^{e}_{k}$ are parameter values. 

\bigskip

\remark{Note that in contrast to conventional analysis the isogeometric BEM works with parameter values instead of real values of the unknown}.

\bigskip

For the excavation methods, that we consider here, the displacements will be unknown. The tractions in patch $e$ can be computed from the stress field that exists in the domain prior to excavation (virgin stress):
\begin{equation}
\label{tdef}
\mathbf{ t}^{e}_{0}(\fieldpt)= \mathbf{ n} (\fieldpt) \myVecGreek{\sigma}_\mathrm{v} (\fieldpt)
\end{equation}
where $ \myVecGreek{\sigma}_\mathrm{v} $ is the (pseudo) virgin stress vector in Voigt notation:
\begin{equation}
\label{Voightnot}
\myVecGreek{\sigma}_{\mathrm{v}}= \left \{\begin{array}{c}\sigma_{x\mathrm{v}} \\\sigma_{y\mathrm{v}} \\\sigma_{z\mathrm{v}} \\\tau_{xy\mathrm{v}} \\\tau_{yz\mathrm{v}}\\\tau_{xz\mathrm{v}}\end{array}\right\}
\end{equation}
and
\begin{equation}
\label{}
\mathbf{ n}(\fieldpt)= \left[\begin{array}{cccccc}n_{x} & 0 & 0 & n_{y} & 0 & n_{z} \\0 & n_{y} & 0 & n_{x} & n_{z} & 0 \\0 & 0 & n_{z} & 0 & n_{y} & n_{x}\end{array}\right]
\end{equation}
where $n_{x}, n_{y}, n_{z}$ are components of the vector normal to the surface.
For excavation problems the first integral in \myeqref{Regular2} can therefore be defined as:
\begin{equation}
\label{ }
\int_{\Gamma_{e}}  \fund{U}(\sourcept_{n},\fieldpt) \ \myVec{\dual}^{e}(\fieldpt)\ d\Gamma_{e}(\fieldpt) = \int_{\Gamma_{e}}  \fund{U}(\sourcept_{n},\fieldpt) \mathbf{ t}_{0}^{e} (\fieldpt) d\Gamma_{e}(\fieldpt)   
\end{equation}

We use the \textit{independent field approximation method}  \cite{Marussig2015}, i.e. we leave the definition of the geometry untouched and only enrich the basis functions for the approximation of the unknown. This means that $\hat{ \NURBS_{i} }(\xi, \eta)$ are basis functions that are obtained by refining the basis functions $\NURBS_{i} (\xi, \eta)$ describing the geometry by order elevation or knot insertion.
For the infinite patches we assume plane strain conditions along the direction to infinity, i.e. we assume that the displacements are constant in this direction.

Inserting \myeqref{uapprox} we obtain for the second integral:
\begin{eqnarray}
\label{ }
\int_{\Gamma_{e}}  \fund{T}(\sourcept_{n},\fieldpt) \ \myVec{\primary}^{e}(\fieldpt)\ d\Gamma_{e}(\fieldpt)  
&=& \sum_{k=1}^{K} \left( \quad \int_{\Gamma_{e}}  \fund{T}(\sourcept_{n},\fieldpt)    \hat{\NURBS}_{k} (\xi,\eta) \  d\Gamma_{e}(\fieldpt)  \right) \myVec{\primary}_{k}^{e}
\end{eqnarray}

\subsection{Collocation points}
The collocation points $\sourcept_n$ are the locations where the unknown parameters reside i.e. the anchors of the basis functions (\myeqref{Greville}).
Since the anchor locations are calculated individually for each patch a unique location can only be achieved when the basis functions for describing the unknown (displacements) match at the locations where patches connect. We call this \textit{continuous collocation}.

However, in the BEM it is possible to waive this restriction. In \textit{discontinuous  collocation} the collocation points are moved slightly inside the patch with an offset. For information on how the optimal location is determined see \cite{Marussig2015}.
We use discontinuous collocation when sharp corners are encountered, for example at intersections or at tunnel ends, thereby also avoiding the stress singularity there.

\subsection{System of equations}
An assembly of all patch contributions leads to the following system of equations to be solved:
\begin{equation}
\label{Syseq}
[\mathbf{ L}] \{\mathbf{ u}\}= \{\mathbf{ r}\}
\end{equation}
where $[\mathbf{ L}]$ is the assembled left hand side matrix, $\{\mathbf{ u}\}$ is a vector that contains the parameter values $\mathbf{ u}_{n}$ at collocation points and $\{\mathbf{ r}\}$ is a right hand side, computed from the known tractions obtained by \myeqref{tdef}.

\subsection{Numerical integration of patch integrals}
For each boundary patch $e$ we have to compute the following integrals:
\begin{eqnarray}
\label{ }
 \mathbf{ U}_{n}^{e}= \int_{\Gamma_{e}}  \fund{U}(\sourcept_{n},\fieldpt(\xi,\eta))   \mathbf{ t}_{0} (\xi,\eta) \  d\Gamma_{e}(\fieldpt) \\
\mathbf{ T}_{nk}^{e} = \int_{\Gamma_{e}}  \fund{T}(\sourcept_{n},\fieldpt(\xi,\eta))    \hat{\NURBS}_{k} (\xi,\eta) \  d\Gamma_{e}(\fieldpt) \\
\mathbf{ T}_{n}^{e} = \int_{\Gamma_{e}}  \fund{T}(\sourcept_{n},\fieldpt(\xi,\eta))   \  d\Gamma_{e}(\fieldpt)
\end{eqnarray} 
The evaluation of the patch integrals is more involved than in the FEM since we are dealing with integrals that contain singular fundamental solutions.
Even though Gauss Quadrature is not really suitable we use it for lack of alternatives.
The integration scheme now depends on the location of the collocation point ($\sourcept_n$) since the fundamental solutions $\fund{U}, \fund{T}$ tends to infinity as $\fieldpt$ approaches $\sourcept_n$. If $\sourcept_n$ is outside the patch we use regular integration otherwise we have to use singular integration. The careful monitoring of the the accuracy of the integration is paramount for the quality of the results. Therefore the number of Gauss points has to be increased or the integration region be subdivided if $\sourcept_n$ is close. For the singular integration we use the method of dividing into triangular subregions, where the Jacobian tends to zero as $\sourcept_n$ is approached. For a detailed description of efficient and accurate integration methods  the reader is referred to \cite{BeerMarussig}.

\subsubsection{Refinement demonstration on a NATM tunnel}
\begin{figure}
\begin{center}
\begin{overpic}[scale=0.7]{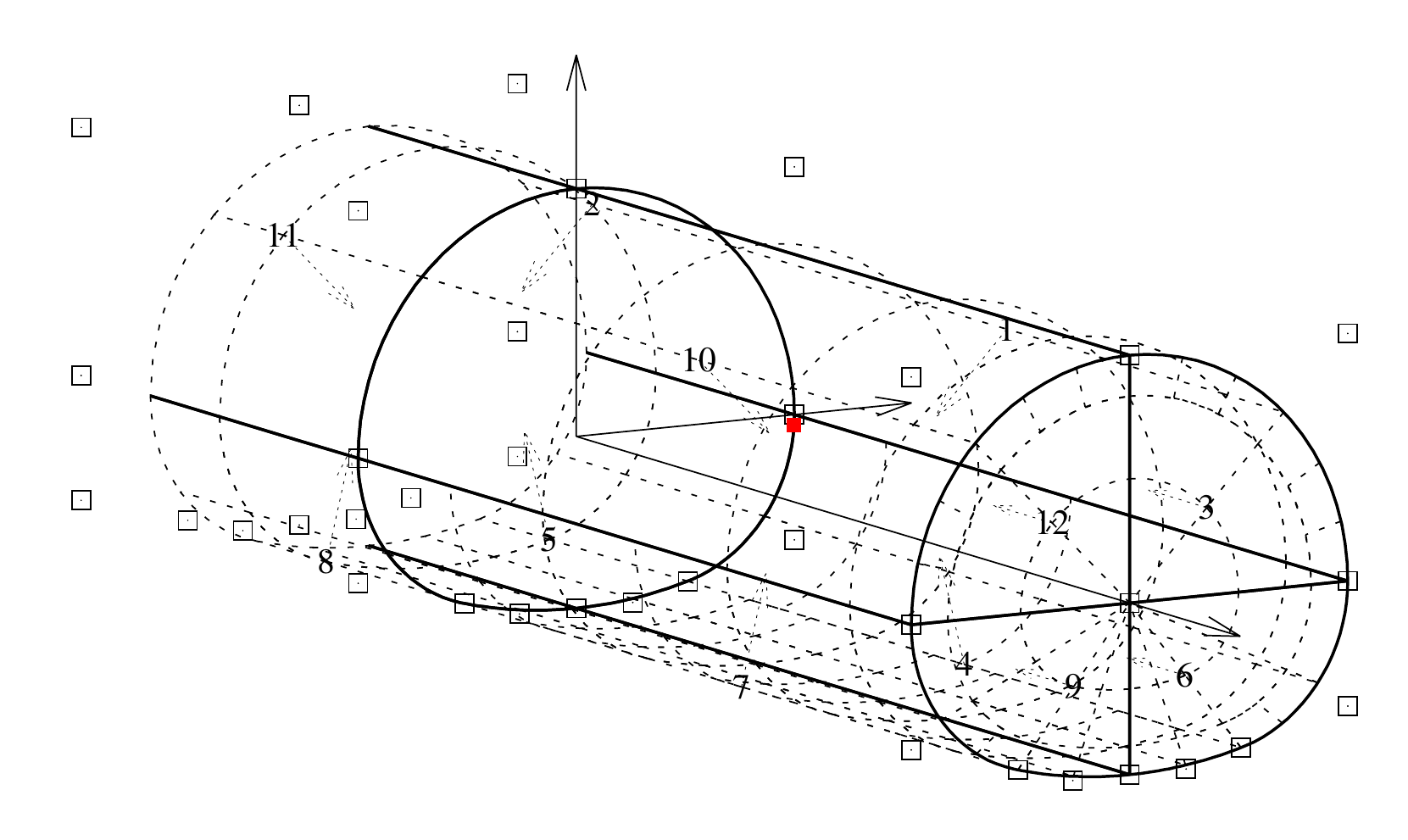}
\put(88,18){x}
\put(65,35){y}
\put(40,62){z}
\end{overpic}
\caption{Definition of semi-infinite NATM tunnel with control points shown as hollow squares. The point where the max. compressive stress is monitored is depicted by a red square.}
\label{NATM3D}
\end{center}
\end{figure}

Here we  demonstrate how refinement can be applied elegantly on the example of a NATM tunnel. The geometry of the semi-infinite tunnel is exactly defined by 4 finite patches, 4 degenerate patches and 4 infinite patches (see \myfigref{NATM3D}).

We investigate various refinement procedures by refining the basis functions that define the geometry along the tunnel (which are linear) and note the effect on the results.
The refinements only have effect on the number of collocation points (and therefore number of degrees of freedom). No re-meshing is necessary. The refinement options are: order elevation, knot insertion and K-refinement (a combination of the two).

\begin{figure}[H]
\begin{center}
\begin{overpic}[scale=0.20]{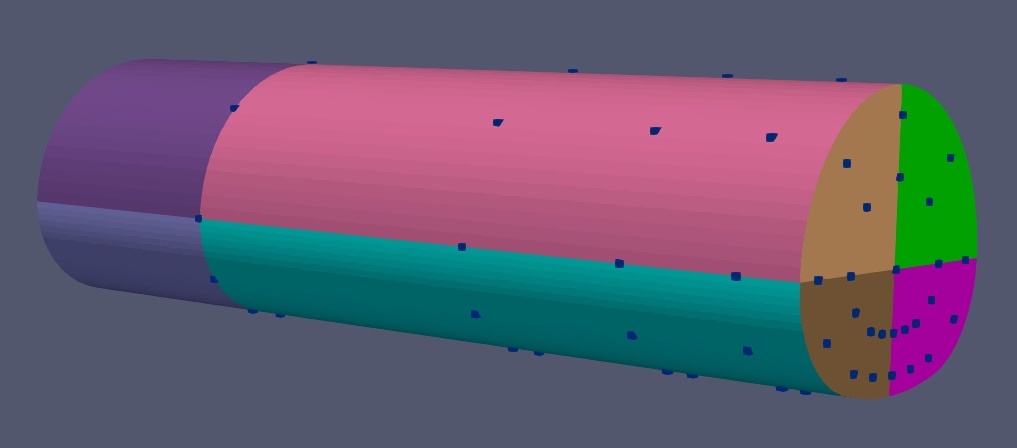}
\put(-30,20){d.o.f.=219}
\put(-30,10){Stress= 3.308}
\end{overpic}
\begin{overpic}[scale=0.5]{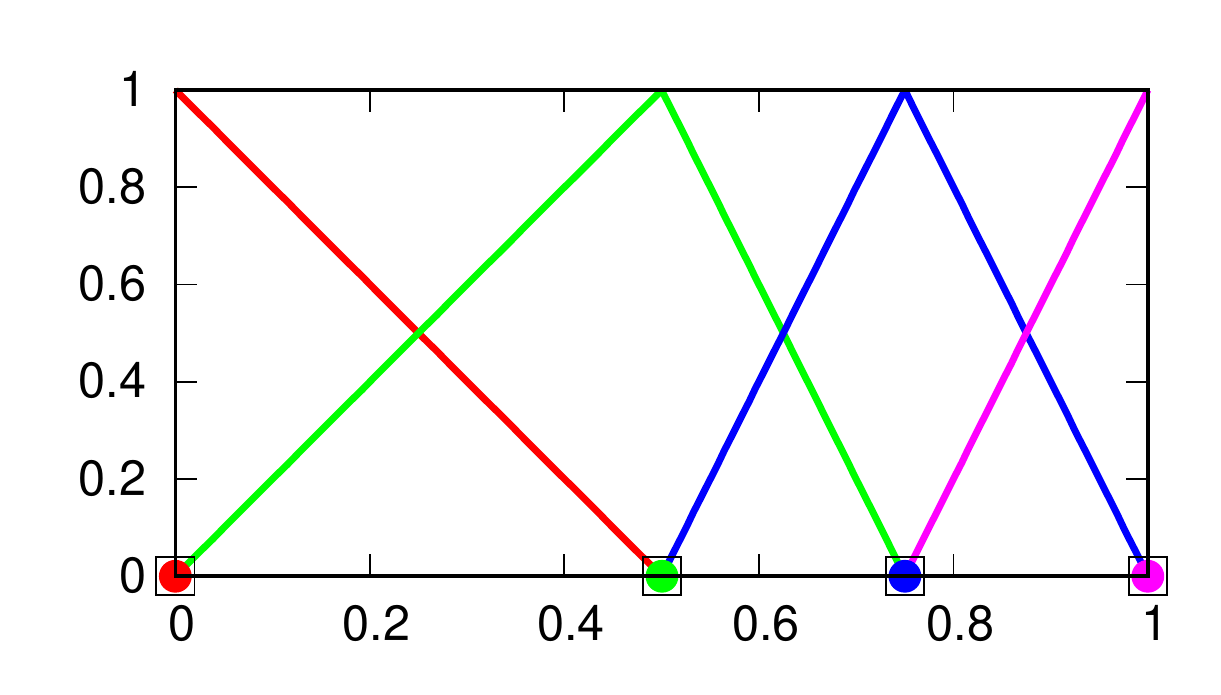}
\end{overpic}
\caption{Refinement of NATM tunnel: Knot insertion at $\eta$= 0.5 and 0.7}
\label{NATM2}
\end{center}
\end{figure}

\begin{figure}[H]
\begin{center}
\begin{overpic}[scale=0.20]{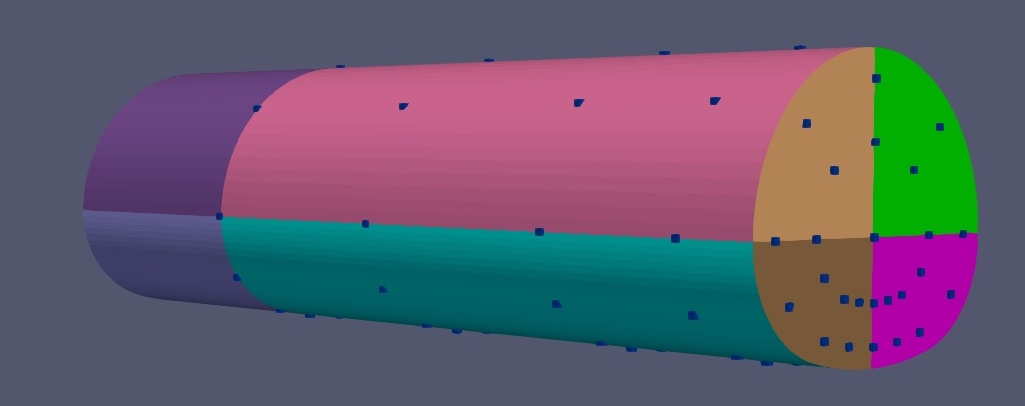}
\put(-30,20){d.o.f.=255}
\put(-30,10){Stress= 3.310}
\end{overpic}
\begin{overpic}[scale=0.5]{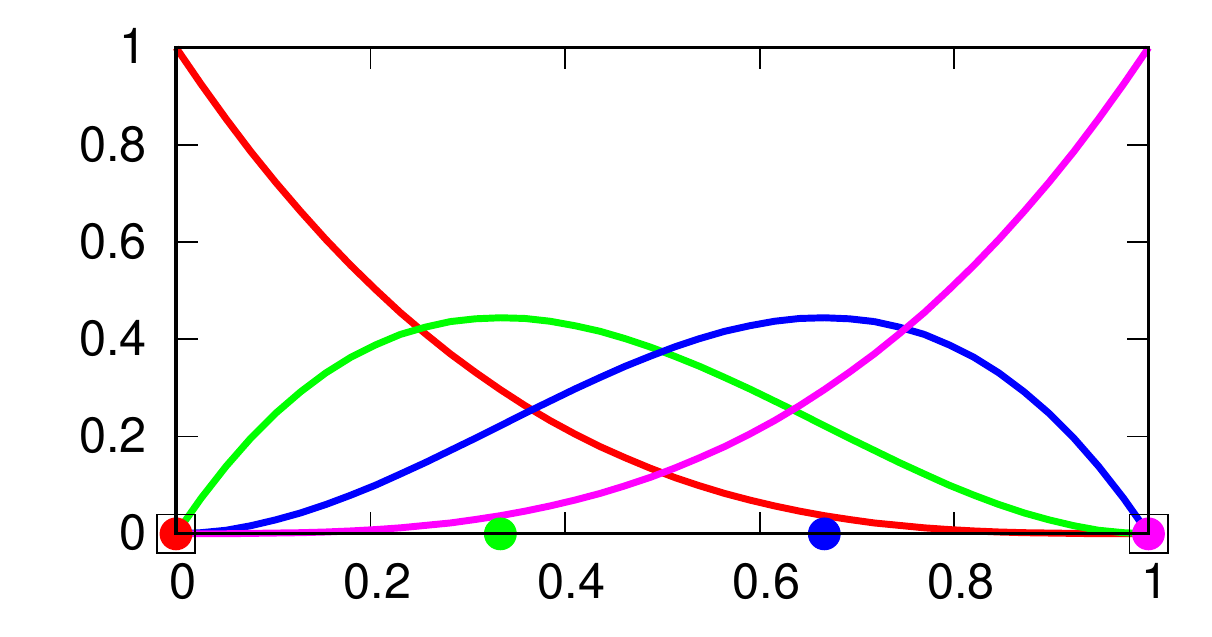}
\end{overpic}
\caption{Refinement of NATM tunnel: Order elevation from linear to cubic}
\label{NATM3}
\end{center}
\end{figure}

\begin{figure}[H]
\begin{center}
\begin{overpic}[scale=0.20]{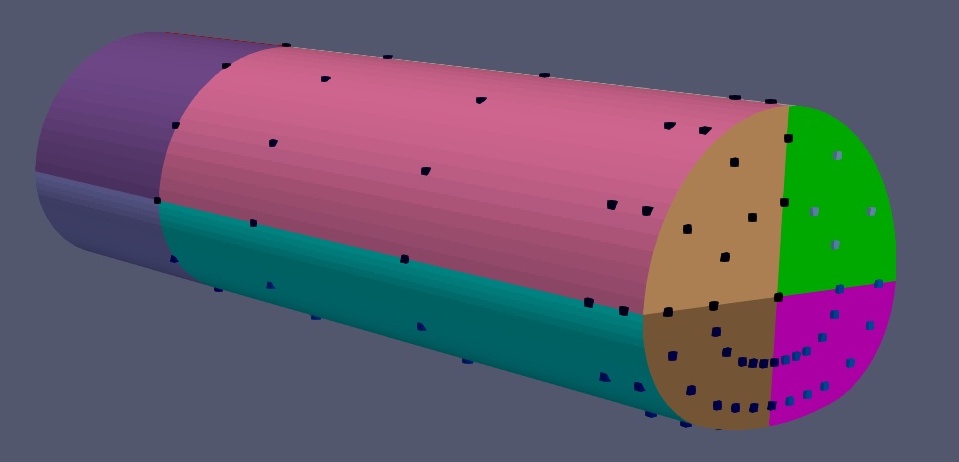}
\put(-35,20){d.o.f.=321}
\put(-35,10){Stress= 3.314}
\end{overpic}
\begin{overpic}[scale=0.5]{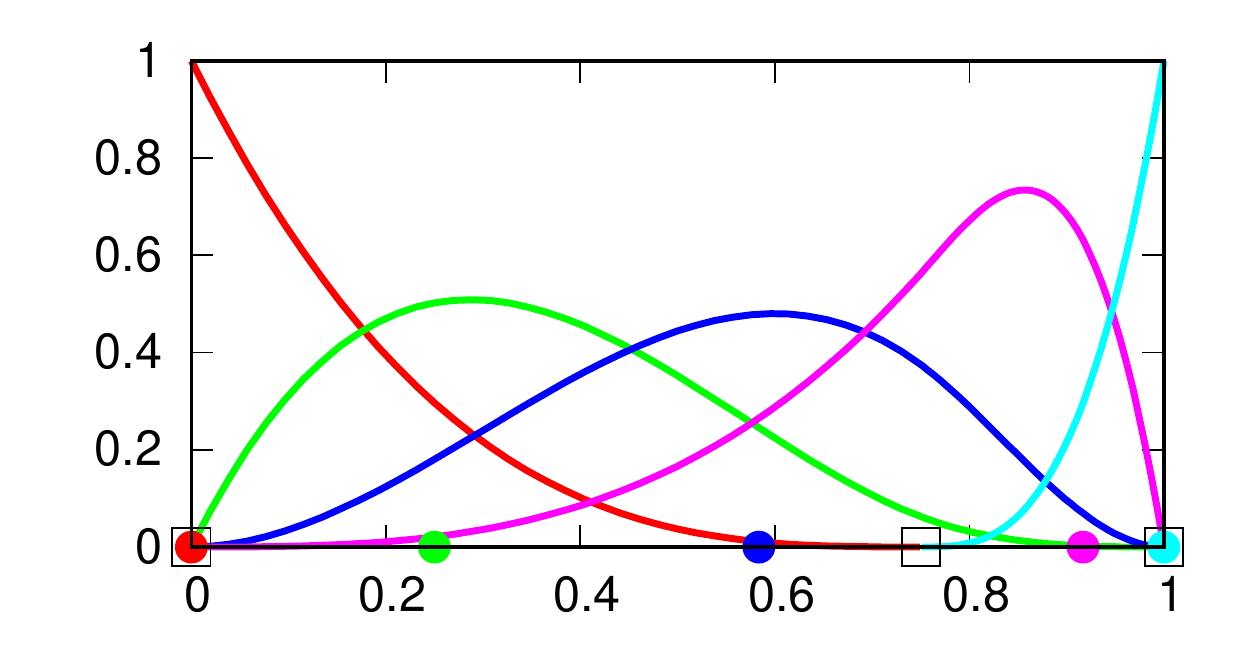}
\end{overpic}
\caption{Refinement of NATM tunnel: Order elevation from linear to cubic and insertion of knot at $\eta$= 0.75}
\label{NATM4}
\end{center}
\end{figure}

At the end of the tunnel we use discontinuous  collocation. This means that the corners are rounded and stress concentrations there do not pollute the results.
For this demonstration we use non-dimensional properties with elastic modulus E=10 and Poisson's ratio $\nu=0$ and a simple loading, namely an excavation in a virgin stress field of 1 (compression) in vertical direction and zero in the horizontal directions. We look at the number of degrees of freedom (d.o.f) and the change in the maximum  stress at the point depicted in \myfigref{NATM3D}. In Figures  \ref{NATM2}, \ref{NATM3} , \ref{NATM4}  we show the change in d.o.f and stress with the refinement of the solution. We also show the basis functions used for the approximation of the unknown displacements in the direction along the tunnel (right) and the resulting collocation points (left).

\subsubsection{Example}
\label{ExBorehole}
\begin{figure}[H]
\begin{center}
\includegraphics[scale=0.20]{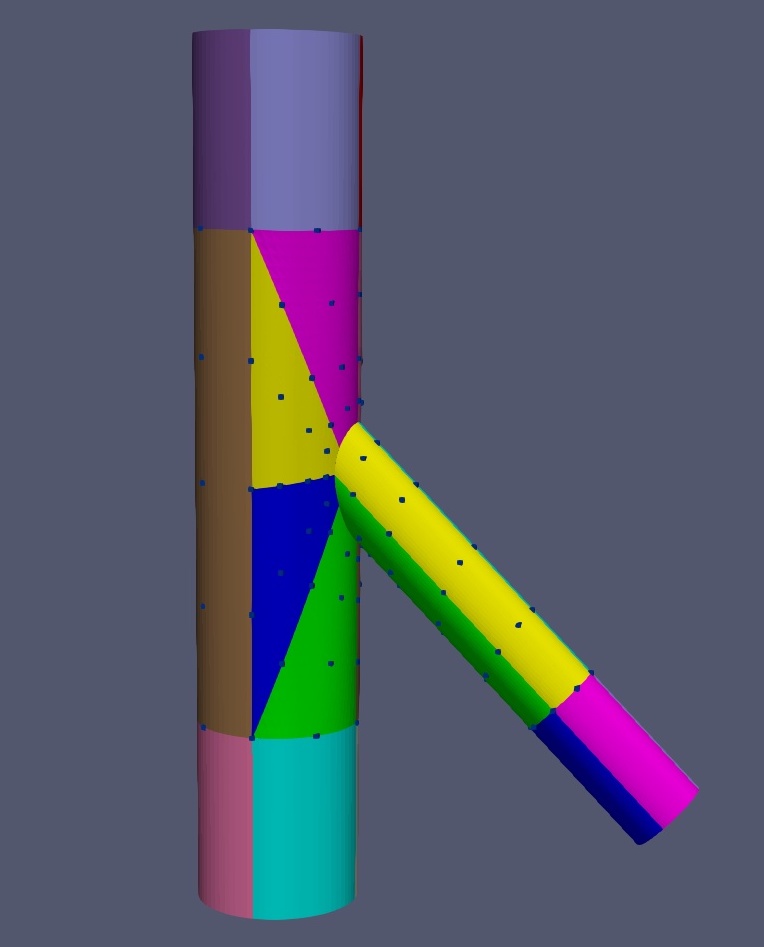}
\includegraphics[scale=0.40]{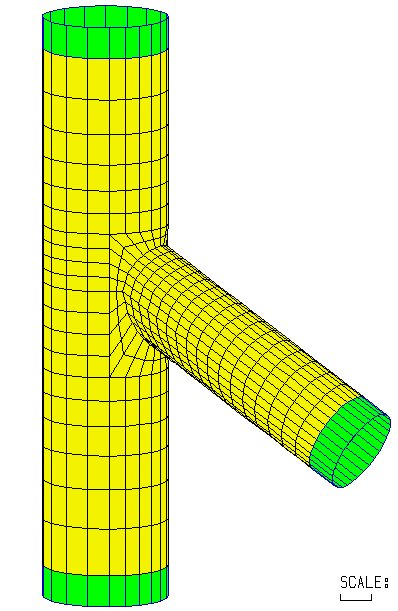}
\caption{Left:  IGABEM discretisation, showing collocation points used for the simulation, Right: Conventional BEM mesh using Serendipity boundary elements}
\label{Bore1}
\end{center}
\end{figure}

The example relates to a real problem in oil extraction. In order to extend the production area, smaller boreholes are sometimes deviated from the main borehole.
In this example we look at such a deviation at 4200 m depth, drilled in sandstone. The main borehole has a diameter of 20 cm, the deviated borehole 10 cm.
The input data are as follows:
\begin{itemize}
  \item Stress regime: $\sigma_\mathrm{v}$ = 130.5 MPa,  $\sigma_\mathrm{h}$= 121.5 MPa, $\sigma_\mathrm{H}$= 112.5 MPa
  \item Well pressure= 72 MPa
  \item Elastic properties: E= 10 GPa, $\nu$= 0.25
\end{itemize}
The stress-subscripts v, h, H mean vertical, horizontal towards the deviation, horizontal normal to it.
The IGABEM discretisation of the boreholes is shown in \myfigref{Bore1} on the left and consists of 8 finite patches, 8 infinite patches and 8 trimmed patches. Special algorithms published  in \cite{Beer2022S} were used to determine the trimmed patches near the intersection. The collocation points shown are a result of K-refinment of the basis functions that describe the geometry. At the intersection discontinuous collocation was used. 
\begin{figure}[H]
\begin{center}
\includegraphics[scale=0.20]{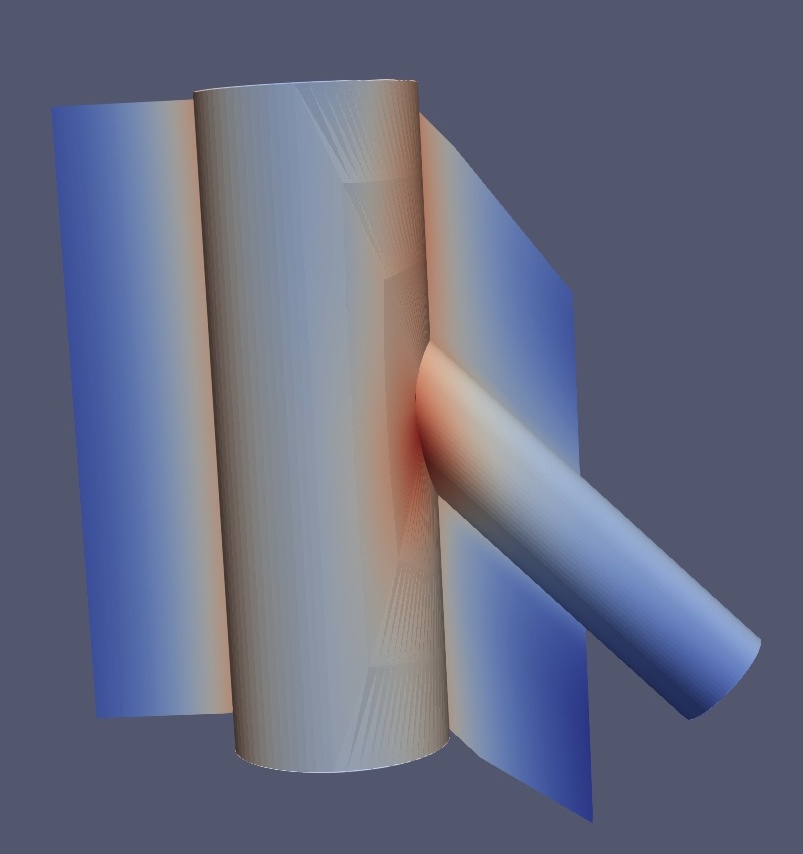}
\caption{Contours of displacements}
\label{BoreD}
\end{center}
\end{figure}

The IGABEM simulation has 357 degrees of freedom. We compare this with the mesh required for a BEM simulation using conventional boundary elements (\myfigref{Bore1} right) which has over 2500 degrees of freedom. It should be noted that the IGABEM discretisation exactly defines the geometry of the two cylinders, only approximating the geometry at the intersection, whereas in the conventional BEM both the geometry of the cylinder and the intersection is approximated.
\myfigref{BoreD} shows one result of the simulation, namely the contours of displacement. It can be seen that despite the fact that collocation is discontinuous at the intersection, no jump in displacements can be observed. This has also been observed in other published simulations using discontinuous collocation (see for example \cite{Marussig2015}).
We will revisit the example later when we introduce plasticity.

\section{The boundary element method with volume effects}
So far we have only considered effects on the boundary and neglected effects that occur inside the volume. This means the simulation can only deal with homogeneous and elastic domains.  We now explain how this restriction can be lifted.

\subsection{Integral equations with volume effects}
To account for inhomogeneous and inelastic behaviour we use the method of initial stress, well known in the FEM community for dealing with plasticity.
The \myeqref{Inteq} has to be amended by adding a volume integral over a specified domain $\domain_0$ (for a derivation see \cite{BeerMarussig}):
\begin{equation}
\label{InteqV}
    \begin{aligned}
\int_{\boundary} \fund{T}(\sourcept_{n},\fieldpt) ( \myVec{u}(\fieldpt) - \myVec{u}(\sourcept_{n})) \ d\boundary(\fieldpt) - \mathbf{ A}_{n} \myVec{u}(\sourcept_{n}) &=& \int_{\boundary} \fund{U}(\sourcept_{n},\fieldpt) \ \myVec{t}(\fieldpt) \ d\boundary(\fieldpt) \\
+  \int_{\domain_{0}} \fund{E} (\sourcept_{n},\fieldptv)  \myVecGreek{\sigma}_{0} (\fieldptv) d \domain_{0} (\fieldptv)  .
\end{aligned}
\end{equation}
In the above  $\fund{E} (\sourcept_{n},\fieldptv)$ is a matrix containing the fundamental solutions for the strain at point $\fieldptv$ due to a source at $\sourcept_n$.  
 $ \myVecGreek{\sigma}_{0} (\fieldptv)$ is the initial stress at point $\fieldptv$ inside $\domain_0$ in Voigt notation:
\begin{equation}
\label{Voightnot}
\myVecGreek{\sigma}_{0}= \left \{\begin{array}{c}\sigma_{x0} \\\sigma_{y0} \\\sigma_{z0} \\\tau_{xy0} \\\tau_{yz0}\\\tau_{xz0}\end{array}\right\}
\end{equation}
The fundamental solution $\fund{E}$ is in index notation:
\begin{equation}
\label{ }
E_{ijk}= \frac{-C}{r^{2}}\left[C_{3}(r_{,k}\delta_{ij} + r_{,j}\delta_{ik}) - r_{,i}\delta_{jk} + C_{4} \ r_{,i} r_{,j} r_{,k}\right]
\end{equation}
where the index $i$ specifies the direction of the source and $jk$ the strain component. The constants are: $C=\frac{1}{16 \pi G (1- \nu})$, $C_{3}=1 - 2\nu$ and $C_{4}=3$. 

To arrive at matrix $\fund{E}$ we convert the tensor $E_{ijk}$ to a matrix:
\begin{equation}
\label{eq_Voigt_convertion}
\fund{E}= \left[\begin{array}{cccccc}E_{111} & E_{122} & E_{133}  & E_{112}  + E_{121}& E_{123} + E_{132} & E_{113} +E_{131}  \\E_{211} & E_{222} & E_{233}  & E_{212} + E_{221}  & E_{223}  + E_{232} & E_{213} + E_{231}\\E_{311} & E_{322} & E_{333}  & E_{312} + E_{321} & E_{323}  + E_{332}& E_{313} + E_{331}\end{array}\right]
\end{equation}

For the volume discretisation we generate a grid of points inside $\domain_0$.
We gather the initial stresses at grid points in a long vector $\{ \myVecGreek{\sigma}_{0}\} $ and expand \myeqref{Syseq} to include volume effects:
\begin{equation}
\label{DisIE}
[\mathbf{ L}] \{\mathbf{u}\} = \{\mathbf{ r}\} + [\mathbf{ B}_{0}] \{ \myVecGreek{\sigma}_{0}\} 
\end{equation} 

The sub-matrices of matrix $ [\mathbf{ B}_{0}]$, related to collocation point $n$ and grid point $i$ are given by
\begin{equation}
\label{SubB}
\mathbf{ B}_{0ni} =\ \int_{\domain_{0}}  \fund{E} (\sourcept_{n},\fieldptv) N_{i}(\fieldptv) d \domain_{0} (\fieldptv) 
\end{equation}
where  $N_{i}(\fieldptv)$ are interpolation functions which will be introduced later.

\subsection{Computation of the initial stress at grid points}

Consider the elastic stresses $\myVecGreek{\sigma}$ due to strains $\myVecGreek{\epsilon} $:
\begin{equation}
\label{ }
 \myVecGreek{\sigma}= \mathbf{ D} \myVecGreek{\epsilon} 
\end{equation} 
where $ \mathbf{ D}$ is the elasticity matrix and
\begin{equation}
\label{Voightnot}
\myVecGreek{\epsilon}= \left \{\begin{array}{c}\epsilon_{x} \\\epsilon_{y} \\\epsilon_{z} \\\gamma_{xy} \\\gamma_{yz}\\\gamma_{xz}\end{array}\right\}
\end{equation}

The initial stresses $ \myVecGreek{\sigma}_{0} $ may arise due to a difference between the elastic properties that are used to compute the fundamental solutions and the properties of $\domain_0$ (referred as inclusion here):
\begin{equation}
\label{InitialS}
\myVecGreek{\sigma}_{0} =( \mathbf{ D} - \mathbf{ D}_{e,incl}) \myVecGreek{\epsilon} =\mathbf{ D}^{\prime}\myVecGreek{\epsilon} 
\end{equation}
where $ \mathbf{ D}_{e,incl}$ is the elasticity matrix for the inclusion.

The initial stresses $ \myVecGreek{\sigma}_{0} $ due to inelastic behaviour are given by:
\begin{equation}
\label{InitialSp}
\myVecGreek{\sigma}_{0} =( \mathbf{ D} - \mathbf{ D}_{ep,incl}) \myVecGreek{\epsilon} =\mathbf{ D}^{\prime}\myVecGreek{\epsilon} 
\end{equation}
where $ \mathbf{ D}_{ep,incl}$ is the elasto-plastic constitutive matrix of the inclusion.

\subsection{Computation of displacements at grid points}

To obtain the initial stresses we need to compute the strains at grid points inside the inclusion. These are obtained as derivatives of the displacements.
The displacement vector at a grid point with the coordinates $\pt{x}_{g}$ is given by:
\begin{equation}
    \begin{aligned}
        \label{eps0}
       \mathbf{ u}(\pt{x_{g}}) &= \int_{\boundary} \left[ \fund{U}(\pt{x}_g,\fieldpt) \   \myVec{\dual} (\fieldpt)  - \fund{T}(\pt{x}_g,\fieldpt) \   \myVec{\primary} (\fieldpt) \right] d\boundary (\fieldpt) \\
%         \nonumber
        &+ \int_{\domain_{0}} \fund{E} (\pt{x}_g,\fieldptv) \myVecGreek{\sigma}_{0} (\fieldptv)  d \domain_{0} (\fieldptv)   
    \end{aligned}
\end{equation}

Displacement vectors at all grid points can be gathered in long vector $\{\mathbf{ u} \}_g$, leading to the following discretised version of equation \ref{eps0}:
\begin{equation}
\label{ }
\{\myVec{u}\}_g= [\hat{\mathbf{ A}}] \{\mathbf{u}\} + \{\bar{\mathbf{ c}}\} + [\bar{\mathbf{ B}}_{0}]\{ \myVecGreek{\sigma}_{0}\}   
\end{equation}
where $[\hat{\mathbf{ A}}]$ is an assembled matrix that multiplies with the unknown $\{\mathbf{u}\}$ and $\{\bar{\mathbf{ c}}\}$ collects the displacements due to given BC's. $ [\bar{\mathbf{ B}}_{0}]$ is similar to $ [\mathbf{ B}_{0}]$ in \myeqref{SubB} except that the grid point coordinates $\pt{x}_{g}$ replace the source point coordinates $\sourcept_{n}$.

Because of the singularity of $\fund{T}$ the displacements can not be computed on the boundary. For a point on a patch boundary with the coordinates $\pt{x}_{k}$ we compute the displacement by:
\begin{equation}
\label{ }
\mathbf{ u}(\pt{x}_{k})=\sum_{i}^{I} R_{i}^{e}(\myVecGreek{\xi}_{k}) \mathbf{ u}_{i}^{e}
\end{equation}
where $R_{i}^{e}(\myVecGreek{\xi})$ are the NURBS basis functions used for approximating the displacements in patch $e$. The superscript $e$ indicates the patch that contains the point $\pt{x}_{k}$ and $\myVecGreek{\xi}_{k}$ are its local coordinates.
The matrix $[\hat{\mathbf{ A}}]$ and the vector $\{\bar{\mathbf{ c}}\} $ have to be modified for these grid points.
 
\subsection{Evaluation of the volume integral}
As mentioned in the introduction, the topic of volume integration has been addressed early in the development of the BEM. Integration via volume cells was proposed, but was not very attractive because of the requirement of generating a mesh, similar to the FEM. So called mesh-less methods were proposed such as the dual reciprocity method that use radial basis functions \cite{Nintcheu}. However, these are not really mesh-less, since points inside the domain had to be specified. Most importantly they do not work for infinite domain problems.
Recently it was proposed to use NURBS-based supercells to define the volume \cite{Beer17}. This works well for simple volume geometries, but is not suitable for more complicated geometries such as in the example presented later.

Here we follow the traditional BEM approach with the main difference that we propose that the cell meshes can be generated completely automatic with minimum user effort. 
We associate the grid points, where the initial stresses $\{ \myVecGreek{\sigma}_{0}\}$ are stored, with nodes of cells.
For the cell geometry we use brick elements used by the FEM community.  Here we discuss the volume integration first before proceeding to the automatic generation of cell meshes, using immersed technology.
\begin{figure}[H]
\begin{center}
\begin{overpic}[scale=0.6]{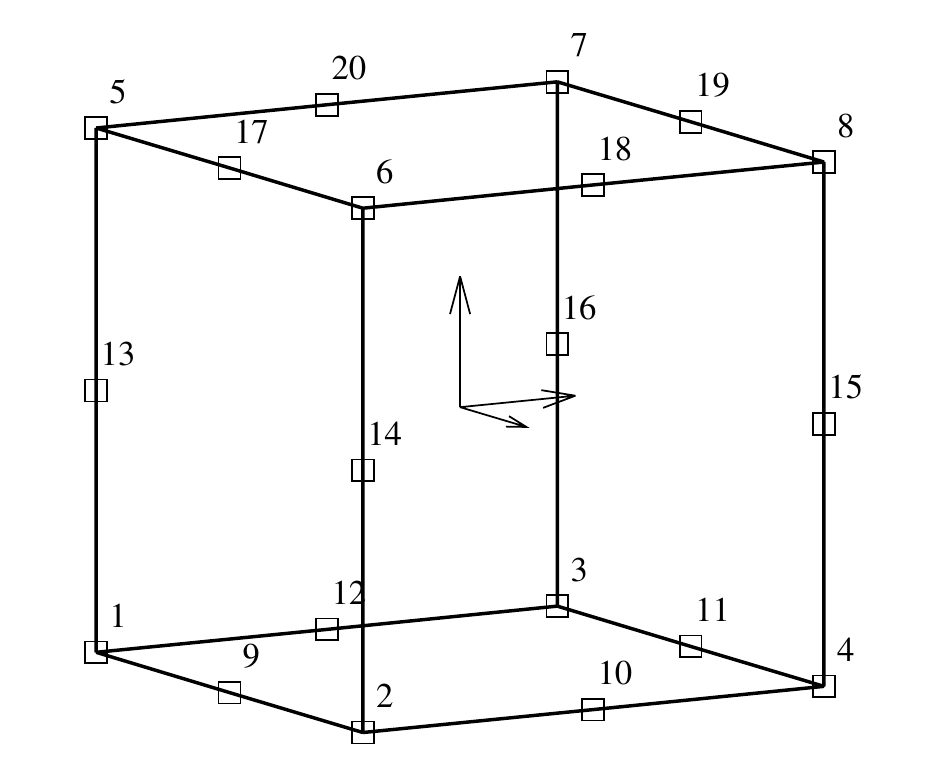}
\put(57,35){$\xi$}
\put(65,42){$\eta$}
\put(50,58){$\zeta$}
\end{overpic}
\caption{A cell in the local $\xi,\eta,\zeta$ coordinate system.}
\label{Cells}
\end{center}
\end{figure}
The geometry of a cell $c$ is defined by:
\begin{equation}
\label{cellgeo}
 \pt{x}^{c}(\myVecGreek{\xi})= \sum_{j=1}^{J} N_{j}(\myVecGreek{\xi})   \pt{x}^{c}_{j}
\end{equation}
where $N_{j}(\myVecGreek{\xi})$ are linear or quadratic Serendipity shape functions of the local coordinate $\myVecGreek{\xi}=(\xi,\eta,\zeta)^{\mathrm{T}}=[-1,1]^3$ and $\pt{x}^{e}_{j}$ are nodal coordinates.  The mid-side nodes can be active or inactive. If the mid-side node is active a quadratic interpolation is assumed along the edge, if inactive there is a linear interpolation.

The same shape functions are used to interpolate the initial stress values between nodal points (\myeqref{SubB}).
The initial stresses at cell nodes are computed from the strains using \myeqref{InitialS} or \myeqref{InitialSp}. 

To obtain the strains at cell nodes we first interpolate the displacements inside a cell:
\begin{equation}
\label{cellinter}
\mathbf{ u}^{c}(\myVecGreek{\xi})= \sum_{j=1}^{J}   N_{j}(\myVecGreek{\xi})  \mathbf{ u}^{c}_{j}
\end{equation}
where $\mathbf{ u}^{c}_{j}$ is the displacement vector at nodal point $j$.

The global derivatives of the displacements are:
\begin{equation}
\label{cellinter}
\frac{\partial \mathbf{ u}^{c}(\myVecGreek{\xi})}{\partial \mathbf{ x}}= \sum_{j=1}^{J} \frac{\partial N_{j}(\myVecGreek{\xi})}{\partial \mathbf{ x}}  \mathbf{ u}^{c}_{j}
\end{equation}
The global derivatives of the shape functions in terms of the local derivatives given by:
\begin{equation}
\label{ }
\frac{\partial N_{j}(\myVecGreek{\xi})}{\partial \mathbf{ x}}= \mathbf{ J}^{-1}  \frac{\partial N_{j}(\myVecGreek{\xi})}{\partial \myVecGreek{\xi}}
\end{equation}
with the Jacobian matrix:
\begin{equation}
\label{ }
\mathbf{ J}= \left[\begin{array}{ccc}\frac{\partial x}{\partial \xi} & \frac{\partial y}{\partial \xi} & \frac{\partial z}{\partial \xi} \\ \\ \frac{\partial x}{\partial \eta} & \frac{\partial y}{\partial \eta}& \frac{\partial z}{\partial \eta} \\ \\ \frac{\partial x}{\partial \zeta} & \frac{\partial y}{\partial \zeta}& \frac{\partial z}{\partial \zeta}\end{array}\right]
\end{equation}
The Jacobian is:
\begin{equation}
\label{ }
J= |\mathbf{ J} |
\end{equation}

In Voight notation we define a strain (pseudo)-vector:
\begin{equation}
\label{ }
 \myVecGreek{\epsilon}= \left(\begin{array}{c} \frac{\partial u_{x}}{\partial x} \\ \\ \frac{\partial u_{y}}{\partial y} \\ \\ \frac{\partial u_{z}}{\partial z} \\ \\ \frac{\partial u_{x}}{\partial y} + \frac{\partial u_{y}}{\partial x} \\ \\ \frac{\partial u_{y}}{\partial z} + \frac{\partial u_{z}}{\partial y} \\ \\ \frac{\partial u_{x}}{\partial z} + \frac{\partial u_{z}}{\partial x} \end{array}\right)
\end{equation}

The strain vectors at all nodal points can be stored into a long vector $\{ \myVecGreek{\epsilon} \}$ and expressed as:
\begin{equation}
\label{ }
\{ \myVecGreek{\epsilon} \}= [\mathbf{ \hat{B}}] \{ \mathbf{ u} \}_g
\end{equation}
where $ [\mathbf{ \hat{B}}]$ is an assembled matrix containing derivatives of shape functions at cell nodes and $ \{ \mathbf{ u} \}_g$ contains the displacement vectors at all cell nodes.

\subsection{Numerical volume integration with cells}
Here we explain how the volume integration is carried out to arrive at matrix $[\mathbf{ B}_0]$ .
The sub-matrix $\mathbf{ B}_{0ni}$ related to collocation point $n$ and grid point $i$ can be expressed as:
\begin{eqnarray}
\label{ }
\mathbf{ B}_{0ni} &=& \sum_{c=1}^{n_c} \mathbf{ B}^{c}_{0ni}  
\end{eqnarray}
 where $c$ is the cell number, $n_c$ is the number of cells and $i$ is the number of the grid point that stores the initial stress that $\fund{E}$ will be multiplied with.

The sub-matrices of $\mathbf{ B}_0$ are computed by:
\begin{equation}
\label{ }
\mathbf{ B}^{c}_{0ni} = \int_{\domain_{e}}  \fund{E} (\sourcept_{n},\fieldptv) N_{k(i)}(\myVecGreek{\xi}(\fieldptv)) d \domain_{e} (\fieldptv)
\end{equation}
where $k(i)$ is the  local node number of cell $e$ associated with grid point $i$ and $(\myVecGreek{\xi}(\fieldptv))$ is the local coordinate associated with the global coordinate $\fieldptv$.
If $\fieldptv$ is outside the cell $N_{k(i)}$ is zero.

Since the fundamental solution $\fund{E}$ tends to infinity with $o(r^{2})$ - where $r$ is the distance between $\sourcept_{n}$ and $\fieldptv$ - care has to be taken with the integration. 
Gauss Quadrature is not really suited for the evaluation of singular integrals but in absence of alternatives it is used widely in the BEM community.
However, if the number of Gauss points is held constant, the precision of integration deteriorates rapidly with the proximity of the singular point, which is defined as the value of distance $R$ to the singular point, relative to the size of the cell ($L$). If $\sourcept_n$ is on one of the surfaces of the cell the integrand is singular and special integration, referred to as \textit{singular integration}, is applied.

Depending on the proximity of $\sourcept_n$ we need to increase the number of Gauss points to achieve adequate precision of the integration. To estimate the number of the required Gauss points, limiting values of $R/L$ for a given number of Gauss points are presented in \cite{BeerMarussig}. In the current implementation we limit the maximum number of Gauss points, that can be used, to 5. If this number is not sufficient, we change from \textit{regular} to \textit{singular} integration.

\newpage

\subsubsection{Regular integration} 
If $\sourcept_n$ is far enough away and the specified maximum number of Gauss points is adequate to achieve the required precision we can use regular integration.

Using Gauss Quadrature to get an approximation of $\mathbf{ B}_{0ni}^{c}$ we obtain:
\begin{equation}
\label{Regular}
\mathbf{ B}^{c}_{0ni} \approx  \sum^{I}_{i=1} \sum^{J}_{j=1} \sum^{K}_{k=1} \ \fund{E} (\sourcept_{n},\fieldpt (\xi_{i},\eta_{j},\zeta_{k})) N_{j(i)}(\xi_{i},\eta_{j},\zeta_{k})) W_{i} W_{j} W_{k} J(\xi_{i},\eta_{j},\zeta_{k}))
\end{equation}
where $I,J,K$ denotes the number of Gauss points in $\xi,\eta,\zeta$- directions chosen depending on the values of $R/L_{\xi},R/L_{\eta},R/L_{\zeta}$, where $R$ is the minimum distance from $\sourcept_n$ to the cell and  $L_{\xi},L_{\eta},L_{\zeta}$ represent the sizes of the cell in $\xi,\eta,\zeta$ directions.
$ W_{i}, W_{j}, W_{k}$ are Gauss weights and $J(\xi_{i},\eta_{j},\zeta_{k})$ is the Jacobian at the Gauss point with the coordinates $\xi_{i},\eta_{j},\zeta_{k}$.

\subsubsection{Singular integration} 
This integration scheme is for the case that $\sourcept$ is very near or on one of the surfaces of a cell.
The idea is to subdivide the cell into tetrahedral subregions and use them for the integration. The aim is for the Jacobian of the subregion to approach zero as the singular point is approached, thereby cancelling out the singularity.
This is explained in \myfigref{Tetra}.

First we compute the local coordinates of the singular point on the cell ($\xi_{s}, \eta_{s}, \zeta_{s}$.)
Next we define the Gauss point coordinates in the $\myVecGreek{\bar{\xi}}=(\bar{\xi},\bar{\eta},\bar{\zeta})^{\mathrm{T}}=[-1,1]^3$ coordinate system of a tetrahedral subregion. We then map to the cell coordinate system $\myVecGreek{\xi}$ using:
\begin{equation}
\label{map}
\myVecGreek{\xi} (\myVecGreek{\bar{\xi}})= \sum_{j=1}^{J} \breve{N}_{j}(\myVecGreek{\bar{\xi}}) \myVecGreek{\xi} _{j}
\end{equation}
where $\breve{N}_{j}$ are linear shape functions and $\myVecGreek{\xi} _{j}$ are node coordinates defining the tetrahedral subregion.
For the example in  \myfigref{Tetra} these coordinates are defined in table \ref{tetcor}.
\begin{figure}
\begin{center}
\begin{overpic}[scale=0.5]{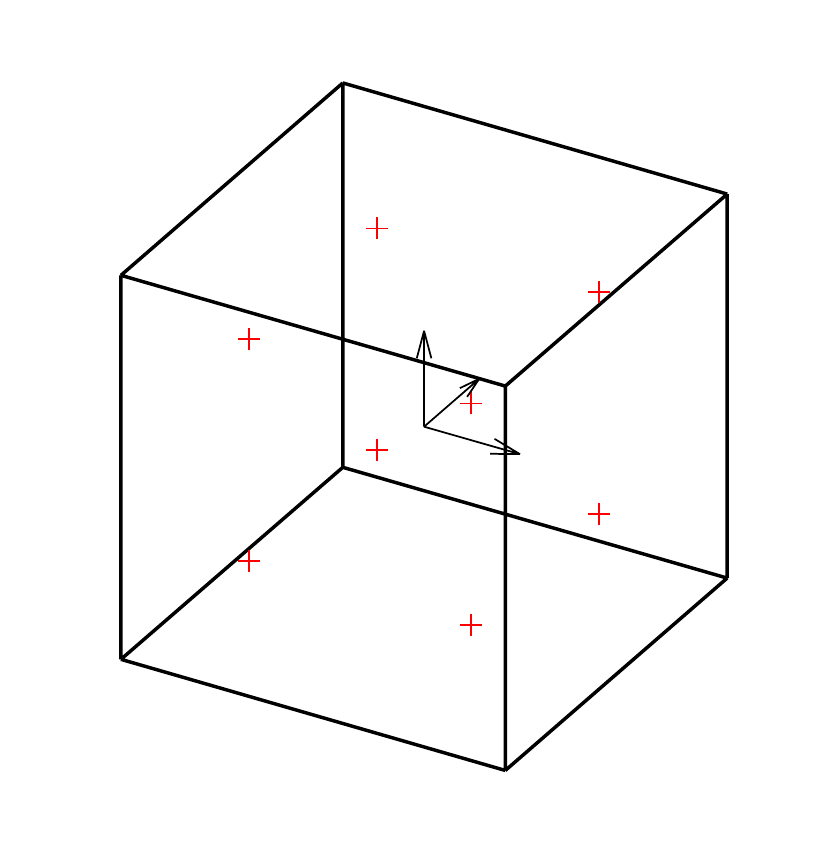}
\put(63,45){$\bar{\xi}$}
\put(60,55){$\bar{\eta}$}
\put(50,65){$\bar{\zeta}$}
\put(50,0){a)}
\end{overpic}
\begin{overpic}[scale=0.5]{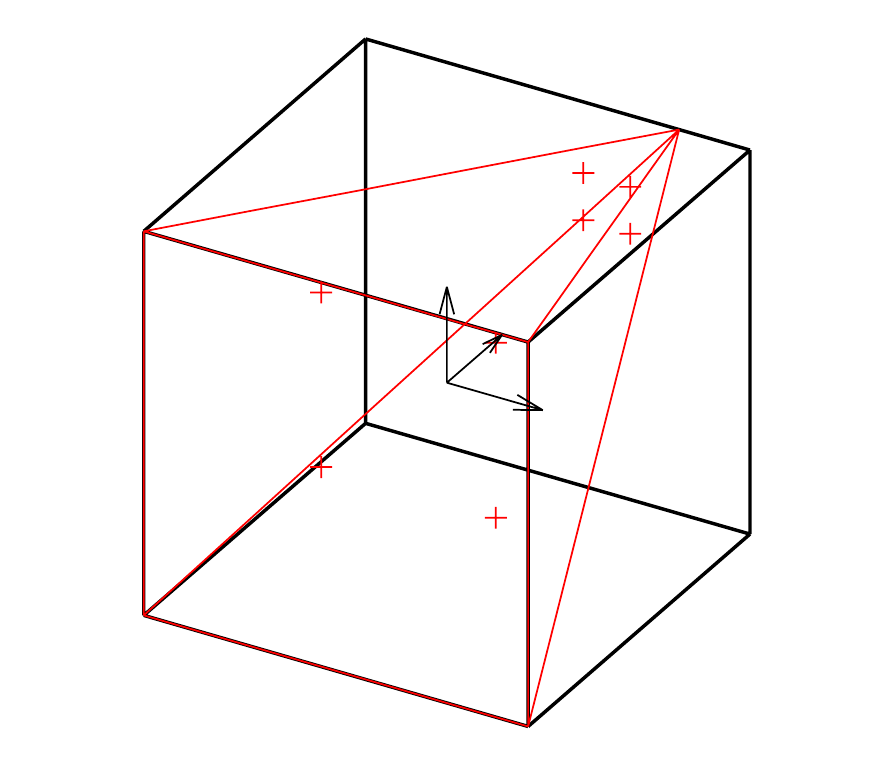}
\put(65,40){$\xi$}
\put(60,50){$\eta$}
\put(50,60){$\zeta$}
\put(75,75){$\sourcept$}
\put(50,0){b)}
\end{overpic}
\begin{overpic}[scale=0.4]{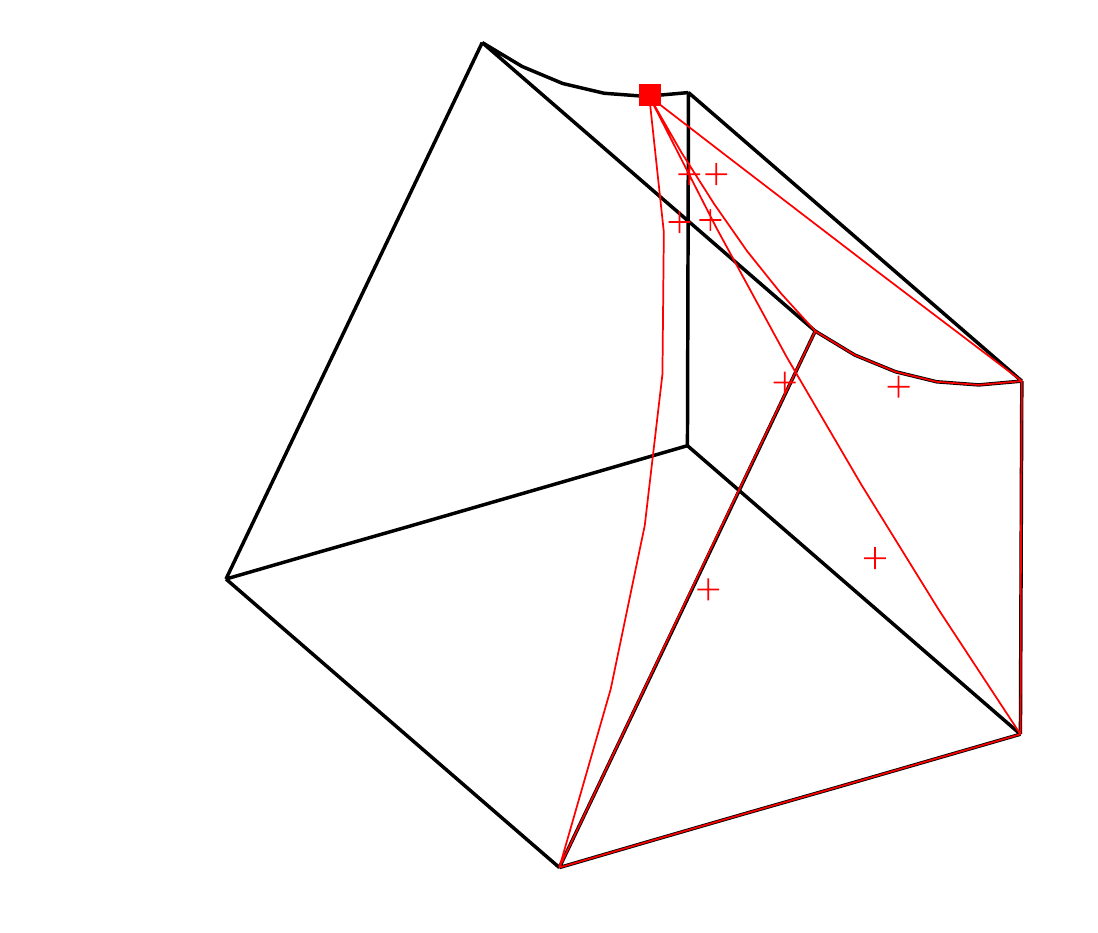}
\put(50,0){c)}
\put(55,80){$\sourcept$}
\end{overpic}
\caption{Mapping of Gauss points to a tetrahedral subregion: Gauss points in a) $\bar{\xi}, \bar{\eta}, \bar{\zeta}$-space b) $\xi, \eta, \zeta$-space c) $x,y,z$- space. $\sourcept$ is the singular point.}
\label{Tetra}
\end{center}
\end{figure}
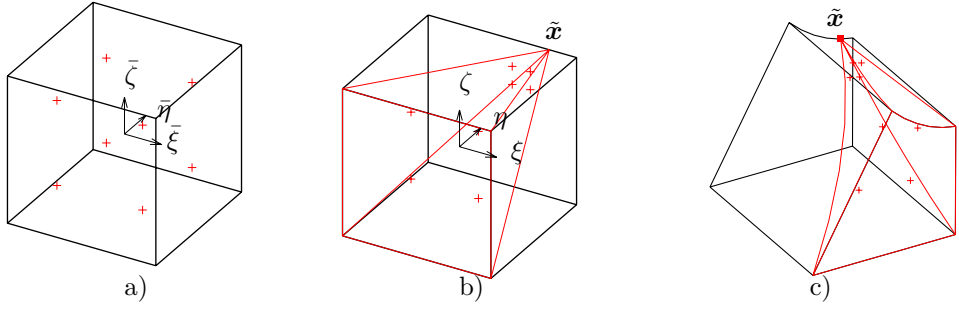
\begin{mytable}
  {H}               %% table position  
  {Coordinates of tetrahedra points.}  %% caption
  {tetcor}  %% label
  {ccccccccc}         %% column layout
  \mytableheader{ i= &1 & 2 & 3 & 4 & 5  & 6 & 7 & 8}  %% header
  %% table content
$\xi$ & \num{-1} & \num{1} & \num{-1} & \num{1} & $\xi_{s}$ & $\xi_{s}$ &  $\xi_{s}$ & $ \xi_{s}$ \\
$\eta$ & \num{-1} & \num{-1} & \num{-1} & \num{-1} & $\eta_{s}$ & $\eta_{s}$ &  $\eta_{s}$ & $ \eta_{s}$ \\
$\zeta$ & \num{-1} & \num{-1} & \num{1} & \num{1} & $\zeta_{s}$ & $\zeta_{s}$ &  $\zeta_{s}$ & $ \zeta_{s}$ \\
\end{mytable}
Because the last 4 node coordinates are the same, the Jacobian of mapping (\ref{map}), $J_{t}$,  tends to zero  with $o(r^{2})$ as the singular point is approached, therefore cancelling out the singularity.

Equation \ref{Regular} is changed to:
\begin{equation}
\label{Singular}
\mathbf{ B}^{c}_{0ni} \approx  \sum^{N_{t}}_{n_t=1} \sum^{I}_{i=1} \sum^{J}_{j=1} \sum^{K}_{k=1} \ \fund{E} (\sourcept_{n},\fieldpt (\xi_{i},\eta_{j},\zeta_{k})) N_{j(i)}(\xi_{i},\eta_{j},\zeta_{k})) W_{i} W_{j} W_{k} J(\xi_{i},\eta_{j},\zeta_{k})) J_{t}
\end{equation}
where $N_{t}$ is the number of tetrahedra, which depends on the location of the singular point (3 if it is on a corner, 4 if it is on an edge as in the example, 5 otherwise).
The computation of $[\bar{\mathbf{ B}}_{0}]$ follows exactly the same procedure, except that the collocation point $n$ is replaced by grid point $i$.

\section{Solution procedure}
The nonlinear solution procedure is explained in detail in \cite{Beer2018a} and involves the solution of:
\begin{equation}
\label{onestep}
[\mathbf{ L}]^{\prime} \{\mathbf{ u}\} = \{\mathbf{ r}\}^{\prime} 
\end{equation}
where
\begin{eqnarray}
[\mathbf{L}]^{\prime} & = & ( [\mathbf{L}] - [\mathbf{B}_0]) ([\mathbf{D}] - [\mathbf{D}_{incl}]) [\mathbf{A}] )   \\
\{\mathbf{ r}\}^{\prime} & = &\{\mathbf{r}\} + [\mathbf{B}_0] 
([\mathbf{D}] - [\mathbf{D}_{incl}]) \{\mathbf{b}\} 
\end{eqnarray}
$[\mathbf{D}]$,$[\mathbf{D}_{incl}]$ are big matrices of size ($6\cdot n_g \times 6\cdot n_g$), where $n_g$ is the number of grid points, that have sub-matrices $\mathbf{D}$, $\mathbf{D}_{incl}$ on the diagonal.
Matrices $[\mathbf{A}]$ and vector $\{\mathbf{b}\}$ are computed using matrices $[\hat{\mathbf{A}}]$,$[\bar{\mathbf{B}}_0]$ and $[\hat{\mathbf{B}}_0]$ (for more details see \cite{BeerMarussig}).
It should be noted that matrices $[\mathbf{A}]$,$[\mathbf{B}_0]$ and vector $\{\mathbf{b}\}$ only depend on geometry, the elastic properties of the infinite domain and the boundary conditions.
Therefore these matrices can be pre-computed. The actual nonlinear simulation just involves matrix multiplications.
Equation \ref{onestep} means that either a Newton-Raphson or a modified Newton-Raphson method can be applied. The left hand side already includes the effect of $\mathbf{D}_{incl}$ or $\mathbf{D}_{ep,incl}$. If the inclusion is elastic then a solution can be obtained without iteration.

\

\remark{\textbf{It should be pointed out that since $\mathbf{D}_{incl}$ is defined at every grid point the method lends itself to the simulation of heterogeneous ground conditions. In fact each grid point can be assigned different material properties}}.

\

Table \ref{tab:Size} shows the sizes of the pre-computed matrices where $n_d$ is the degrees of freedom.
\begin{mytable}
  {H}               %% table position  
  {Size of matrices.}  %% caption
  {tab:Size}  %% label
  {cccccc}         %% column layout
  \mytableheader{ Matrix: & $[\mathbf{L}] $ & $[\hat{\mathbf{A}}] $ & $[\mathbf{ B}_0]$ & $[\bar{\mathbf{B}}_0]$ & $[\hat{\mathbf{ B}}]$}  %% header
  %% table content
matrix size: &  $n_d\times n_d$ & $3\cdot n_g \times n_d$ & $3 \cdot n_d \times 6 \cdot n_g$ & $3\cdot n_g \times 6 \cdot n_g$ & $ 6 \cdot n_g \times 3 \cdot n_g$ \\
\end{mytable}%

\section{The immersed boundary element method}
\label{sec:IBEM}
Here we present a method that allows cell meshes to be automatically generated. The basic idea comes from the FEM community and procedures that allow to mould a regular FEM mesh to specified surfaces. 
Immersed FEM technology is not directly applicable to the BEM since the task is somewhat different: The aim of the immersed FEM is to very accurately follow the defined surfaces. This importance is diminished for the main purpose that we use the cells, namely the numerical evaluation of volume integrals.  Furthermore, connectivity - vitally important for FEM work - can be abandoned. On the other hand one aspect, not of importance to FEM work, is important here: Since volume integration is very compute intensive it is paramount that the cell mesh is well designed, i.e. we should be aiming at a small number of cells, avoiding tetrahedral cells and an intelligent grading toward the boundary.
It should be mentioned that cells have a multiple purpose: The interpolation of displacements and the computation of strains, the interpolation of the initial stresses and the volume integration. It is also important to realise that the fundamental solution $\fund{E}$, which appears in the integral to be evaluated, decays with $o(r^{-2})$, where $r$ is the distance between $\sourcept_n$ and $\fieldpt$. This means that a fine mesh is required near the boundary but this can be coarsened further away.

Cell mesh generation starts with a definition of $\domain_{0}$, which may be a region of the domain that has different material properties than the infinite domain or a region where elasto-plastic behaviour is expected.

\

\remark{The beauty of the non-linear IGABEM is, that instead of always having to consider the whole domain as in the FEM, we can concentrate on a portion of the domain where we think the non-linear effects are significant (i.e. the near field)}

\

Here we present two approaches: Scanning on a rectangular grid and growing cells from the boundary. The reason for this is that, depending on the shape of the boundary surface, the scanning option may not result in an optimal cell mesh.

\subsection{Option 1: Scanning}
The first step is to generate a volume grid consisting of lines in three orthogonal directions. The grid can be quite coarse but has to be fine enough so that enough lines intersect the boundary.  We compute intersection points using an algorithm that determines the intersection between a line and a surface published in \cite{Beer2022S}. Next, we define cells within the grid and detect cells that have intersection points. 
For these cells the portion that intrudes into void space is trimmed. The steps are summarised in algorithm \ref{alg:Imersed}.
 \begin{myalgorithm}{Algorithm for scanning}{alg:Imersed}
	\REQUIRE  Definition of $\domain_0$
	\STATE  Generate a volume grid and define cells within the grid.
	\STATE  Compute intersections of the grid lines with boundary patches.
	\STATE  Assign intersection points to cells
	\FOR {all intersected cells}
        \STATE  Trim cells
	\ENDFOR	
	\STATE  Subdivide cells as necessary to achieve required accuracy
	\STATE  Deactivate cells that are inside the excavation
	\RETURN Generated cell mesh
\end{myalgorithm}

\subsubsection{Trimming of cells}
\begin{figure}[H]
\begin{center}
\includegraphics[scale=0.15]{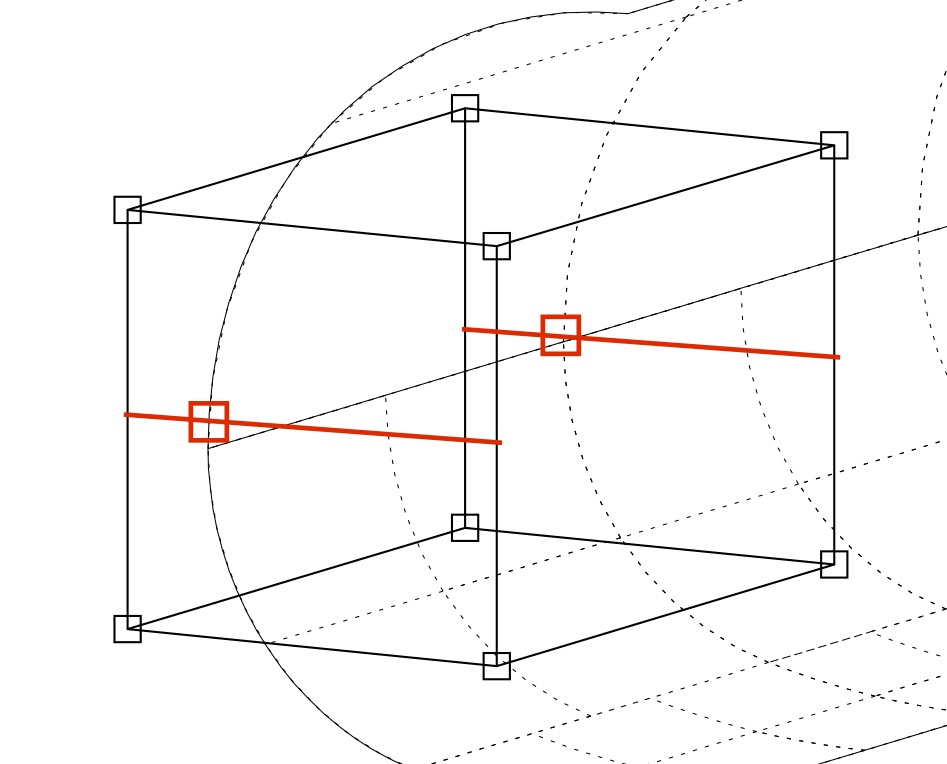}
\includegraphics[scale=0.15]{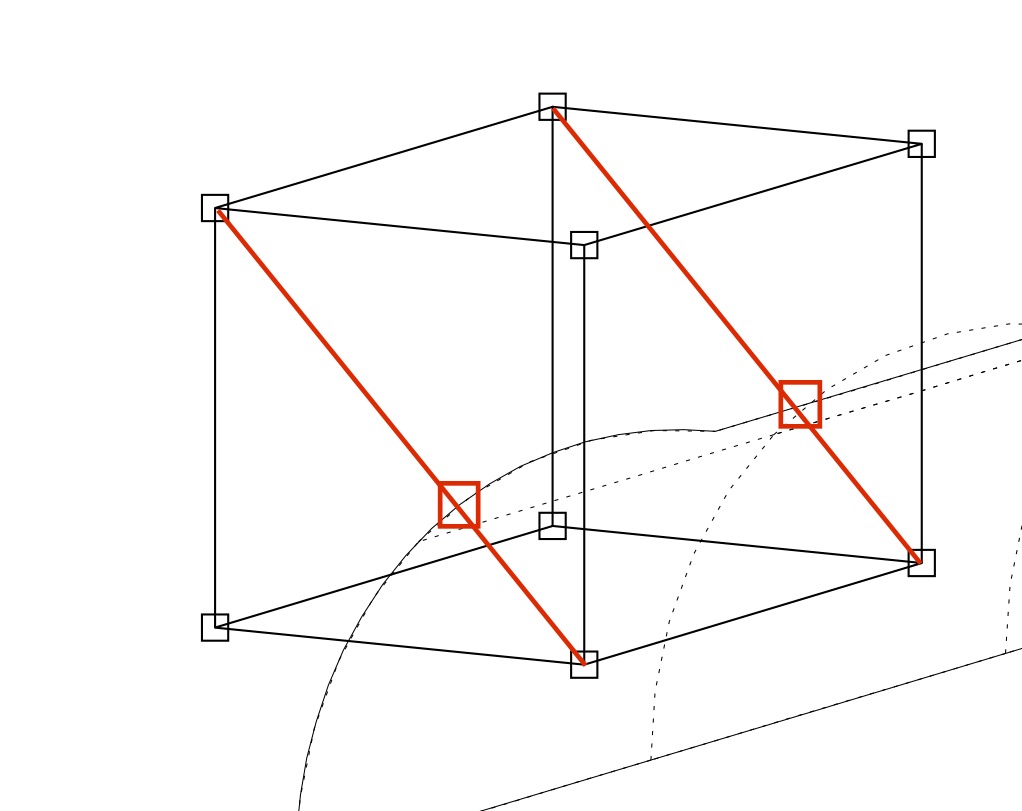}
\caption{Determination of additional points for trimming. Left: Cut through cell, Right: Cut is at a corner of cell}
\label{Trim1}
\end{center}
\end{figure}

For the trimming we must first establish which part of the cell intrudes into void space. This is done by comparing the gridline direction with the outward normal to the boundary. We distinguish between the case where a cut goes through a cell and one where only a portion is cut off. In the first case we move the corner nodes to the already computed intersection points. Next, we define lines at mid-side between the corner nodes and compute intersection  points (see \myfigref{Trim1} on the left). We then move the mid-side nodes to the computed intersection points. In the second case, we define diagonal lines as shown in  \myfigref{Trim1} on the right and compute intersection points. Using these points and the previously computed intersection points, we subdivide the cell into two sub-cells. Next we proceed as in case 1 to determine the mid-side node locations.

\newpage

\begin{figure}[H]
\begin{center}
\includegraphics[scale=0.20]{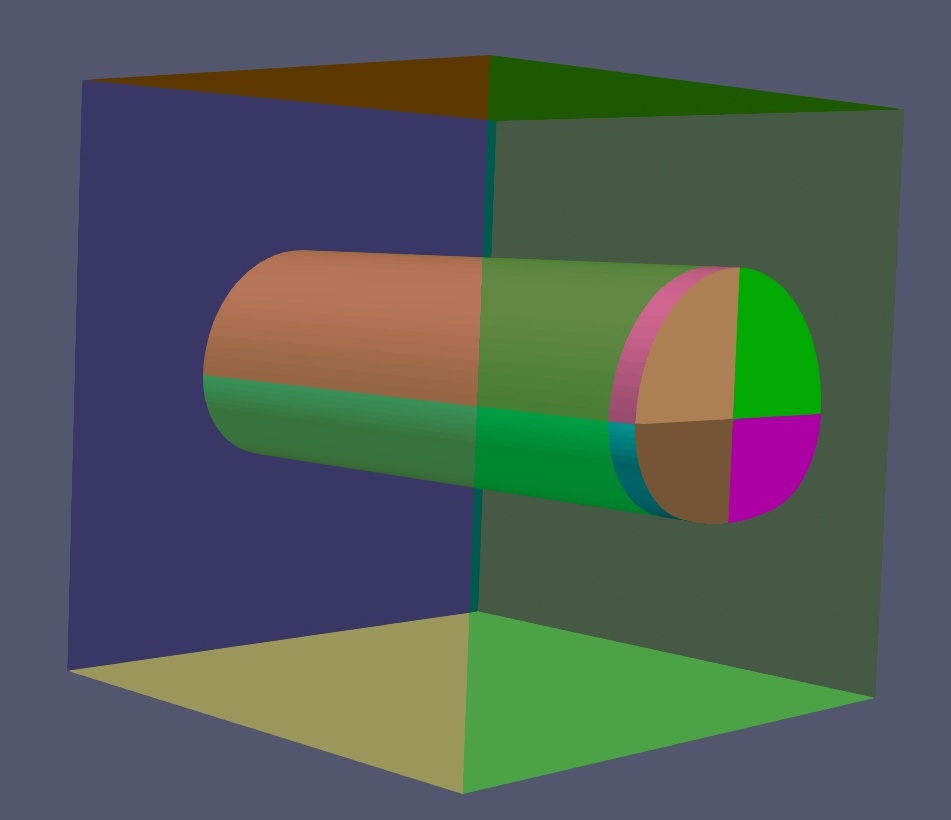}
\caption{Define  $\domain_0$ }
\label{Step1}
\end{center}
\end{figure}

\begin{figure}[H]
\begin{center}
\includegraphics[scale=0.20]{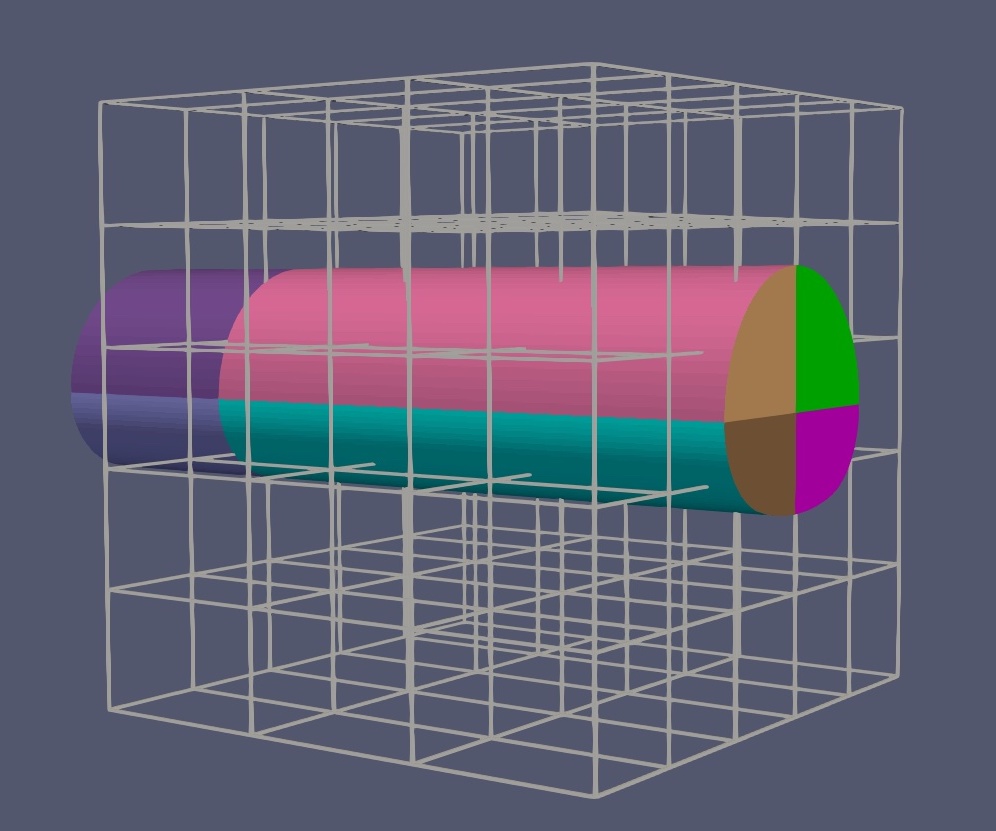}
\includegraphics[scale=0.20]{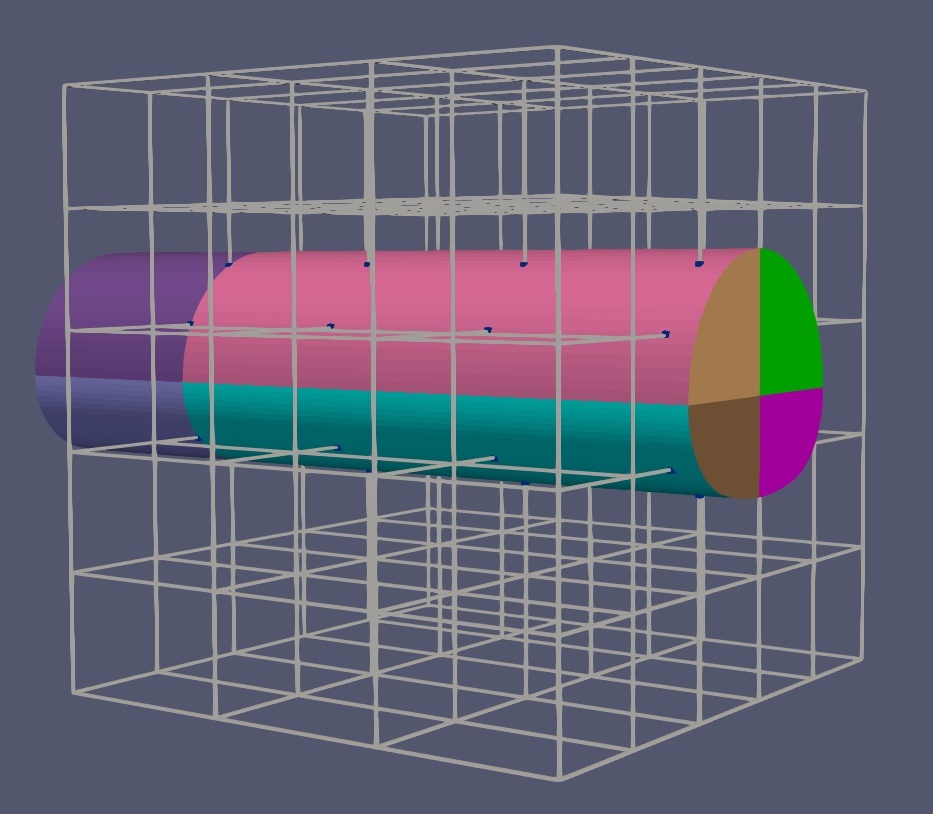}
\caption{Establish volume grid and intersect.}
\label{Step23}
\end{center}
\end{figure}

\begin{figure}[H]
\begin{center}
\includegraphics[scale=0.25]{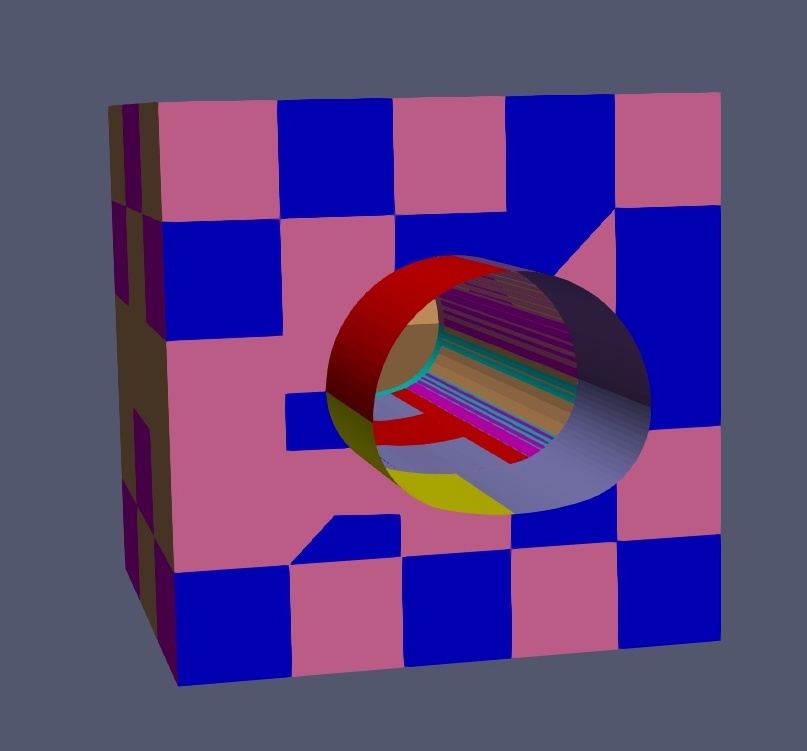}
\caption{Generated cell mesh.}
\label{Mesh}
\end{center}
\end{figure}

We explain this procedure on the example of the NATM tunnel in Figures \ref{Step1}, \ref{Step23}, \ref{Step45} and arrive at the cell mesh in \ref{Mesh}.
\begin{figure}[H]
\begin{center}
\includegraphics[scale=0.20]{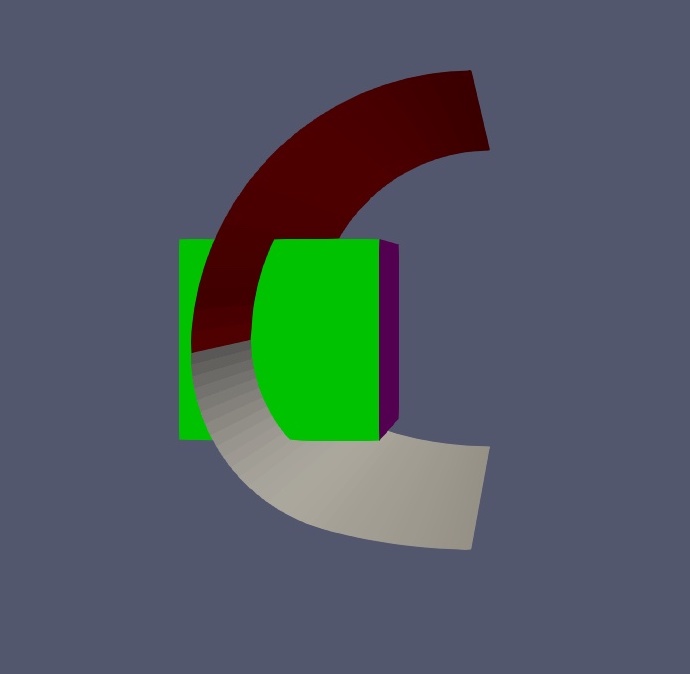}
\includegraphics[scale=0.184]{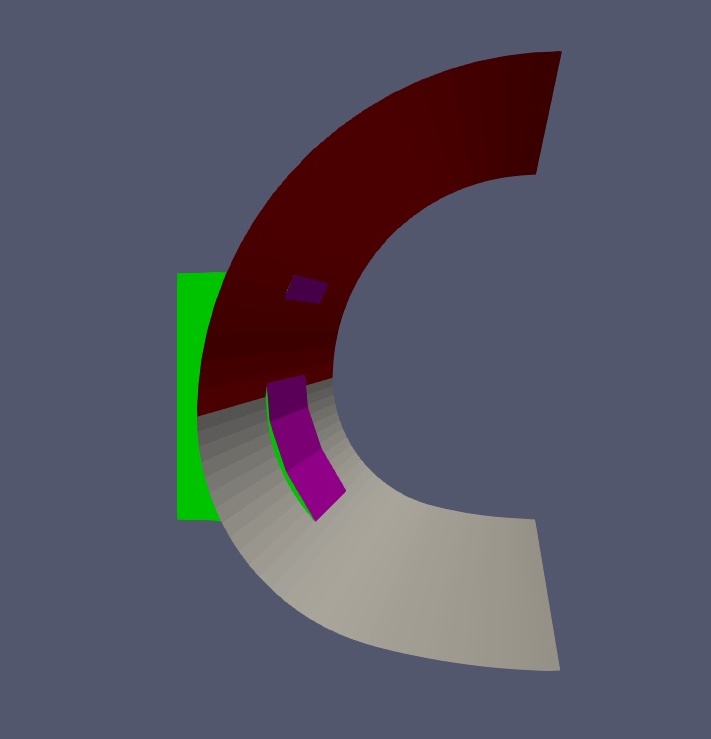}
\includegraphics[scale=0.155]{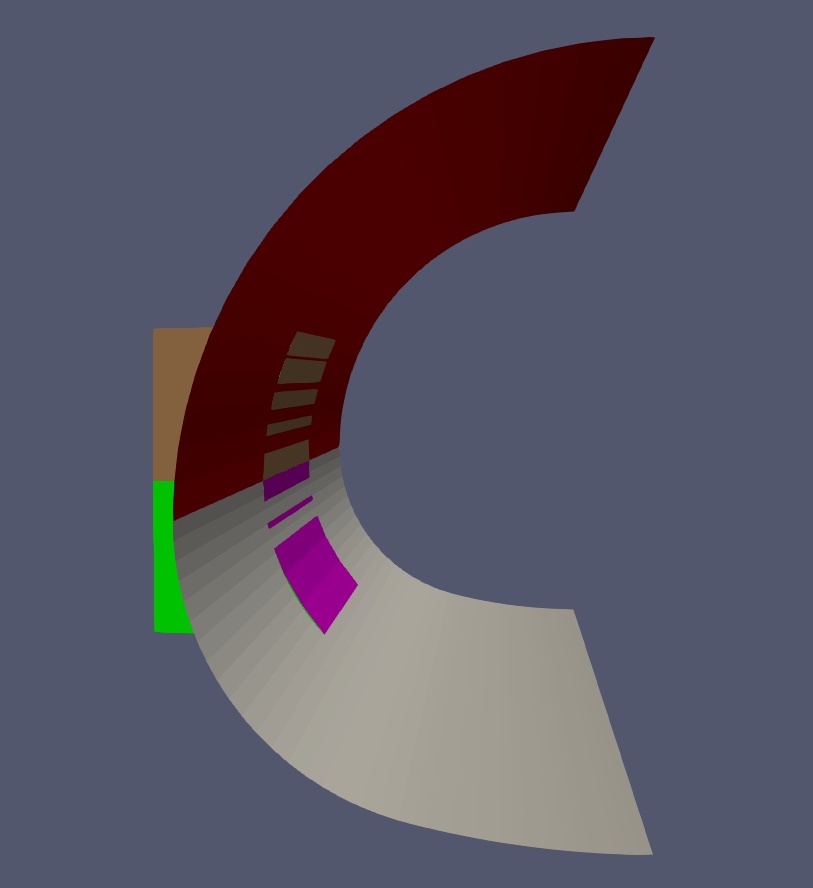}
\caption{Detect intersecting cells, trim and refine. Case 1: Cut goes through cell}
\label{Step45}
\end{center}
\end{figure}
\begin{figure}[H]
\begin{center}
\includegraphics[scale=0.20]{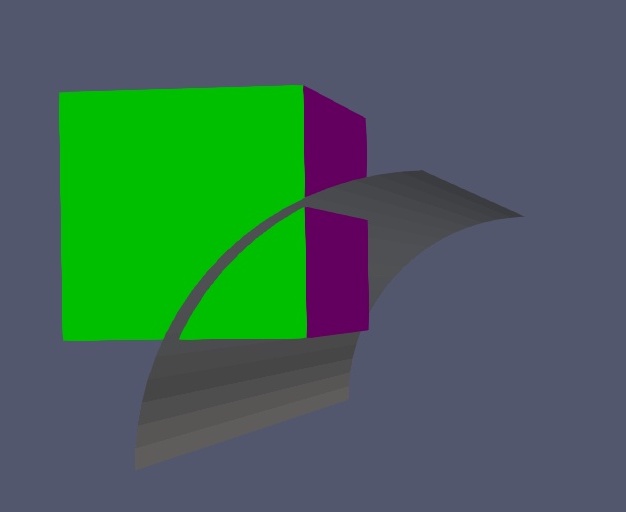}
\includegraphics[scale=0.2]{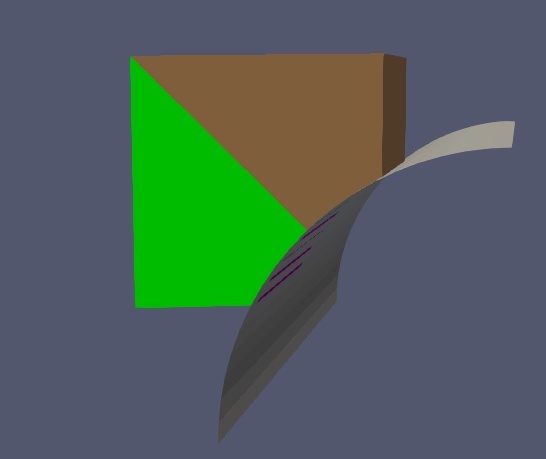}
\caption{Detect intersecting cells and trim. Case 2: Cut is at corner }
\label{Step45}
\end{center}
\end{figure}

\subsection{Option 2: Growing mesh from the boundary}
\begin{figure}[H]
\begin{center}
\includegraphics[scale=0.20]{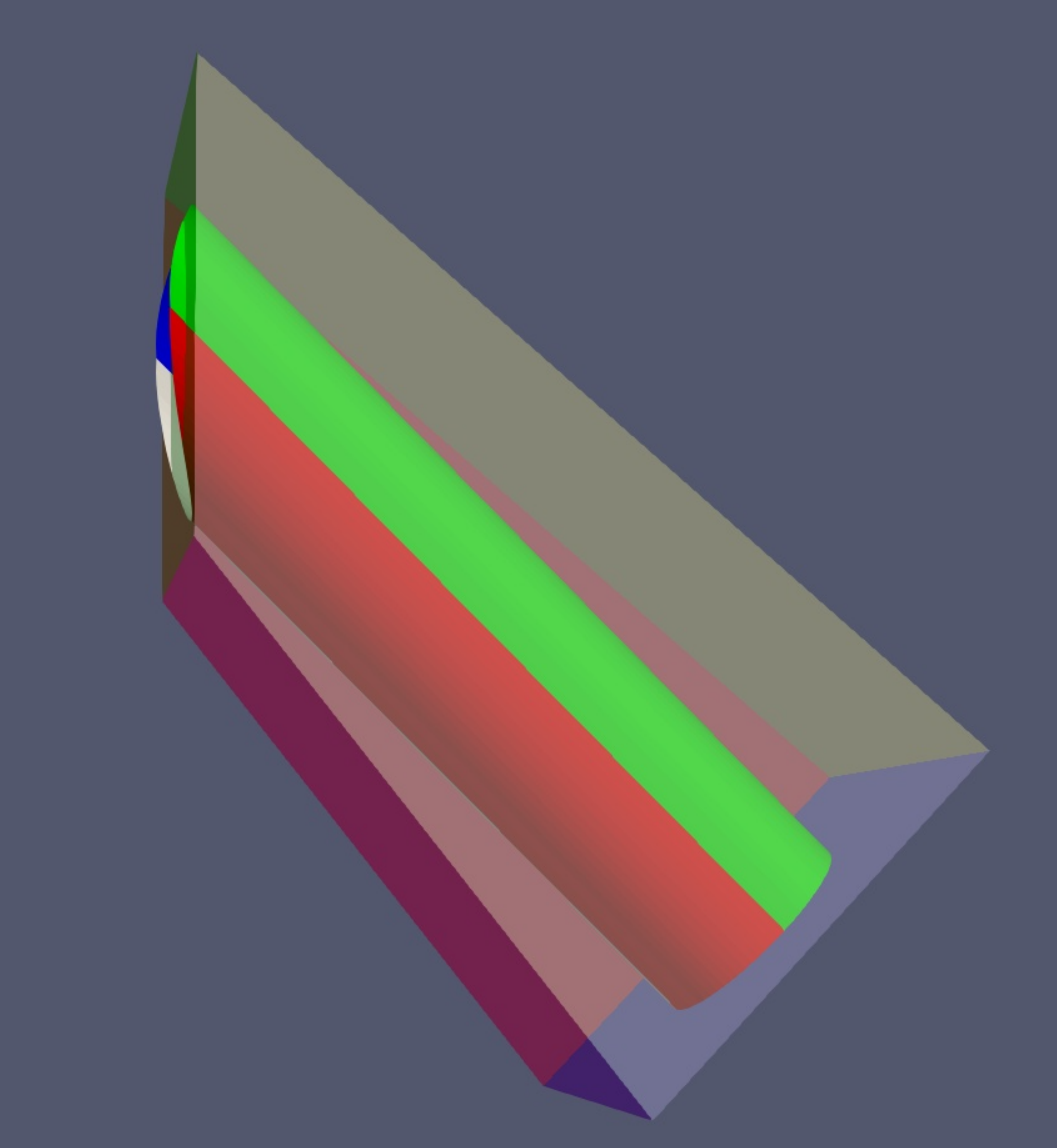}
\caption{Growing cells, step1: Definition of super-cell }
\label{scewbox}
\end{center}
\end{figure}

In this option we first define a super-cell that surrounds the boundary patches that we want the mesh to grow from (see \myfigref{scewbox}).

\begin{figure}[H]
\begin{center}
\includegraphics[scale=0.60]{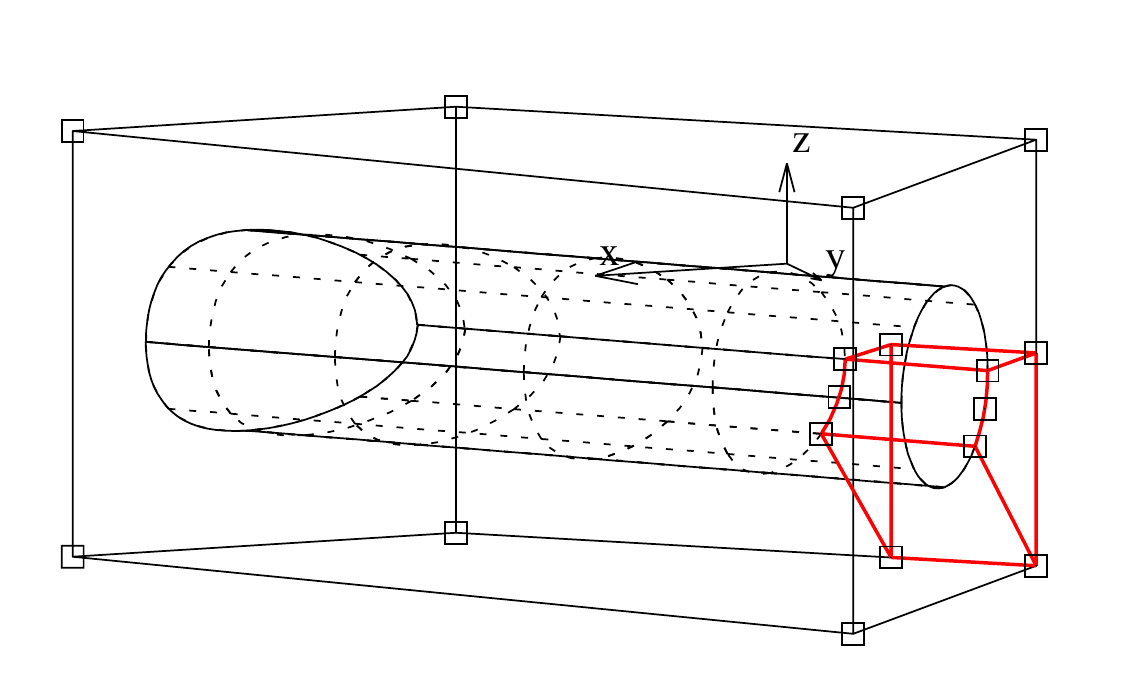}
\includegraphics[scale=0.15]{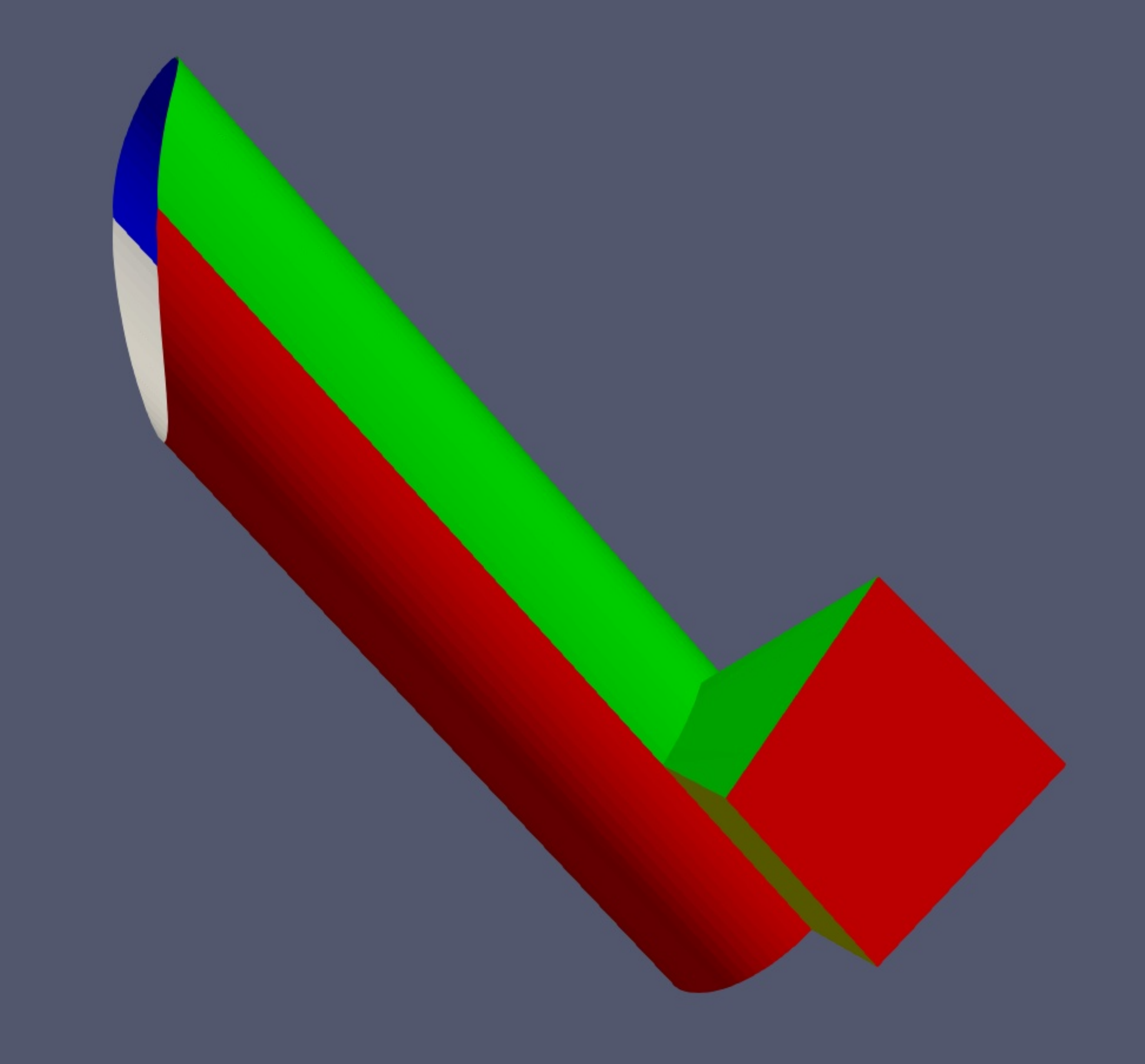}
\caption{Growing cells, step 2: Get subdivision points on patch surface and grow cell. }
\label{cellone}
\end{center}
\end{figure}

\begin{figure}[H]
\begin{center}
\includegraphics[scale=0.15]{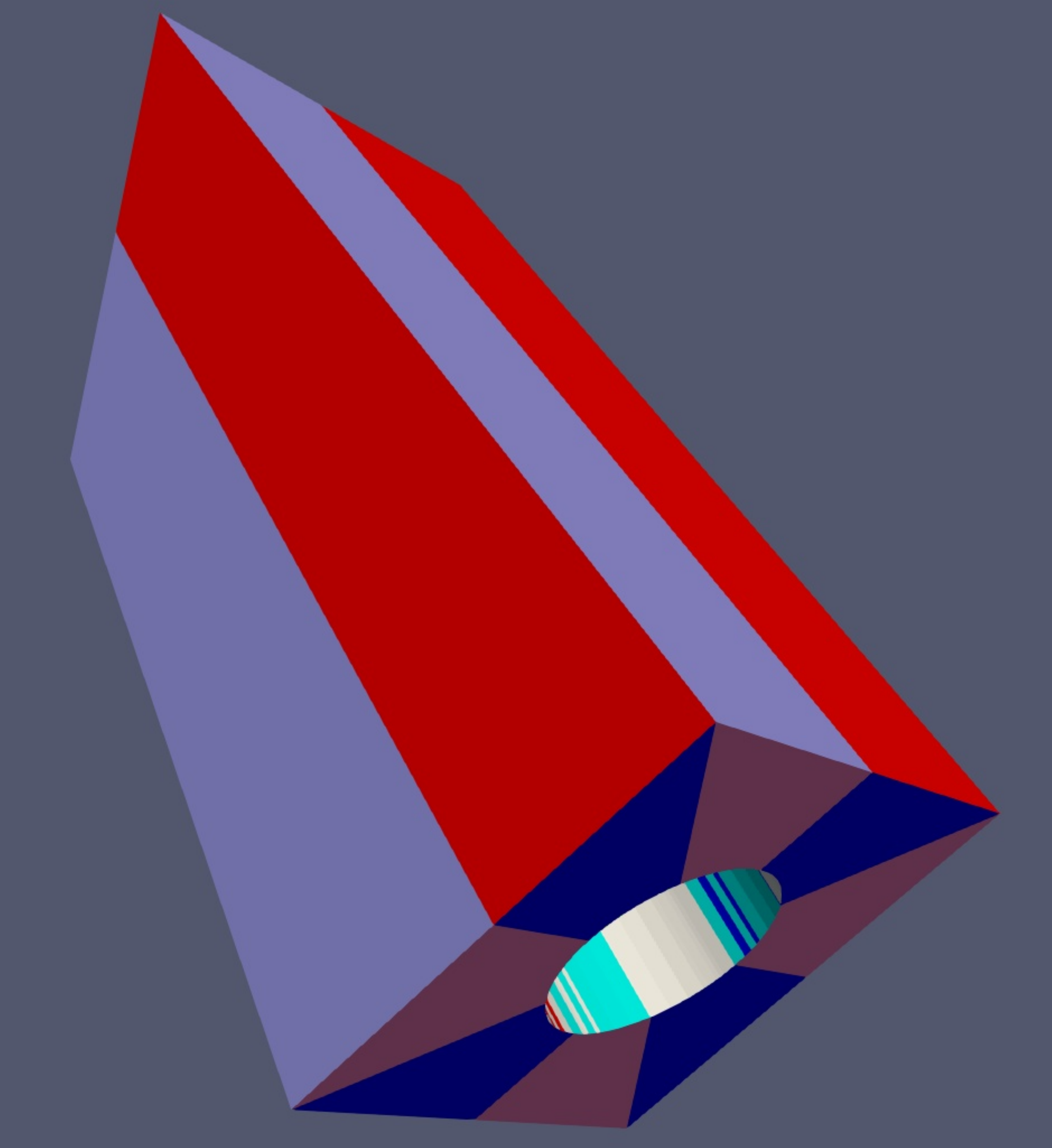}
\caption{Final grown cell mesh}
\label{cellall}
\end{center}
\end{figure}

Next we subdivide the surface of the patches depending on the number of subdivisions required, calculate the coordinates of subdivision points using \myeqref{mapatch} and connect the points on the patch to the super-cell nodes (see \myfigref{cellone}).
The generated cell mesh is shown in (see \myfigref{cellall})

\section{Example}
We revisit the example in section \ref{ExBorehole} but this time we perform a nonlinear simulation.
The additional input data are as follows:
\begin{itemize}
  \item Plastic properties: Mohr-Coulomb, $\phi$= 30° , c=16.6 MPa, Pore pressure= 72 MPa
\end{itemize}
The main aim of the simulation was to investigate the plastic zone near the deviation.

The generation of a cell mesh consists of the following steps:
\begin{enumerate}
  \item Generate cell mesh for the main borehole only, using the scanning option
  \item Identify the cell at the junction between the main and the deviated borehole
  \item Change this cell to a super-cell by skewing it in the direction of the deviated borehole
  \item Grow cells from the deviated borehole
\end{enumerate}
As we change the cell to a supercell and skew it, the adjoining cells in vertical direction have to be skewed also to prevent overlap. Since there is no continuity requirement the cells in the other direction are not affected.
\begin{figure}[H]
\begin{center}
\includegraphics[scale=0.25]{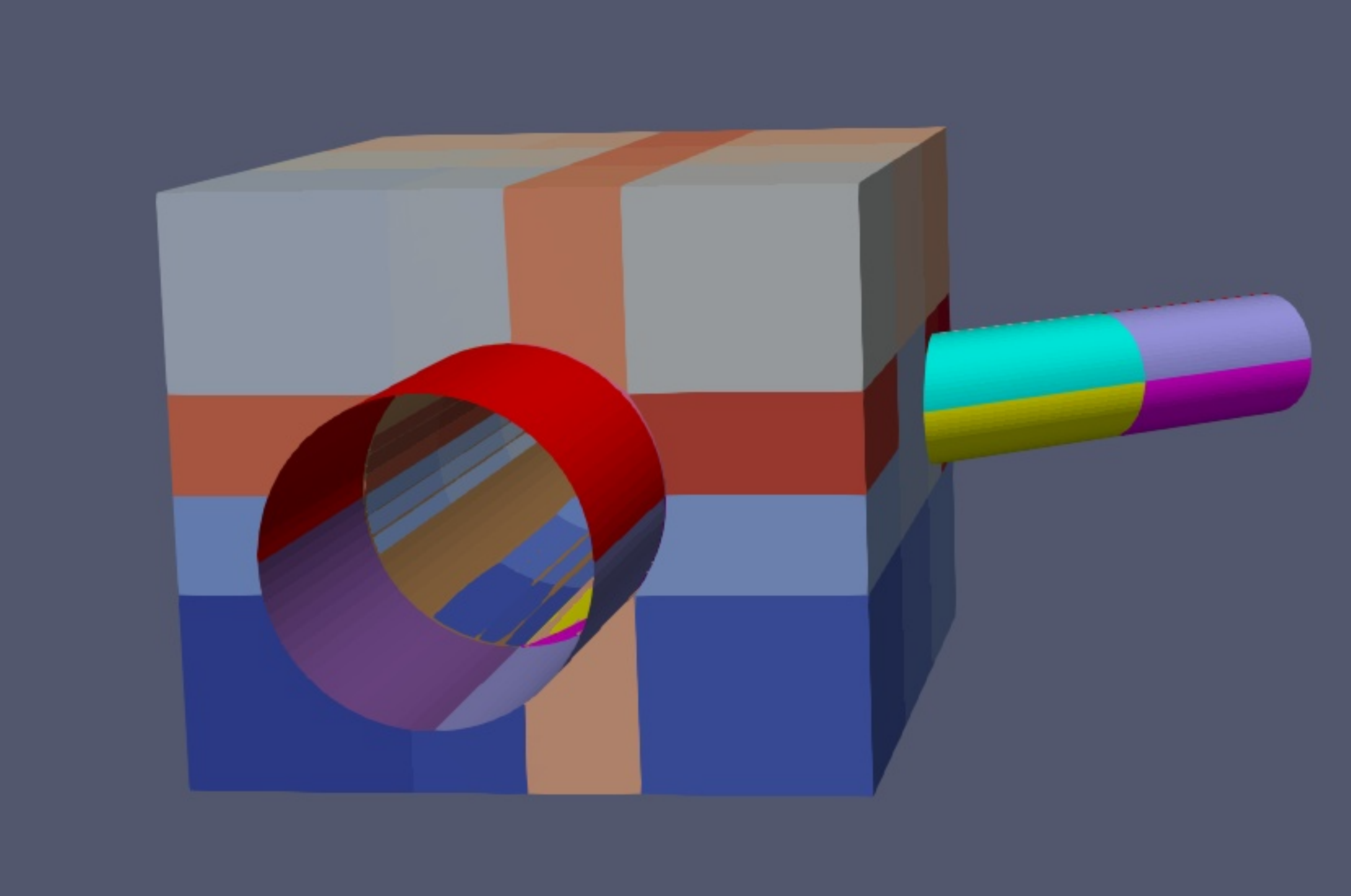}
\caption{Step 1: Generation of cell mesh around main borehole.}
\label{mesh1}
\end{center}
\end{figure}
\begin{figure}[H]
\begin{center}
\includegraphics[scale=0.15]{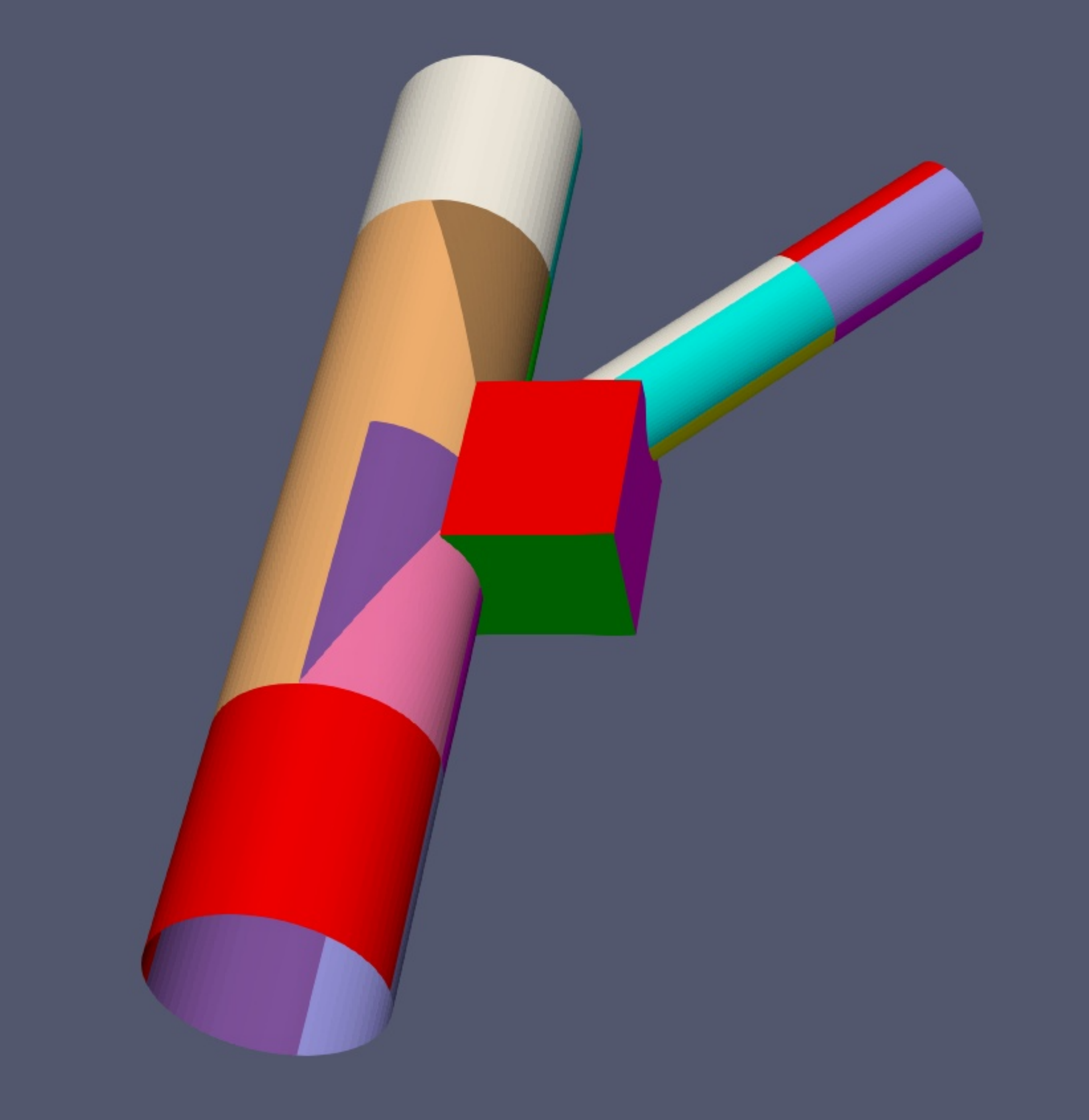}
\includegraphics[scale=0.15]{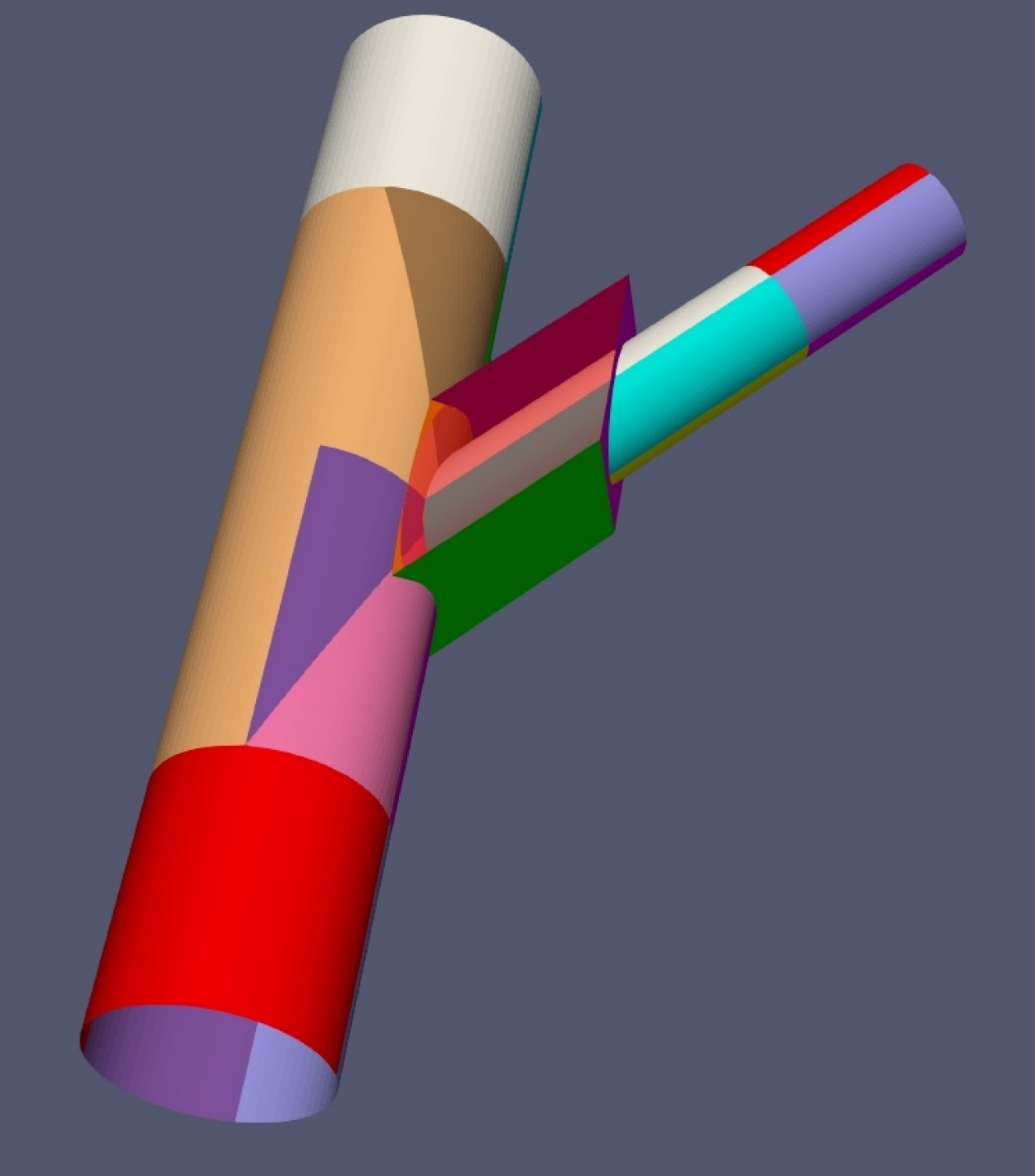}
\includegraphics[scale=0.15]{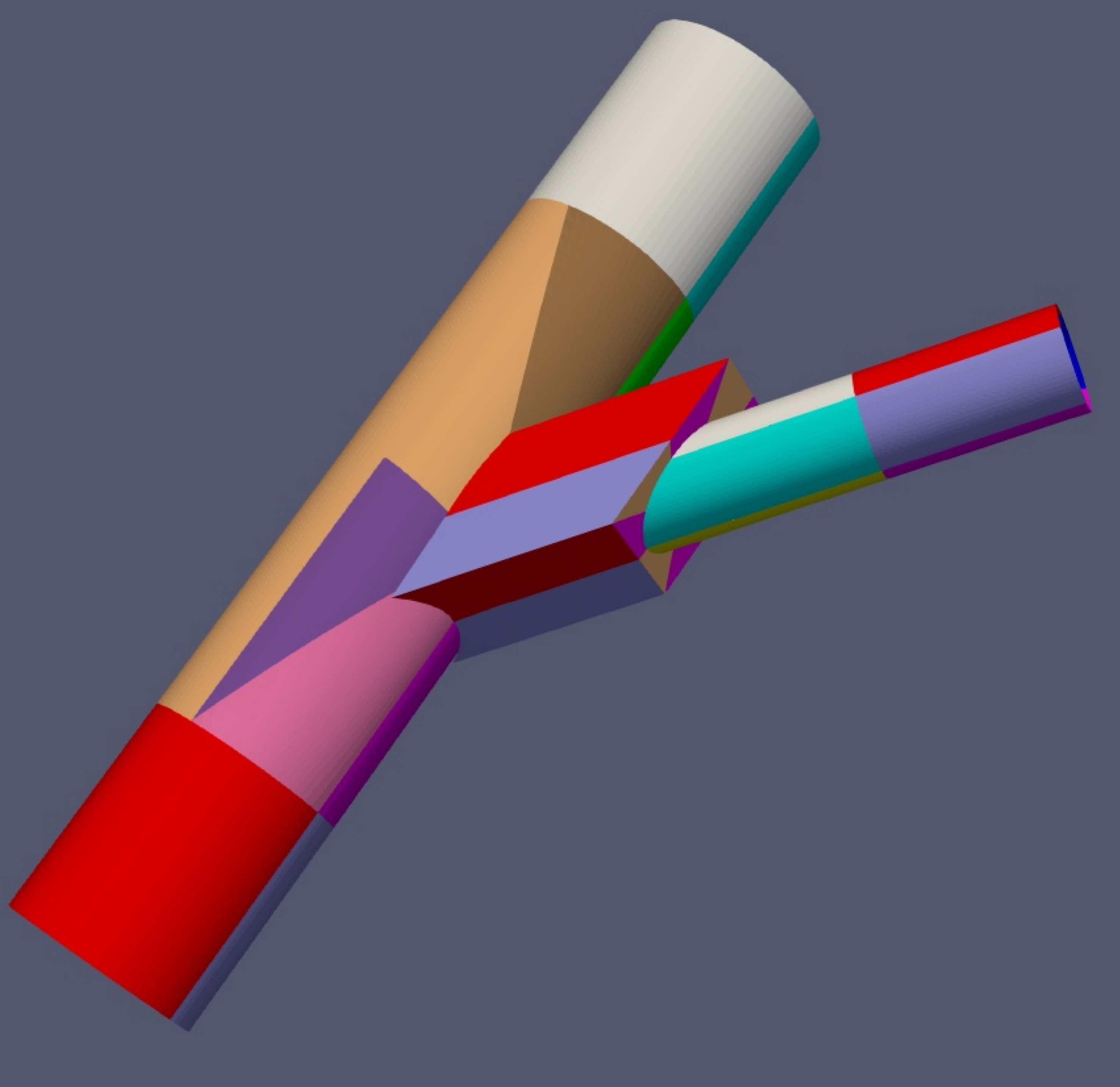}
\caption{Steps 2 to 4: Identify cell at intersection, generate super-cell and grow cells.}
\label{mesh23}
\end{center}
\end{figure}
\begin{figure}[H]
\begin{center}
\includegraphics[scale=0.2]{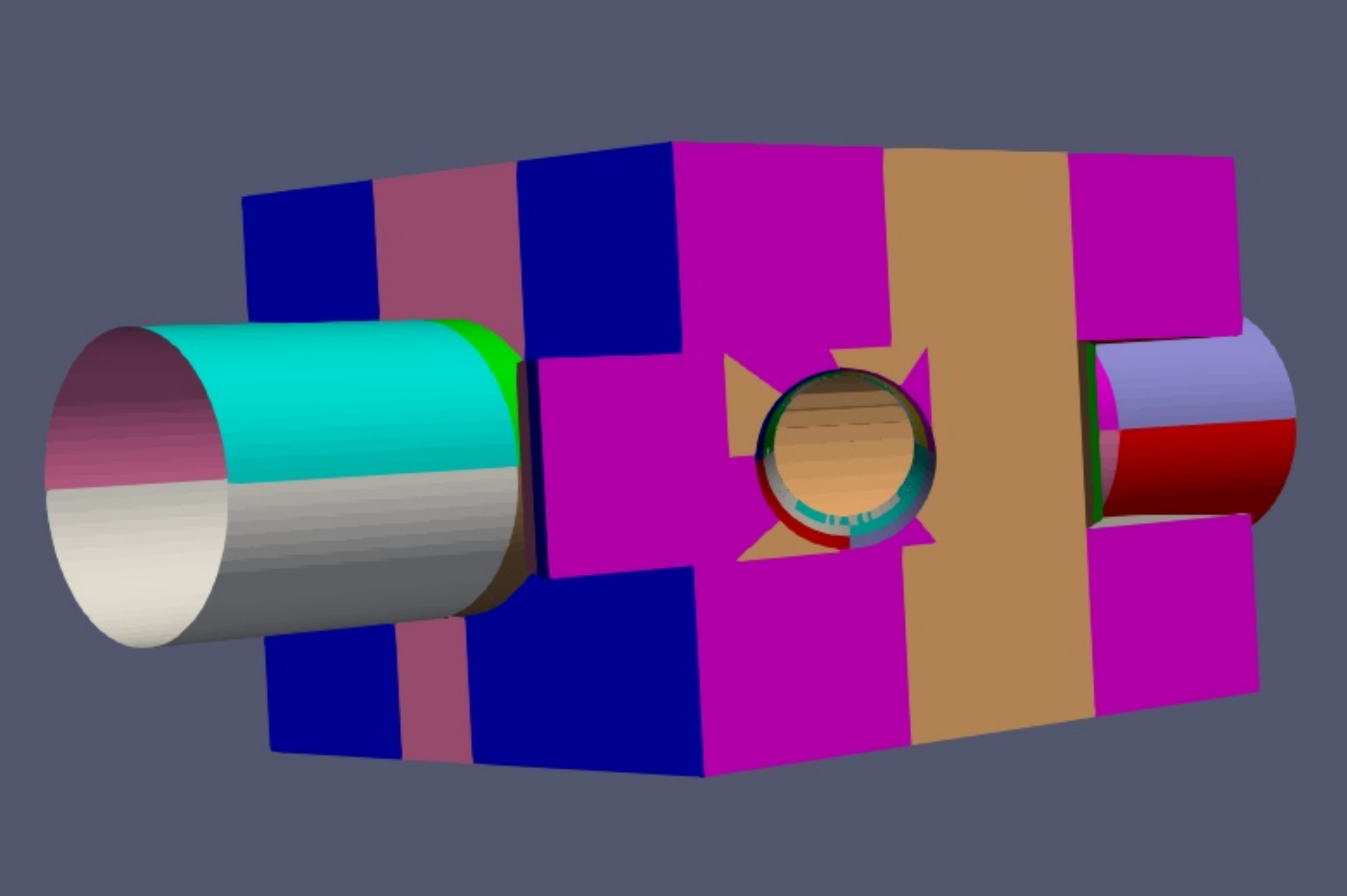}
\caption{Final cell mesh.}
\label{meshfinal}
\end{center}
\end{figure}

The generation of the cell mesh is explained in Figures \ref{mesh1} and \ref{mesh23}. The final mesh is shown in \myfigref{meshfinal} and has 76 20-node cells and 551 grid points.

\

\remark{In this first implementation some user interaction is required. However, the aim in the future is that cell mesh generation is made completely automatic}.

\

\begin{figure}[H]
\begin{center}
\includegraphics[scale=0.35]{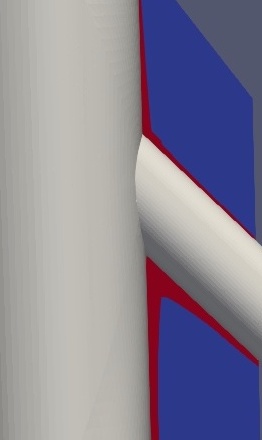}
\includegraphics[scale=0.25]{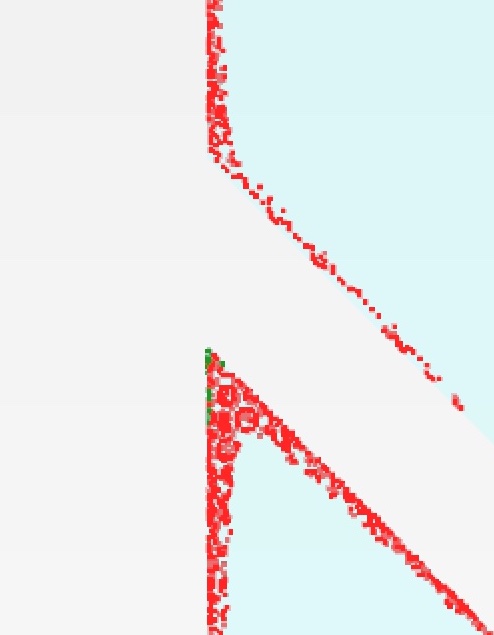}
\caption{Extent of plastic zone (marked in red) from the IGABEM simulation and (right) marked plastic points from the FEM simulation.}
\label{Bore4}
\end{center}
\end{figure}

The results were compared with a 3-D FEM simulation with Program PLAXIS, which had about one million degrees of freedom. In contrast the BEM simulation had only 376 degrees of freedom. 
In \myfigref{Bore4} we show one result of the simulation namely the extent of the plastic zone near the deviation and a comparison with plastic points of the FEM simulation. Good agreement can be found.

Table \ref{tab:Run} shows the run times for precomputing the matrices for the nonlinear solution on a MacBook pro laptop. It also shows the size of the matrices. The software was written in MATLAB using the NURBS toolkit \cite{Spink}.
\begin{mytable}
  {H}               %% table position  
  {Run times for pre-computing matrices.}  %% caption
  {tab:Run}  %% label
  {cccccc}         %% column layout
  \mytableheader{ Matrix/Vector: & $[\mathbf{L}],\{ \mathbf{r} \}$ & $[\hat{\mathbf{A}}], \{\bar{\mathbf{c}} \} $ & $[\mathbf{ B}_0]$ & $[\bar{\mathbf{B}}_0]$ & $[\hat{\mathbf{ B}}]$}  %% header
  %% table content
matrix size: & $376\times 376$ & $1653\times 376$ & $ 376 \times 3306$ & $1653 \times 3306$ & $ 3306 \times 1653$ \\
run time (sec): & 66& 128 & 91 & 355 & 1  \\
\end{mytable}%

The actual nonlinear solution took only 4 seconds. It should be noted that the evaluation of the matrices that use volume integration is well suited for parallel computation, since integration of each cell is independent of other cells. A proper implementation in C++ or Fortran with parallel processing features would drastically reduce run times for the pre-computation.
Again it has to be stressed that once the matrices are pre-computed, simulations with different inclusion properties can be carried out very fast.
A re-computation of only $\{\bar{\mathbf{c}} \}$ and $\{ \mathbf{r} \}$ is necessary if the loading has changed. A re-computation of the other matrices is only necessary if the geometry or the elastic properties of the infinite domain, with which the fundamental solutions are computed, have changed.
\begin{figure}[H]
\begin{center}
\includegraphics[scale=0.2]{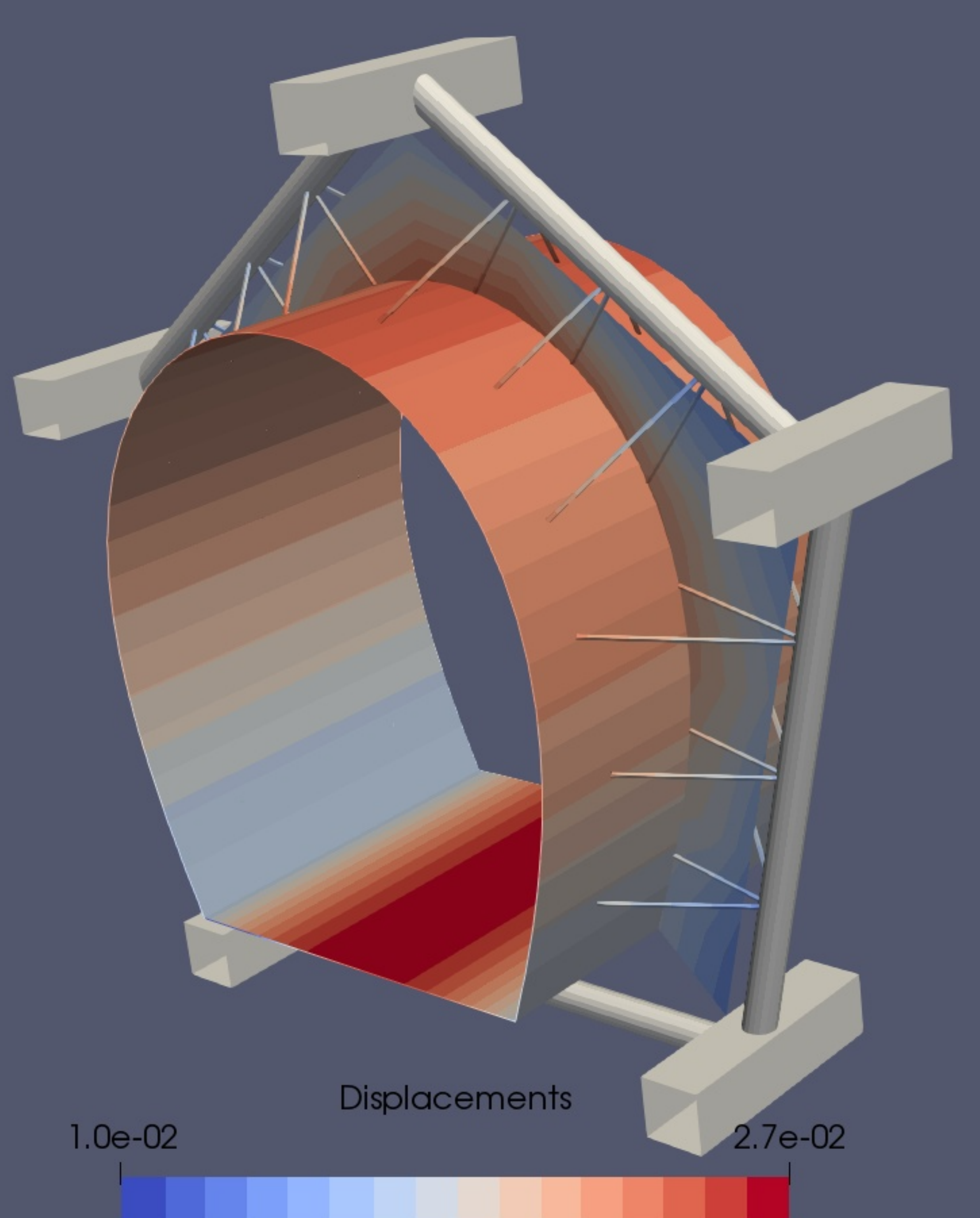}
\caption{BEM simulation of the effect of rock bolts on the displacement of an excavation. Taken from  \cite{Beer2018a}  }
\label{RibinRock}
\end{center}
\end{figure}

\section{Summary and Outlook}
In this paper a novel approach to the simulation of underground excavations was presented. The BEM is ideally suited for such problems that involve infinite or semi-infinite domains. However, there has been limited acceptance of the method and the FEM still dominates the geomechanics simulation software available.  Progressive software companies use BEM software, but only for elastic, piecewise homogeneous domains.
The limited acceptance is due to the fact that, compared to the FEM, very few researchers are working on the BEM. In addition, features important for geomechanics such as ground support, heterogeneous ground conditions and elasto-plasticity were not properly addressed.

The aim of our research over the last decades was to remedy the problem and to increase the applicability of the method in geomechanics. 
The work started with the introduction of isogeometric methods that greatly simplified the definition of excavation geometries, in essence avoiding mesh generation. Using this approach, geometry information can taken directly from Computer Aided Design (CAD) programs, but CAD data have to be post-processed before being applied to simulation. The definition of the geometry by users, without CAD, is also greatly simplified. Geometries can be defined very accurately (and sometimes even exact) with few parameters.
The use of NURBS also makes the refinement of the solution very user friendly, as no re-meshing is involved. In this paper we hope to have demonstrated the beauty of the IGABEM.

The most time consuming research effort, however, was concerned with the treatment of heterogeneous ground conditions and plasticity.  A major stumbling block was the requirement of the generation of a volume mesh. To leave this up to the user was unacceptable. In this paper we have shown a possibility of making cell generation automatic. On a practical example we have shown that a combination of scanning and cell growing leads an acceptable cell mesh, in terms of computation time and accuracy. Results of the IGABEM simulation are comparable with the ones form a FEM simulation, but obtained with considerable less effort.
Further research is needed to develop software that automatically generates cell meshes, that are robust for all geometries, without any user intervention. Maybe this needs to also include using 3-D NURBS, that have been presented previously in \cite{Beer2016}.  Indeed, the ultimate aim should be a situation where the existence of a volume mesh is invisible to the user.

Finally, a reminder that in previous papers \cite{Beer2018a,BEER2020113409} we have successfully tackled the problem of modelling ground support (rock bolts, concrete lining). See for example \myfigref{RibinRock}). This means that the IGABEM can perform all the simulation tasks (at least for underground excavations) that the FEM can, albeit with considerable less effort.

A lot needs to be done to pack these developments into a suitable software product, starting from the migration into a more efficient language than MATLAB and parallelisation. 
Also, intensive testing is required to ensure that the software is stable and leads to acceptable results for all possible cases.
Unfortunately, due to lack of funding, we are not able to continue this work and we hope that this paper serves as an inspiration to readers. Maybe some progressive software company would like to turn this into a next generation simulation software for geomechanics. Our assistance is assured.

%\bibliographystyle{apalike}
%\bibitem[Author(year)]{label}
\bibliographystyle{apalike}\biboptions{authoryear}
\bibliography{bookbib}

\end{document}